\documentclass[a4paper,11pt,reqno,oneside,abstracton]{scrartcl}
\usepackage{scalerel,stackengine,relsize}
\stackMath
\usepackage{amssymb}
\usepackage{bbm}
\usepackage{ifthen}
\usepackage{tikz}
\usepackage{tikz-3dplot}
\usepackage{bm}
\usepackage[hyphens]{url} 
\usepackage{amsmath,amssymb,amsfonts,amsthm,amsopn,dsfont,mathrsfs}
\usepackage{esint}
\usepackage{graphicx}
\usepackage{enumitem}
\usepackage{bigints}
\usepackage{fancyhdr}
\usepackage[T1]{fontenc}
\usepackage{mathtools}
\usepackage[utf8]{inputenc}
\usepackage{verbatim} 
\usepackage[backend=bibtex,style=numeric,bibencoding=utf8,language=auto,autolang=other,maxnames=10]{biblatex}
\usepackage[a4paper, total={6in, 9in}]{geometry}

\theoremstyle{plain}
\newtheorem{thm}{Theorem}[section]
\newtheorem*{thm*}{Theorem}
\newtheorem*{conj*}{Conjecture}
\newtheorem{prop}{Proposition}[section]
\newtheorem{lem}{Lemma}[section]
\newtheorem{cor}{Corollary}[section]
\theoremstyle{definition}

\newtheorem{dfn}{Definition}[section]
\theoremstyle{remark}
\newtheorem{rem}{Remark}[section]
\numberwithin{equation}{section}

\renewcommand{\r}{\mathbb{R}}
\newcommand{\R}{\mathcal{R}}
\newcommand{\s}{\mathbb{S}}

\newcommand{\z}{\mathbb{Z}}
\newcommand{\eps}{\varepsilon}
\newcommand{\D}{\mathcal{D}}
\newcommand{\m}{\mathbb{M}}

\newcommand{\n}{\mathbb{N}}
\renewcommand{\le}{\leqslant}\renewcommand{\leq}{\leqslant}
\renewcommand{\ge}{\geqslant}\renewcommand{\geq}{\geqslant}

\renewcommand{\H}{\mathscr{H}}
\renewcommand{\L}{\mathcal{L}}
\DeclareRobustCommand{\rchi}{{\mathpalette\irchi\relax}}
\newcommand{\irchi}[2]{\raisebox{\depth}{$#1\chi$}}
\DeclareMathOperator{\vol}{vol}

\DeclareMathOperator{\im}{Im}
\DeclareMathOperator{\dist}{dist}

\DeclareMathOperator{\spt}{spt}

\DeclareMathOperator{\dom}{Dom}

\newcommand{\res}{\mathbin{\vrule height 1.6ex depth 0pt width
		0.13ex\vrule height 0.13ex depth 0pt width 1.3ex}}
\newcommand\reallywidehat[1]{%
	\savestack{\tmpbox}{\stretchto{%
			\scaleto{%
				\scalerel*[\widthof{\ensuremath{#1}}]{\kern-.6pt\bigwedge\kern-.6pt}%
				{\rule[-\textheight/2]{1ex}{\textheight}}
			}{\textheight}%
		}{0.5ex}}%
	\stackon[1pt]{#1}{\tmpbox}%
}
\newcommand{\Addresses}{{
  \bigskip
  \footnotesize

  R.~Caniato, \textsc{Department of Mathematics, ETH, R\"amistrasse 101, CH-8093 Z\"urich, Switzerland}\par\nopagebreak
  \texttt{riccardo.caniato@math.ethz.ch}

  \medskip

  F.~Gaia, \textsc{Department of Mathematics, ETH, R\"amistrasse 101, CH-8093 Z\"urich, Switzerland}\par\nopagebreak
  \texttt{filippo.gaia@math.ethz.ch}

}}

\title{Weak and strong $L^p$-limits of vector fields with finitely many integer singularities in dimension $n$}
\author{Riccardo Caniato \and Filippo Gaia}
\date{\today}

\begin{document}
\begin{titlepage}
\maketitle
\begin{abstract}
    For every given $p\in [1,+\infty)$ and $n\in\mathbb{N}$ with $n\ge 1$, the authors identify the strong $L^p$-closure $L_{\mathbb{Z}}^p(D)$ of the class of vector fields having finitely many integer topological singularities on a domain $D$ which is either bi-Lipschitz equivalent to the open unit $n$-dimensional cube or to the boundary of the unit $(n+1)$-dimensional cube. Moreover, for any $n\in \mathbb{N}$ with $n\geq 2$ the authors prove that $L_{\mathbb{Z}}^p(D)$ is weakly sequentially closed for every $p\in (1,+\infty)$ whenever $D$ is an open domain in $\mathbb{R}^n$ which is bi-Lipschitz equivalent to the open unit cube. As a byproduct of the previous analysis, a useful characterisation of such class of objects is obtained in terms of existence of a (minimal) connection for their singular set.\\ 
    
    \vspace{1mm}
    \hspace{-3.25mm}\textbf{MSC (2020)}: 46N20 (primary); 58E15 (secondary).
\end{abstract}

\tableofcontents
\end{titlepage}
\section{Introduction}
Throughout the following paper we will always use the notations listed in subsection 1.4.
%
\subsection{Statement of the results and motivation}
Given any smooth map $u:X\to Y$ between closed, oriented and connected $(n-1)$-dimensional manifolds, the \textit{degree} of $u$ is a measure of how many times $X$ wraps around $Y$ under the action of $u$. It can be defined as follows\footnote{An alternative definition can be given in terms of the orientations of the preimages of regular points of $u$, see for instance \cite[Chapter 7]{berger-gostiaux}}:
\begin{align*}
    \text{deg}(u)=\int_{X}u^*\omega,
\end{align*}
with $\omega\in\Omega^{n}(Y)$ being a renormalized volume form on $Y$, i.e. $\omega$ is nowhere vanishing and \begin{align*}
    \int_Y\omega=1.
\end{align*}
Now let $D\subset\r^{n}$ be an open and bounded Lipschitz domain and consider a map $u\in W^{1,p}(D,Y)$ for some $p\ge 1$ being \textit{smooth up to finitely many point singularities}, which simply means that $u\in C^{\infty}(D\smallsetminus S_u,Y)$ for some finite set $S_u\subset D$. In this case we write $u\in R^{1,p}(D,Y)$. 
We define the \textit{degree of $u$ at some singular point $x\in S_u$} as
\begin{align}\label{definition of degree at some point}
    \text{deg}(u,x):=\deg(u|_{\partial D'})=\int_{\partial D'}u^*\omega\in\z,
\end{align}
where $D'\subset\subset D$ is any open, piecewise smooth domain in $D$ such that $\overline{D'}\cap S_u=\{x\}$. Notice that Definition \eqref{definition of degree at some point} is independent from the choice of the set $D'$.\\
If $\text{deg}(u,x)\neq 0$ for some $x\in S_u$, then we say that $x$ is a \textit{topological singularity} of $u$ and we refer to the subset of $S_u$ made of the topological singularities of $u$ as the \textit{topological singular set of} $u$, which we denote by $S_u^{top}$.\\
Notice that if $u\in R^{1,n-1}(D)$ then $ u^\ast\omega\in \Omega_1^{n-1}(D)$, moreover by \eqref{definition of degree at some point} we see that the $u^*\omega$ "detects" the topological singularities of $u$, in the sense that
\begin{align}
\label{Equation: degree detects topological singularities}
    \int_{\partial D'}u^*\omega=\sum_{x\in S_u^{top}\cap D'}\deg(u,x)
\end{align}
for every open, piecewise smooth domain $D'\subset\subset D$ such that $\partial D'\cap S_u=\emptyset$. From \eqref{Equation: degree detects topological singularities} one can deduce that
\begin{align*}
    *d(u^*\omega)=\sum_{x\in S_u^{top}}\text{deg}(u,x)\delta_x\text{ in }\mathcal{D}'(D).
\end{align*}

\begin{rem}
Sobolev maps smooth up to a finite set of topological singularities arise frequently as solutions of variational problems in critical or supercritical dimension. For example, this is the best regularity which is possible to guarantee for energy minimizing harmonic maps in $W^{1,2}(B^3,\s^2)$ (see \cite[Theorem II]{schoen-uhlenbeck}). Again, quite recently the second author and T. Rivi\`ere considered the following "weak $\s^1$-harmonic map equation"
\begin{align*}
    \text{div}(u\wedge \nabla u)=0 \qquad\, \mbox{ in } \mathcal{D}'(B^2)
\end{align*}
and gave a completely variational characterization of the solutions in $R^{1,p}(B^2,\s^1)$ with finite "renormalized Dirichlet energy" for $p>1$ (see \cite[Theorem I.3]{gaia-riviere} for further details).

We also remark that the presence of topological singularities is deeply linked to fundamental questions concerning the strong $W^{1,p}$-approximability through smooth maps of elements in $W^{1,p}(D,Y)$ (see \cite{bethuel-coron-demengel-helein}, \cite{hardt-riviere-approximation}).
\end{rem}

The previous discussion motivates the following general definition.
\begin{dfn}
	Let $p\in [1,\infty]$. Let $F\in\Omega_p^{n-1}(D)$. We say that $F$ has \textit{finitely many integer singularities} if there exists a finite set of points $S\subset D$ such that $F\in\Omega^{n-1}(D\smallsetminus S)$ and 
	\begin{align*}
	*dF=\sum_{x\in S}a_x\delta_{x},
	\end{align*}
	where $a_x\in\z$ for every $x\in S$.
	The class of $L^p$ integrable $(n-1)$-forms on $D$ having finitely many integer singularities will be denoted by $\Omega_{p,R}^{n-1}(D)$.
\end{dfn}
As we have seen above, $u^\ast\omega\in \Omega_{1,R}^{n-1}(D)$ for every $u\in R^{1,n-1}(D,Y)$, for any closed, oriented and connected $n-1$-manifold $Y$. Other simple examples of elements of $\Omega_{p,R}^{n-1}(D)$ can be constructed as follows.
Let $\sigma: D\to \mathbb{R}$ denote the fundamental solution of the Laplace equation, i.e.
\begin{align*}
    \sigma(x)=\begin{cases}
        -\lvert x\rvert\quad & \text{if } n=1,\\
        -\frac{1}{2\pi}\log\lvert x\rvert\quad & \text{if } n=2\\
        \frac{1}{n(n-2)\alpha(n)}\frac{1}{\lvert x\rvert^{n-2}}\quad & \text{if } n\geq 3,
    \end{cases}
\end{align*}
where $\alpha(n)$ denotes the volume of the unit ball in $\mathbb{R}^n$. 
Then $\ast d\sigma\in \Omega_{p,R}^{n-1}(D)$ for any $p\in [1,\frac{n}{n-1})$. In fact $\ast d(\ast d\sigma)=\Delta \sigma=\delta_0$.\\
Clearly any finite linear combination with integer coefficients of translations of $\ast d\sigma$ also belongs to $\Omega_{p,R}^{n-1}(D)$. In fact one can show that any element $F$ of $\Omega_{p,R}^{n-1}(D)$ can be decomposed as such a linear combination plus some $\tilde{F}\in \Omega_p^{n-1}(D)$ with $\ast d \tilde F=0$. In particular, if $p\geq \frac{n}{n-1}$ then $\ast d F=0$.
Thus the class $\Omega_{p,R}^{n-1}(D)$ is relatively simple from an analytical point of view and so it is natural to ask which forms in $\Omega_p^{n-1}(D)$ can be approximated by elements in $\Omega_{p,R}^{n-1}(D)$.
The main purpose of the present paper consists in giving a description of the strong and weak closure of the class $\Omega_{p,R}^{n-1}(D)$ for any open domain in $\r^n$ which is bi-Lipschitz equivalent to the open unit $n$-cube $Q_1(0)\subset\r^n$.\\
First we will address the strong closure in the case of the open, unit $n$-dimensional cube $Q_1(0)\subset\r^n$. To this end we introduce the class of $(n-1)$-forms with integer-valued fluxes.\\
For any $F\in \Omega_{L^p_{loc}}^{n-1}$, for any $x_0\in Q_1(0)$ let $\tilde{R}_{F,x_0}\subset (0,r_0)$ be the set of radii $\rho\in (0,r_0)$ such that
\begin{enumerate}
	    \item the hypersurface $\partial Q_{\rho}(x_0)$ consists $\H^{n-1}$-a.e. of Lebesgue points of $F$,
	    \item there holds $|F|\in L^p\big(\partial Q_{\rho}(x_0),\H^{n-1}\big)$.
\end{enumerate}
One can check that $\L^1\big((0,r_{x_0})\smallsetminus \tilde{R}_{F,x_0}\big)=0$.
\begin{dfn}\label{Definition: integer valued fluxes on unit cube}
	Let $p\in [1,\infty]$, let $F\in \Omega_{L^p_{loc}}^{n-1}\cap \Omega_p^{n-1}(Q_1(0),\mu)$ for some Radon measure $\mu$, for any $x_0\in Q_1(0)$ let $\tilde{R}_{F, x_0}$ be defined as above.
    We say that $F$ has \textit{integer-valued fluxes} if for any $x_0\in Q_1(0)$, for $\L^1-$a.e. $\rho\in \tilde{R}_{F, x_0}$ there holds\footnote{Notice that for the associated vector field $V=(*F)^\flat$ condition \eqref{equation: integer fluxes condition, definition} reads
	    \begin{align*}
	        \int_{\partial Q_{\rho}(x_0)}V\cdot \nu_{\partial Q_{\rho}(x_0)}d\mathscr{H}^{n-1}\in \mathbb{Z}.
	    \end{align*}} 
	    \begin{align}\label{equation: integer fluxes condition, definition}
 	            \int_{\partial Q_{\rho}(x_0)}i_{\partial Q_{\rho}(x_0)}^*F\in\z.
	    \end{align}
	The space of $L^p(\mu)$-vector fields with integer valued fluxes will be denoted by $\Omega_{p,\mathbb{Z}}^{n-1}(Q_1(0),\mu)$. The set of radii $\rho\in \tilde{R}_{F, x_0}$ for which \eqref{equation: integer fluxes condition, definition} holds will be denoted by $R_{F, x_0}$
\end{dfn}
%
We will always write $\Omega_{p,\mathbb{Z}}^{n-1}(Q_1(0))$ for $\Omega_{p,\mathbb{Z}}^{n-1}(Q_1(0),\L^{n})$, where $\L^n$ denotes the $n$-dimensional Lebesgue measure.

First of all we observe that (\ref{Equation: degree detects topological singularities}) implies that $\Omega_{p, R}^{n-1}(Q_1(0))\subset \Omega_{p,\mathbb{Z}}^{n-1}(Q_1(0))$.
More general examples of forms in $\Omega_{p,\mathbb{Z}}^{n-1}(Q_1(0))$ can be constructed as follows.
Let again $Y$ be a smooth, closed, oriented and connected $n-1$-dimensional manifold. Let $u\in W^{1, n-1}(Q_1(0), Y)$. Then for any $x_0\in Q_1(0)$, for a.e. $\rho\in (0, 2\dist_\infty(x_0, \partial Q_1(0)))$, $u\big\vert_{\partial Q_\rho(0)}\in W^{1,n-1}(\partial Q_\rho(0),Y)$. Therefore
for any such $\rho$
\begin{align}
\label{Equation: u*omega satisfy Def I.2}
    \int_{\partial Q_\rho (0)}i^\ast_{\partial Q_\rho(x_0)}(u^\ast \omega)=\deg \left(u\big\vert_{\partial Q_\rho(x_0)}\right)\in \mathbb{Z}
\end{align}
Notice that $\deg \left(u\big\vert_{\partial Q_\rho(x_0)}\right)$ is well defined (by means of approximation by functions in $W^{1,\infty}(\partial Q_\rho(x_0), Y)$, see \cite{brezis-nirenberg}, Section I.3).\\
We will show that in fact the closure of ${\Omega_{p,R}^{n-1}(Q_1(0)})$ in $\Omega_{p}^{n-1}(Q_1(0))$ is exactly $\Omega_{p, \mathbb{Z}}^{n-1}(Q_1(0))$. More precisely we have
\begin{thm}
\label{Theorem: strong approximation for vector fields}
Let $n\in \mathbb{N}$ and let $p\in [1,\infty)$. Let $F\in \Omega_{p,\mathbb{Z}}^{n-1}(Q_1(0),\mu)$. Then we have
\begin{enumerate}
    \item if $q\in[0,1]$ and $\displaystyle{p\in\Big[1,\frac{n}{n-1}\Big)}$, then there exists a sequence $\{F_k\}_{k\in \mathbb{N}}$ in $\Omega_{p,R}^{n-1}(Q_1(0))$ such that $F_k\to F$ in $\Omega_{p}^{n-1}(Q_1(0))$ as $k\to\infty$.
    \item if $q\in(-\infty,0]$ and $\displaystyle{p\in\Big[\frac{n}{n-1},+\infty\Big)}$, then $\ast dF=0$.
\end{enumerate}
\end{thm}
The reason why we have introduced the weighted measures $\mu=f\,\L^n$ for $q\neq 0$ is that forms belonging to $\Omega_{p,\z}^{n-1}(Q_1(0),\mu)$ appear naturally in the proof of Corollary \ref{Corollary: strong approximation on boundary of cube}. Nevertheless, we advise the reader to assume $q=0$ (i.e. $\mu=\L^n)$ throughout section 2 at a first reading of the present paper. This allows to skip many technicalities without losing formality, since all the results of this paper are independent on Corollary \ref{Corollary: strong approximation on boundary of cube}.

With the help of Theorem \ref{Theorem: strong approximation for vector fields} we will get another characterization of the $L^p$-closure of $\Omega_{p,R}^{n-1}(Q_1(0))$. For this we recall the following definition (compare with \cite[Section II]{brezis-coron-lieb}):
\begin{dfn}[Connection and minimal connection]
Let $M\subset\r^n$ be any embedded Lipschitz $m$-dimensional submanifold of $\r^n$ (with or without boundary) such that $\overline{M}$ is compact as a subset of $\r^n$.  
 
  A $1$-dimensional current $I\in\mathcal{R}_1(M)$ is said to be a \textit{connection} for (the singular set of) $F$ if $\m(I)<+\infty$ and $\partial I=*dF$ in $(W^{1,\infty}_0(M))^\ast$.
  
 A $1$-dimensional current $L\in\mathcal{R}_1(M)$ is said to be a \textit{minimal connection} for (the singular set of) $F$ if it is a connection for $F$ and
 \begin{align*}
     \mathbb{M}(L)=\inf_{\substack{T\in\mathcal{D}_1(M)\\\partial T=*dF}}\mathbb{M}(T).
 \end{align*}
\end{dfn}
We will see in Corollary \ref{Corollary: existence of connections implies existence of minimal connection} that $F$ admits a connection if and only if it admits a minimal connection.\\
Here is the characterization of $\Omega_{p,\z}^{n-1}(Q_1(0))$ in terms of minimal connections:
\begin{thm}\label{Theorem: characterization of the integer valued fluxes class, introduction}
	Let $n\in \mathbb{N}_{>0}$, let $p\in [1,+\infty)$. Let $F\in\Omega_p^{n-1}(Q^n_1(0))$. Then, the following are equivalent:
	\begin{enumerate}
		\item there exists $L\in\R_1(Q_1(0))$ such that $\partial L=*dF$ in $(W^{1,\infty}_0(Q_1(0)))^\ast$.
		\item for every Lipschitz function $f:\overline{Q_1(0)}\rightarrow [a,b]\subset\r$ such that $f\vert_{\partial Q_1(0)}\equiv b$, we have 
		\begin{align*}
		\int_{f^{-1}(t)}i_{f^{-1}(t)}^*F\in\z, \qquad\mbox{ for $\L^1$-a.e. } t\in [a,b];
		\end{align*}
		\item $F\in\Omega_{p,\z}^{n-1}(Q_1(0))$.
	\end{enumerate}
\end{thm}
In other words, $F\in\Omega_{p,\z}^{n-1}(Q_1(0))$ if and only if $F$ admits a (minimal) connection. This characterization  allows to generalize the definition of the class $\Omega_{p,\z}^{n-1}(Q_1(0))$ to general Lipschitz domains:
\begin{dfn}\label{Definition: integer valued fluxes on generic domain}
Let $M\subset\r^n$ be any embedded Lipschitz $m$-dimensional submanifold of $\r^n$ (with or without boundary).
We define
\begin{align*}
    \Omega_{p,\z}^{n-1}(M):=\{F\in\Omega_p^{n-1}(M) \mbox{ s.t. } \exists\, L\in\mathcal{R}_1(M) \mbox{ connection for } F\}.
\end{align*}
\end{dfn}
Notice that if $M=Q_1(0)$, Definition \ref{Definition: integer valued fluxes on unit cube} and Definition \ref{Definition: integer valued fluxes on generic domain} coincide by Theorem \ref{Theorem: characterization of the integer valued fluxes class, introduction}.\\
We will deduce from the previous results that the approximation result can be extended to any open domain which is bi-Lipschitz equivalent to $Q_1(0)$ of $\partial Q_1(0)$ (see Theorem \ref{Theorem: strong approximation for general domains}).\\
%
We mention here two other corollaries of Theorem \ref{Theorem: strong approximation for vector fields}.
\begin{cor}
\label{Corollary: approximation of partial I, introduction}
Let $n\in \mathbb{N}$. Let $I\in \mathcal{R}_1(Q_1^n(0))$ be an integer rectifiable one-current. Then there exists a one form $\omega\in \Omega_{1,\z}^{n-1}(Q_1^n(0))$ such that $\ast d\omega=\partial I$ and $\partial I$ can be approximated in $(W^{1,\infty}_0(Q_1^n(0)))^\ast$ by finite sums of Dirac-deltas with integer coefficients. More precisely, there exist sequences $(P_i)_{i\in \mathbb{N}}$ and $(N_i)_{i\in \mathbb{N}}$ of points in $Q_1^n(0)$ such that
\begin{align*}
    \partial I= \sum_{i\in \mathbb{N}}(\delta_{P_i}-\delta_{N_i})\,\text{ in }(W^{1,\infty}_0(Q_1^n(0)))^\ast\text{ and }\sum_{i\in \mathbb{N}}\lvert P_i-N_i\rvert<\infty.
\end{align*}
Moreover if $I$ is supported on a Lipschitz submanifold $M$ of $\mathbb{R}^n$ compactly contained in $Q_1(0)$, the points in the sequences $(P_i)_{i\in \mathbb{N}}$ and $(N_i)_{i\in \mathbb{N}}$ can be chosen to belong to $M$.
\end{cor}
The next corollary was obtained first by R. Schoen and K. Uhlenbeck (\cite{schoen-uhlenbeck-2}, Section 4) and F. Bethuel and X. Zheng (\cite[Theorem 4]{bethuel-zheng}).
\begin{cor}
\label{Corollary: Bethuel Zheng, introduction}
Let $Q_1(0)\subset \mathbb{R}^2$ be the unit cube in $\mathbb{R}^2$. Let $u\in W^{1,p}(Q_1(0), \mathbb{S}^1)$ for some $p\in (1,\infty)$.\\
If $p\geq 2$, then $u$ can be approximated in $W^{1,p}$ be a sequence of functions in $C^\infty(Q_1(0), \mathbb{S}^1).$\\
If $p<2$, then $u$ can be approximated in $W^{1,p}$ by a sequence of functions in
\begin{align*}
    \mathcal{R}:=\left\{v\in W^{1,p}(Q_1(0),\mathbb{S}^1); v\in C^\infty(Q_1(0)\smallsetminus A, \mathbb{S}^1),\text{ where }A\text{ is some finite set}\right\}.
\end{align*}
\end{cor}
In the second part of the paper we turn our attention to the weak closure of the space $\Omega_{p,R}^{n-1}(D)$ for a domain $D\subset\mathbb{R}^n$ which is bi-Lipschitz equivalent to $Q_1(0)$ (or equivalently of $\Omega_{p,\mathbb{Z}}^{n-1}(D)$). We will show the following.

\begin{thm}[Weak closure]\label{Theorem: weak closure}
    Let $n\in \mathbb{N}_{\geq 2}$, $p\in(1,+\infty)$ and $D\subset\r^n$ be any open and bounded domain in $\r^n$ which is bi-Lipschitz equivalent to the $n$-dimensional unit cube $Q^n_1(0)$. 
    Then, the space $\Omega_{p,\z}^{n-1}(D)$ is weakly sequentially closed in $\Omega_p^{n-1}(D)$.
\end{thm}
Notice that, by Theorem \ref{Theorem: strong approximation for vector fields} (or more generally by Theorem \ref{Theorem: strong approximation for general domains}), the statement of Theorem \ref{Theorem: weak closure} is trivial for $p\in[n/(n-1),+\infty)$. Thus, we just need to provide a proof in case $p\in(1,n/(n-1))$. 

We first treat the case of the open unit $n$-cube $Q_1(0)\subset\r^n$ by exploiting the characterization of  $\Omega_{p,\z}^{n-1}(Q_1(0))$ given by Theorem \ref{Theorem: characterization of the integer valued fluxes class, introduction} and a suitable slice distance \`a la Ambrosio-Kirchheim (see \cite{ambrosio-kirchheim} and \cite{hardt-riviere}). We then address the general case by standard arguments (see Remark \ref{Remark: weak closure for bounded domains}).\\
We remark that the case $n=1$ is different. In fact for any interval $I\subset \mathbb{R}$ there holds $\overline{\Omega_{p,R}^0(I)}^{\Omega_p^0(I)}=\Omega_p^0(I)$ (see Lemma \ref{Lemma: weak closure for n=1}).\\

Our main motivation to look at forms (instead of maps) with finitely many integer topological singularities is the need of developing geometric measure theory for principal bundles in order to face the still deeply open questions arising in the study of $p$-Yang-Mills lagrangians.

Let $p\in[1,+\infty)$ and $G$ be any compact matrix Lie group. Consider the trivial $G$-principal bundle on $B^n$ given by $\text{pr}_1:P:=B^n\times G\to B^n$, where $\text{pr}_1$ is the canonical projection on the first factor. The $p$\textit{-Yang-Mills lagrangian} on $P$ is given by
\begin{align*}
    \text{YM}_p(A):=\int_{B^n}\lvert F_A\rvert^p\, d\L^n=\int_{B^n}\lvert dA+A\wedge A\rvert^p\, d\L^n,\qquad\forall\,A\in\Omega^1(B^n,\mathfrak{g}),
\end{align*}
where $\mathfrak{g}$ denotes the Lie algebra of $G$.

As it is described in \cite{kessel-riviere}, the reason why we aim to extend the set of the by now classical Sobolev connections is purely analytic and justified by issues arising in the application of the direct method of calculus of variations to $p$-Yang-Mills Lagrangians. On the other hand, the need to extend the notion of bundles in order to allow more and more singularities to appear has already been faced in many geometric applications, which brought to the introduction of coherent and reflexive sheaves in gauge theory (see \cite{kobayashi-nomizu-1},\cite{kobayashi-nomizu-2}).\\

Notice that all results mentioned above can be formulated in terms of vector fields: for any $F\in \Omega_{p}^{n-1}(D)$ we can consider the associated vector field $V_F:=(\ast F)^\flat$. In fact for the proof of some of the results we preferred to work with vector fields instead of $n-1$-forms.\\

\subsection{Related literature and open problems}
Theorem \ref{Theorem: strong approximation for vector fields} was firstly announced to hold for a $3$-dimensional domain in \cite{kessel-riviere} and a full proof of the $3$-dimensional case was eventually given by the first author in \cite{caniato}. Some form of the $2$-dimensional case was treated in \cite{petrache-2d}, where M. Petrache proved that a strong approximability result holds for $1$-forms admitting a connection both on the $2$-dimensional disk and on the $2$-sphere. In both cases, the proof that we give here is more general and simple. The $3$-dimensional version of Theorem \ref{Theorem: weak closure} was already treated by M. Petrache and T. Rivi\`ere in \cite{petrache-riviere-abelian}. Nevertheless, here we took the opportunity to present the arguments in a more detailed and complete way. Both Theorem \ref{Theorem: strong approximation for vector fields} and Theorem \ref{Theorem: weak closure} in dimension $n\neq 2,3$ are instead completely new. 

The first open problems that relate directly to our results are linked to the celebrated Yang-Mills Plateau problem. Indeed, the weak sequential closedness of the class $\Omega_{p,\z}^2(B^3)$ implies that such forms behave well-enough to be considered as suitable "very weak" curvatures for the resolution of the $p$-Yang-Mills Plateau problem for $U(1)$-bundles on $B^3$ (see the introduction of \cite{petrache-riviere-abelian} for further details). The question to address would be if and how we can exploit the same kind of techniques in order to face the existence and regularity issues linked to the so called "non abelian case" (i.e. the case of bundles having a non abelian structure group) in supercritical dimension. An interesting proposal in this sense is due to M. Petrache and T. Rivi\`ere and can be found in \cite{petrache-riviere-non-abelian}, where a suitable class of weak connections in the supercritical dimension $5$ is introduced and studied.  

One could also hope that the technique presented in this paper could be adapted to show $L^p$-closedness (weak and strong) of classes of differential forms exhibiting "integer fluxes" properties similar to the one described in Definition \ref{Definition: integer valued fluxes on unit cube}. As an example we define here the class $\Omega_{p,H}^n(Q_1^{2n}(0))$ of differential forms with "Hopf singularities".\\
Recall that for any $n\in \mathbb{N}_{\geq 1}$, for any smooth map $f:\mathbb{S}^{2n-1}\to\mathbb{S}^n$ the \textit{Hopf invariant} of $f$ is defined as follows: let $\omega$ be the standard volume form on $\mathbb{S}^n$. Let $\alpha\in \Omega^{n-1}(\mathbb{S}^{2n-1})$ be such that $f^\ast \omega=d\alpha$. Then the Hopf invariant of $f$ is given by
\begin{align*}
    H(f):=\int_{\mathbb{S}^{2n-1}}\alpha\wedge d\alpha.
\end{align*}
One can show that $H(f)\in \mathbb{Z}$ and that it is independent of the choice of $\alpha$ (see \cite{bott-tu}, Proposition 17.22).
In the spirit of Definition \ref{Definition: integer valued fluxes on unit cube} we say that a form $F\in \Omega_{p}^{n}(Q_1^{2n}(0))$ belongs to $\Omega_{p,H}^n(Q_1^{2n}(0))$ for some $\displaystyle{p\geq 2-\frac{1}{n}}$ if there exists $A\in\Omega_{W^{1,p}}^{n-1}(Q_1^{2n}(0))$ such that $dA=F$ and if for every $x_0\in Q_1^{2n}(0)$ there exists a set $R_{F,x_0}\subset (0,r_{x_0})$, with $r_{x_0}:=2\dist_{\infty}(x_0,\partial Q_1^{2n}(0))$ such that:
\begin{enumerate}
	\item $\L^1\big((0,r_{x_0})\smallsetminus R_{F,x_0}\big)=0$;
	\item for every $\rho\in R_{F,x_0}$, the hypersurface $\partial Q_{\rho}^{2n}(x_0)$ consists $\H^{2n-1}$-a.e. of Lebesgue points of $F$, $A$ and $\nabla A$ (the matrix of all the partial derivatives of the components of $A$);
	\item for every $\rho\in R_{F,x_0}$ we have $|F|,\lvert A\rvert,\lvert\nabla A\rvert\in L^p\big(\partial Q_{\rho}^{2n}(x_0),\H^{n-1}\big)$;
	\item for every $\rho\in R_{F,x_0}$ it holds that 
	    \begin{align*}
	        \int_{\partial Q_{\rho}^{2n}(x_0)}i_{\partial Q_{\rho}^{2n}(x_0)}^*(A\wedge F)\in\z.
	    \end{align*}
\end{enumerate}
Notice that if $u\in W^{1,2n-1}(Q_1^{2n}(0),\mathbb{S}^n)$, then $u^\ast\omega\in \Omega_{p,H}^n(Q_1^{2n}(0))$.\\

\subsection{Organization of the paper}
The paper is organized as follows. Section 2 is dedicated to the strong $L^p$-closure of $\Omega_{p,\mathbb{Z}}^{n-1}(Q_1(0),\mu)$. First we present some preliminary and rather technical lemmata (Sections 2.1-2.3), then we give a proof of Theorem \ref{Theorem: strong approximation for vector fields} (Section 2.4). In section 2.5 we show the characterization of $\Omega_{p,\mathbb{Z}}^{n-1}(Q_1(0))$ in terms of (minimal) connections (Theorem \ref{Theorem: characterization of the integer valued fluxes class, introduction}). In Section 2.6 we exploit this result to extend the approximation result to other Lipschitz manifolds, and in particular to $\partial Q^{n}_1(0)$. Finally in Section 2.7 we prove Corollary \ref{Corollary: approximation of partial I, introduction} and Corollary \ref{Corollary: Bethuel Zheng, introduction}.\\
In section 3 we discuss the weak $L^p$ closure of $\Omega_{p,R}^{n-1}(Q_1(0))$. First we will introduce a slice distance \`a la Ambrosio-Kirchheim, first on spheres (Section 3.1) and then on cubes (Section 3.2). In Section 3.3 we discuss some of the properties of the slice distance and in Section 3.4 we use it to obtain a proof of Theorem \ref{Theorem: weak closure}. We will also discuss briefly the special case $n=1$.
\subsection{Aknowledgements}
We are grateful to prof. Tristan Rivi\`ere for encouraging us to work on this subject, for his insights and for the helpful discussions. We would also like to thank Federico Franceschini for the interesting discussions.\\
This work has been supported by the Swiss National Science Foundation (SNF 200020\_192062).

\subsection{Notation}
Let $M^m\subset\r^n$ be any $m$-dimensional, embedded Lipschitz submanifold of $\r^n$ (with or without boundary) such that $\overline{M}$ is compact as a subset of $\r^n$.
\begin{itemize}
    \item We denote by $i_M:M\hookrightarrow\r^n$ the usual inclusion map. 
    \item We always assume that $M$ is endowed with the $L^\infty$-Riemannian metric given by $g_M:=i_M^*g_e$, where $g_e$ denotes the standard euclidean metric on $\r^n$. 
    \item For every $k=1,...,m$, we define the following spaces of \textit{smooth $k$-forms} on $M$:
    \begin{align*}
        \Omega^k(M)&:=\{i_M^*\omega \mbox{ s.t. } \omega \mbox{ is a smooth } k \mbox{-form on } \r^n\},\\
        \D^k(M)&:=\{\omega\in\Omega^k(M) \mbox{ s.t. } \spt(\omega)\subset\subset M\}.
    \end{align*}
    For every $p\in[1,+\infty]$, we denote by $\Omega_p^k(M)$ and $\Omega_{W^{1,p}}^k(M)$ the completions of $\Omega^k(M)$ with respect to the usual $L^p$-norm and $W^{1,p}$-norm respectively. We call $\Omega_p^k(M)$ the space of \textit{$L^p$ $k$-forms} on $M$ and $\Omega_{W^{1,p}}^k(M)$ the space of \textit{$W^{1,p}$ $k$-forms} on $M$. 
    \item By the symbol "$\ast$", we denote the Hodge star operator associated with the metric $g_M$ on $M$. By "$\flat$" and "$\sharp$" we denote the usual musical isomorphisms associated with the metric $g_M$. Recall that, under this notation, the map
    \begin{align}\label{musical isomorphism}
        \Omega_p^{n-1}(M)\ni\omega\mapsto(\ast\omega)^{\sharp}\in L^p(M,\r^n)
    \end{align}
    gives an isomorphism onto its image. Exploiting this fact, we frequently identify $(n-1)$-forms with vector fields on $M$. 
\end{itemize}
Let $x_0\in\r^n$ and $\rho>0$. 
\begin{itemize}
    \item We denote by $\lVert\,\cdot\,\rVert_{\infty}:\r^n\rightarrow[0,+\infty)$ the following norm on $\r^n$:
    \begin{align*}
        \lVert x\rVert_{\infty}:=\max_{j=1,...,n}\lvert x_j\rvert.
    \end{align*}
    We denote by $d_{\infty}:\r^n\times\r^n\to\r$ the distance associated to $\lVert\,\cdot\,\rVert_{\infty}$.
    \item We let
    \begin{align*}
        Q_{\rho}^n(x_0):=\bigg\{x\in\r^n \mbox{ s.t. } \lVert x-x_0\rVert_{\infty}<\frac{\rho}{2}\bigg\}
    \end{align*}
    be the open cube in $\r^n$ centered at $x_0$ and having edge-length $\rho$.
    We will sometimes omit the $n$ when the dimension is clear from the context.
    
    We denote by $B^n_{\rho}(x_0)$ the open ball in $\r^n$ centered at $x_0$ with radius $\rho$ (here again we will sometime omit the $n$).
    \item We define 
    \begin{align*}
        \D_k(M)&:=\{k\mbox{-currents on } M\},\\
        \mathcal{M}_k(M)&:=\{T\in\D_k(M) \mbox{ s.t. } \m(T)<+\infty\},\\
        \mathcal{N}_k(M)&:=\{T\in\D_k(M) \mbox{ s.t. } \m(T),\m(\partial T)<+\infty\},\\
        \R_k(M)&:=\{T\in\D_k(M) \mbox{ s.t. } T \mbox{ is integer-multiplicity rectifiable}\},
    \end{align*}
    where $\m(\,\cdot\,)$ denotes the mass of a current (see \cite{krantz_parks-geometric_integration_theory} for further explanations). Moreover, given any $F\in\Omega_1^k(M)$ we define the $(m-k)$\textit{-current associated to} $F$ by
    \begin{align*}
        \langle T_F,\omega\rangle:=\int_{M}F\wedge\omega, \qquad\forall\omega\in\D^{m-k}(M).
    \end{align*}
\end{itemize}
\section{The strong $L^p$-approximation Theorem}
In this section we provide a proof of Theorem \ref{Theorem: strong approximation for vector fields}.
In an attempt to make the proof more accessible, we reformulate the Theorem in terms of vector fields.
For any Radon measure $\mu:=f\L^n$ with $f=\left(\frac{1}{2}-\lVert\,\cdot\,\rVert_\infty\right)^q$, with $q\in (-\infty,1]$ let
\begin{align*}
    L_{R}^p(Q_1(0), \mu):=\{V\in L^p(Q_1(0),\mu) \mbox{ vector field s.t.}\, \ast\hspace{-0.5mm}V^{\flat}\in\Omega_{p,R}^{n-1}(Q_1(0))\}
 \end{align*}
and let
\begin{align*}
    L_{\z}^p(Q_1(0),\mu):=\{V\in L^p(Q_1(0),\mu) \mbox{ vector field s.t. }\ast V^\flat\in \Omega_{p,\mathbb{Z}}^{n-1}(Q_1(0))\}.
\end{align*}
We will sometime write $L^p_{\mathbb{Z}}(Q_1(0))$ for $L^p_{\mathbb{Z}}(Q_1(0),\L^{n})$, where $\L^n$ denotes the $n$-dimensional Lebesgue measure.
\begin{thm}\label{Theorem: strong approximation for vector fields, now really for vector fields} 
Let $V\in L_{\z}^p(Q_1^n(0),\mu)$. The following facts hold:
\begin{enumerate}
    \item if $q\in[0,1]$ and $p\in\big[1,n/(n-1)\big)$, then there exists a sequence $\{V_k\}_{k\in\n}\subset L_{R}^p(Q_1(0),\mu)$ such that $V_k\rightarrow V$ strongly in $L^p(Q_1^n(0),\mu)$;
    \item if $q\in(-\infty,0]$ and $p\in\big[n/(n-1),+\infty\big)$, then $\operatorname{div}(V)=0$ distributionally on $Q_1^n(0)$.
\end{enumerate}
\end{thm}
The case $n=1$ is particularly easy and is treated in Lemma \ref{Lemma: caso n=1}. For the proof in the case $n\ge 2$ we follow the ideas of \cite{petrache-2d} and \cite{caniato}. We present here a plan of the proof, reducing to the case $q=0$ for simplicity and without losing generality.

First of all we show that for any $\varepsilon>0$ it is possible to decompose $Q_1(0)$ into cubes $Q$ of size $\varepsilon$ (plus a negligible rest) so that $$\int_{\partial Q}V\cdot\nu_{\partial Q}\, d\H^{n-1}\in \mathbb{Z}$$ and so that the number of cubes where the integral is different from zero is controlled (Section 2.1). We will then show that $V$ can be approximated on the boundaries of the small cubes $Q$ by smooth vector fields $(V_\varepsilon)_{\varepsilon>0}$ with similar properties (Section 2.2). In Section 2.3 we show that the vector fields $V_\varepsilon$ can be extended inside the cubes $Q$ in such a way that the extension $\tilde{V}_\varepsilon$ has a finite number of singularities in $Q$ (more precisely $\tilde{V}_\varepsilon\vert_{Q}\in L^p_R(Q)$) and is close to $V$ in $L^p(Q)$. In section 2.4 we will combine the previous elements to show that the approximating fields constructed above (up to some shifting and smoothing) satisfy the claim of the Theorem.
\subsection{Choice of a suitable cubic decomposition}
Fix any $\eps\in(0,1/4)$ and $a\in Q_{\eps}(0)$. Let
\begin{align*}
    q_{\eps}&:=\max\,\{q\in\n \mbox{ s.t. } \eps q\le 1-\eps\},\\
    C_{\eps}&:=\bigg\{\bigg(j+\frac{1}{2}\bigg)\eps-\frac{1}{2}, \mbox{ with } j=1,...,q_{\eps}-1\bigg\}^n,\\
    \mathscr{C}_{\eps,a}&:=\big\{Q_{\eps}(x)+a,\, \mbox{ with } x\in C_{\eps}\big\}.
\end{align*}
We say that $\mathscr{C}_{\eps,a}$ is the \textit{cubic decomposition} of $Q_1(0)$ with origin in $a$ and mesh thickness $\eps$.

Let
\begin{align*}
    \mathscr{F}_{\eps,a}&:=\big\{F \,\vert\, F \mbox{ is an } (n-1)\mbox{-dimensional face of } \partial Q \mbox{, for some open cube } Q\in\mathscr{C}_{\eps,a}\big\},\\
    S_{\eps,a}&:=\bigcup_{F\in\mathscr{F}_{\eps,a}}F.
\end{align*}
We say that $S_{\eps,a}$ is the $(n-1)$\textit{-skeleton} of the cubic decomposition $\mathscr{C}_{\eps,a}$.
\allowdisplaybreaks
\begin{lem}[Choice of the cubic decomposition]\label{Lemma: cubic decomposition}
Let $n\in \mathbb{N}_{>0}$. Let $V\in L_{\z}^p(Q_1(0),\mu)$ where $\mu:=f\L^n$ with 
\begin{align*}
    f:=\bigg(\frac{1}{2}-\lVert\,\cdot\,\rVert_{\infty}\bigg)^q
\end{align*}
for some $q\in(-\infty,1]$. Then, there exists a subset $E_V\subset (0,1)$ satisfying the following properties:
\begin{enumerate}
    \item $\L^1\big((0,1)\smallsetminus E_V\big)=0$;
    \item for every $\eps\in E_V$, there exists $a_{\eps}\in Q_{\eps}(0)$ such that $\eps\in R_{V,c_Q}$ for every $Q\in\mathscr{C}_{\eps,a_{\eps}}$ and
    \begin{align}\label{estimate cubic decomposition}
        \lim_{\substack{\eps\in E_V\\ \eps\rightarrow 0^+}}\eps\bigg(\sum_{Q\in\mathscr{C}_{\eps,a_{\eps}}}\int_{\partial Q}\lvert V-(V)_Q\rvert^pf(c_Q)\, d\H^{n-1}\bigg)=0,
    \end{align}
    where $(V)_Q=\displaystyle{\fint_Q V}d\L^n$.
\end{enumerate}
\begin{proof}
For any $x\in Q_1(0)$ let $R_{V,x}:=R_{V^\flat,x}$ (see Definition \ref{Definition: integer valued fluxes on unit cube}).
By assumption
\begin{align}\label{equation: bad set has zero mass}
    \int_0^1 \int_{Q_{1-\rho}(0)}\mathbbm{1}_{\rho\in R_{V,x}}\,d\L^n\,d\rho=&\int_{Q_1(0)}\L^1(R_{V,x})\,d\L^n=\int_{Q_1(0)}2\dist(x,\partial Q_1(0))\,d\L^n\\\nonumber=&\int_0^1\L^n(Q_{1-\rho}(0))\,d\rho.
\end{align}
For any $\rho\in (0,1)$ let
\begin{align*}
    X_\rho=\{x\in Q_{1-\rho}(0) \mbox{ : }\rho\notin R_{V,x}\},
\end{align*}
then by \eqref{equation: bad set has zero mass} $\L^n(X_\rho)=0$ for a.e. $\rho\in (0,1)$. Now notice that for any $\rho\in (0,1)$
\begin{align*}
    \L^n(X_\rho)\geq \sum_{c\in C_\rho}\int_{Q_\rho(0)}\mathbbm{1}_{X_\rho}(x+c)\,d\L^n=\int_{Q_\rho(0)}\sum_{c\in C_\rho}\mathbbm{1}_{X_\rho}(x+c)\,d\L^n,
\end{align*}
thus for a.e. $\rho\in (0,1)$ we have that for a.e. $a_\rho\in Q_\rho(0)$ there holds $\rho\in R_{V,c_Q}$ for any $Q\in \mathscr{C}_{\rho, a_\rho}$. Let $E_V$ be the set of all such $\rho\in (0,1)$.\\
Now let $\varepsilon\in E_V$.
We claim that
\begin{align}\label{equation: claim n-1 decay}
    I_{\eps}:=\int_{Q_{\eps}(0)}\sum_{Q\in\mathscr{C}_{\eps,a}}\int_{\partial Q}\lvert V(x)-(V)_Q\rvert^pf(c_Q)\, d\H^{n-1}(x)\, d\L^n(a)=o\big(\varepsilon^{n-1}\big), \qquad\mbox{ as } \eps\rightarrow 0^+\text{ in }E_V.
\end{align}
Indeed, let $\mathscr{F}$ be the set of the faces of the cube $Q_1(0)\subset\r^n$ and notice that
   \begin{align*}
       I_{\eps}&=\int_{Q_{\eps}(0)}\sum_{F_0\in\mathscr{F}}\sum_{c\in C_{\eps}}\int_{\eps F_0}\bigg\lvert V(x+c+a)-\fint_{Q_{\eps}(c+a)}V\bigg\rvert^pf(c+a)\, d\H^{n-1}(x)\, d\L^n(a)\\ 
       &=\sum_{F_0\in\mathscr{F}}\sum_{c\in C_{\eps}}\int_{\eps F_0}\int_{Q_{\eps}(0)}\bigg\lvert V(x+c+a)-\fint_{Q_{\eps}(c+a)}V\bigg\rvert^pf(c+a)\, d\L^n(a)\, d\H^{n-1}(x)\\ 
       &=\sum_{F_0\in\mathscr{F}}\sum_{c\in C_{\eps}}\int_{\eps F_0}\int_{Q_{\eps}(c)}\bigg\lvert V(x+y)-\fint_{Q_{\eps}(y)}V\bigg\rvert^pf(y)\, d\L^n(y)\, d\H^{n-1}(x).
   \end{align*}
Observe that for any $c\in C_\varepsilon$, $x\in \partial Q_\varepsilon (0)$
\begin{align*}
    & \int_{Q_\varepsilon(c)}\left\lvert V(x+y)-\fint_{Q_\varepsilon(y)}V(z)\right\rvert^pf(y)d\L^n(y)\\ \leq & \fint_{Q_\varepsilon(0)}\int_{Q_\varepsilon(c)}\lvert V(x+y)-V(z+y)\rvert^p f(y) d\L^n(y)d\L^n(z).
\end{align*}
Thus for any $F_0\in \mathscr{F}$
\begin{align*}
    &\sum_{c\in C_{\eps}}\int_{\eps F_0}\int_{Q_{\eps}(c)}\bigg\lvert V(x+y)-\fint_{Q_{\eps}(y)}V\bigg\rvert^pf(y)\, d\L^n(y)\, d\H^{n-1}(x)\\
    \leq &\int_{\varepsilon F_0}\fint_{Q_\varepsilon(0)}\int_{Q_{1-2\varepsilon}(0)}\lvert V(x+y)-V(z+y)\rvert^p f(y)d\L^n(y)d\L^n(z)d\H^{n-1}(x)\\
    \leq &\,2^{p-1}\int_{\varepsilon F_0}\int_{Q_{1-2\varepsilon}(0)}\lvert V(x+y)-V(y)\rvert^p f(y) d\L^n(y)d\H^{n-1}(x)\\&+2^{p-1}\varepsilon^{n-1}\fint_{Q_\varepsilon (0)}\int_{Q_{1-2\varepsilon}(0)}\lvert V(z+y)-V(y)\rvert^p f(y) d\L^n(y)d\L^n(z)\\
    \leq &\,2^p\varepsilon^{n-1}\sup_{\alpha\in Q_\varepsilon(0)}\lVert V-V(\cdot-\alpha)\rVert^p_{L^p(Q_{1-2\varepsilon}, \mu)}.
\end{align*}
Since $C_c^0(Q_1(0))$ is dense in $L^p(Q_1(0),\mu)$ (see \cite[Theorem 4.3]{maggi}), given any $\delta>0$ we can find $\tilde V\in C_c^0(Q_1(0))$ such that
\begin{align*}
    \int_{Q_1(0)}\lvert V-\tilde V\rvert^p\, d\mu\le\delta.
\end{align*}
Notice that by Taylor's Theorem
\begin{align}\label{estimate on the difference measure}
    \left\lvert\frac{f(x+\alpha)-f(x)}{f(x)}\right\rvert=&\left\vert\frac{(\frac{1}{2}-\lVert x+\alpha\rVert_{\infty})^q-(\frac{1}{2}-\lVert x\rVert_{\infty})^q}{(\frac{1}{2}-\lVert x\rVert)^q}\right\rvert\leq q\lVert \alpha\rVert_\infty \frac{(\frac{1}{2}-\lVert x\rVert_\infty-\frac{\varepsilon}{2})^{q-1}}{(\frac{1}{2}-\lVert x\rVert_\infty)^q}\\\nonumber\leq& q\varepsilon 2^{1-q}\left(\frac{1}{2}-\lVert x\rVert_\infty\right)^{-1}\leq C
\end{align}
for every $x\in Q_{1-\eps}(0)$, $\alpha\in Q_\varepsilon(0)$ and for some constant $C>0$ depending only on $q$. 
Thus
\begin{align*}
    \lVert V-V(\cdot-\alpha)\rVert_{L^p(Q_{1-\varepsilon}(0), \mu)}^p&\le 4^{p-1}\bigg(\int_{Q_{1-\eps}(0)}\lvert V-\tilde V\rvert^p\, d\mu+\int_{Q_{1-\eps}(0)}\lvert\tilde V-\tilde V(\cdot-\alpha)\rvert^p\, d\mu\\
    &+\int_{Q_{1-\eps}(0)}\lvert\tilde V(\cdot-\alpha)-V(\cdot-\alpha)\rvert^p\, d\mu\bigg)\\
    &=4^{p-1}\bigg(2\int_{Q_{1}(0)}\lvert V-\tilde V\rvert^p\, d\mu+\int_{Q_{1-\eps}(0)}\lvert\tilde V-\tilde V(\cdot-\alpha)\rvert^p\, d\mu\\
    &+\int_{Q_{1-\eps}(-\alpha)}\lvert\tilde V-V\rvert^p\,\frac{f(\cdot+\alpha)-f}{f}\, d\mu\bigg)\\
    &\le 4^{p-1}(2+C)\delta+4^{p-1}\int_{Q_{1-\eps}(0)}\lvert\tilde V-\tilde V(\cdot-\alpha)\rvert^p\, d\mu.
\end{align*}
As $\tilde V\in C^0_c(Q_1(0))$
\begin{align*}
    \sup_{\alpha\in Q_\varepsilon(0)}\lVert\tilde V-\tilde V(\cdot-\alpha)\rVert_{L^p(Q_{1-\varepsilon}(0), \mu)}^p\to 0 \text{ as }\varepsilon\to 0^+,
\end{align*}
we have
\begin{align*}
    \limsup_{\eps\to 0^+}\sup_{\alpha\in Q_\varepsilon(0)}\lVert V-V(\cdot-\alpha)\rVert_{L^p(Q_{1-\varepsilon}(0), \mu)}^p\le4^{p-1}(2+C)\delta.
\end{align*}
By letting $\delta\to 0^+$ in the previous inequality we get
\begin{align}\label{Equation: Continuity wrt argument for weighted measure}
    \sup_{\alpha\in Q_\varepsilon(0)}\lVert V- V(\cdot-\alpha)\rVert_{L^p(Q_{1-\varepsilon}(0), \mu)}^p\to 0 \text{ as }\varepsilon\to 0^+,
\end{align}
and claim \eqref{equation: claim n-1 decay} follows.

By Fubini's theorem, for every fixed $\varepsilon\in E_V$ there exists some non-negligible subset  $T_{\varepsilon}\subset Q_{\eps}(0)$ such that for any $a\in T_\eps$ $\varepsilon\in R_{V, c_Q}$ for any $Q\in \mathscr{C}_{\varepsilon,a}$ and
\begin{align*}
    \sum_{Q\in\mathscr{C}_{\eps,a}}\int_{\partial Q}\lvert V-(V)_Q\rvert^pf(c_Q)\ d\H^{n-1}&\le\frac{1}{\varepsilon^n}\int_{Q_{\eps}(0)}\sum_{Q\in\mathscr{C}_{\eps,a}}\int_{\partial Q}\lvert V-(V)_Q\rvert^p\, f(c_Q) d\H^{n-1}\, d\L^n(a)\\
    &=\frac{1}{\varepsilon^n}I_\varepsilon.
\end{align*}
By \eqref{equation: claim n-1 decay}, for every $a\in T_{\varepsilon}$ we have
\begin{align*}
    \sum_{Q\in\mathscr{C}_{\eps,a}}\int_{\partial Q}\lvert V-(V)_Q\rvert^pf(c_Q)\, d\H^{n-1}=o(\eps^{-1})
\end{align*}
as $\eps\rightarrow 0^+$ in $E_V$. The statement follows. 
\end{proof} 
\end{lem}
Fix any $V\in L_{\z}^p(Q_1(0))$ and $\eps\in E_V$. From now on, we will denote simply by $\mathscr{C}_{\eps}$ the cubic decomposition $\mathscr{C}_{\eps,a_{\eps}}$ provided by Lemma \ref{Lemma: cubic decomposition}. Accordingly, the subscript "$a_{\eps}$" will be omitted in any writing referring to such a cubic decomposition. 

Given any $Q\in\mathscr{C}_{\eps}$, we say that $Q$ is a \textit{good cube} if 
\begin{align*}
    \int_{\partial Q}V\cdot\nu_{\partial Q}=0
\end{align*}
and that $Q$ is a \textit{bad cube} otherwise. We denote by $\mathscr{C}_{\eps}^g$ the subfamily of $\mathscr{C}_{\eps}$ made of all the good cubes and by $\mathscr{C}_{\eps}^b$ the one made of all the bad cubes. Moreover, we let
\begin{align*}
\Omega_{\eps}&:=\bigcup_{Q\in\mathscr{C}_{\eps}}Q,\quad\Omega_{\eps}^g:=\bigcup_{Q\in\mathscr{C}_{\eps}^g}Q,\quad\Omega_{\eps}^b:=\bigcup_{Q\in\mathscr{C}_{\eps}^b}Q,\\
S_{\eps}&:=\bigcup_{Q\in\mathscr{C}_{\eps}}\partial Q,\quad S_{\eps}^g:=\bigcup_{Q\in\mathscr{C}_{\eps}^g}\partial Q,\quad S_{\eps}^b:=\bigcup_{Q\in\mathscr{C}_{\eps}^b}\partial Q.
\end{align*}
\begin{lem}\label{Lemma: volume of the bad cubes vanishes at the limit}
Assume that $n\geq 2$. Then, we have
\begin{align*}
    \lim_{\substack{\eps\in E_V\\\varepsilon\rightarrow 0^+}}\eps^n\sum_{Q\in\mathscr{C}_{\eps}}f(c_Q)=0.
\end{align*}
In particular if $q=0$ (and thus $\mu=\L^n$) we have
\begin{align*}
    \lim_{\substack{\eps\in E_V\\\varepsilon\rightarrow 0^+}}\L^n(\Omega_\varepsilon^b)=0.
\end{align*}
\begin{proof}
Notice that by estimate \eqref{estimate on the difference measure}
\begin{align}\label{estimate on f}
    \frac{\lvert f-f(c_Q)\rvert}{f}\le C \qquad\mbox{ on } Q \mbox{ for every } Q\in\mathscr{C}_{\eps} 
\end{align}
for some universal constant $C>0$.
For every bad cube $Q\in\mathscr{C}_{\eps}^b$, it holds that
\begin{align*}
    1\le\bigg\lvert\int_{\partial Q}V\cdot\nu_{\partial Q}\, d\H^{n-1}\bigg\rvert\le\int_{\partial Q}\lvert V\rvert\, d\H^{n-1}.
\end{align*}
By multiplying the previous inequality by $f(c_Q)$ and summing over all the bad cubes, we get
\begin{align*}
    \sum_{Q\in\mathscr{C}_{\eps}^b}f(c_Q)&\le\sum_{Q\in\mathscr{C}_{\eps}^b}\int_{\partial Q}\lvert V\rvert f(c_Q)\, d\H^{n-1}\\
    &\le \bigg(\sum_{Q\in\mathscr{C}_{\eps}^b}\int_{\partial Q}\lvert V\rvert^pf(c_Q)\, d\H^{n-1}\bigg)^{\frac{1}{p}} \bigg(\sum_{Q\in\mathscr{C}_{\eps}^b}\int_{\partial Q}f(c_Q)\, d\H^{n-1}\bigg)^{\frac{1}{p'}}\\
    &= (2n)^{\frac{1}{p'}}\eps^{\frac{n-1}{p'}}\bigg(\sum_{Q\in\mathscr{C}_{\eps}^b}\int_{\partial Q}\lvert V\rvert^p f(c_Q)\, d\H^{n-1}\bigg)^{\frac{1}{p}}\bigg(\sum_{Q\in\mathscr{C}_{\eps}^b}f(c_Q)\bigg)^{\frac{1}{p'}},
\end{align*}
which is equivalent to
\begin{align*}
    \sum_{Q\in\mathscr{C}_{\eps}^b}f(c_Q)&\le (2n)^{p-1}\eps^{(p-1)(n-1)}\sum_{Q\in\mathscr{C}_{\eps}^b}\int_{\partial Q}\lvert V\rvert^p f(c_Q)\, d\H^{n-1}.
\end{align*}
Hence, by the triangle inequality, we get
\begin{align}\label{estimate number of bad cubes}
    \nonumber
    \sum_{Q\in\mathscr{C}_{\eps}^b}f(c_Q)&\le (4n)^{p-1}\eps^{(p-1)(n-1)-1}\bigg(\eps\sum_{Q\in\mathscr{C}_{\eps}^b}\int_{\partial Q}\lvert V-(V)_Q\rvert^p f(c_Q)\, d\H^{n-1}\\
    \nonumber
    &\quad+2n\sum_{Q\in\mathscr{C}_{\eps}^b}\int_{Q}\lvert V\rvert^pf(c_Q)\, d\L^n\bigg)\\
    \nonumber
    &\le (4n)^{p-1}\eps^{(p-1)(n-1)-1}\Bigg(\eps\sum_{Q\in\mathscr{C}_{\eps}}\int_{\partial Q}\lvert V-(V)_Q\rvert^pf(c_Q)\, d\H^{n-1}\\
    \nonumber
    &\quad+2n\sum_{Q\in\mathscr{C}_{\eps}^b}\int_{Q}\lvert V\rvert^p f\, d\L^n+2n\sum_{Q\in\mathscr{C}_{\eps}^b}\int_{Q}\lvert V\rvert^p \frac{f(c_Q)-f}{f}f\, d\L^n\Bigg)\\
    &\le (4n)^{p-1}\eps^{(p-1)(n-1)-1}\Bigg(\eps\sum_{Q\in\mathscr{C}_{\eps}}\int_{\partial Q}\lvert V-(V)_Q\rvert^pf(c_Q)\, d\H^{n-1}\\
    \nonumber
    &\quad+2n(1+C)\int_{\Omega_{\eps}^b}\lvert V\rvert^p f\, d\L^n\Bigg).
\end{align}
Therefore
\begin{align*}
    \eps^n\sum_{Q\in\mathscr{C}_{\eps}^b}f(c_Q)&\le (4n)^{p-1}\eps^{p(n-1)}\Bigg(\eps\sum_{Q\in\mathscr{C}_{\eps}}\int_{\partial Q}\lvert V-(V)_Q\rvert^pf(c_Q)\, d\H^{n-1}\\
    \quad&+2n(1+C)\int_{Q_1(0)}\lvert V\rvert^p f\, d\L^n\Bigg)
\end{align*}
and the statement follows from \eqref{estimate cubic decomposition} (here we need the assumption $n>1$).
\end{proof}
\end{lem}
\begin{rem}\label{Remark: no bad cubes if p is big enough}
  Assume that $p\in\big[n/(n-1),+\infty\big)$. In this case $\varepsilon^{(p-1)(n-1)-1}$ remains bounded as $\varepsilon\to 0^+$. Now by Lemma \ref{Lemma: cubic decomposition}
  \begin{align*}
      \eps\sum_{Q\in\mathscr{C}_{\eps}}\int_{\partial Q}\lvert V-(V)_Q\rvert^p f(c_Q)\, d\H^{n-1}\to 0\quad\text{ as }\varepsilon\to 0^+\text{ in }E_V.
  \end{align*}
 Moreover, by Lemma \ref{Lemma: volume of the bad cubes vanishes at the limit}, we have
  \begin{align*}
      \int_{\Omega_{\eps}^b}f\, d\L^n=(1+C)\eps^n\sum_{Q\in\mathscr{C}_{\eps}^b}f(c_Q)\to 0\quad\text{ as }\varepsilon\to 0^+\text{ in }E_V.
  \end{align*}
  This implies $\L^n(\Omega_\varepsilon^b)\to 0$ as $\varepsilon\to 0^+$ in $E_V$, therefore
   \begin{align*}
      \int_{\Omega_{\eps}^b}\lvert V\rvert^pf\, d\L^n\to 0\quad\text{ as }\varepsilon\to 0^+\text{ in }E_V
  \end{align*}
  by absolute continuity of the integral. Thus it follows from \eqref{estimate number of bad cubes} that 
  \begin{align*}
    \sum_{Q\in\mathscr{C}_{\eps}^b}f(c_Q)\to 0^+\quad\text{ as }\varepsilon\to 0^+\text{ in }E_V.
    \end{align*}
  Let $N_{\eps}^b$ be the number of bad cubes in $\mathscr{C}_{\eps}$. Notice that for $q\le 0$, we have $f\ge 2^{-q}$ on $Q_1(0)$. This implies
    \begin{align*}
    N_{\eps}^b\le 2^{q}\sum_{Q\in\mathscr{C}_{\eps}^b}f(c_Q)\to 0^+\quad\text{ as }\varepsilon\to 0^+\text{ in }E_V.
    \end{align*}
  Since $N_\varepsilon^b\in \mathbb{Z}$ for any $\varepsilon\in E_V$, $N_{\eps}^b=0$ for every $\eps\in E_V$ small enough. Hence, whenever $p\in\big[n/(n-1),+\infty\big)$ and $q\le 0$ we will assume, without losing generality, that there are no bad cubes in our chosen decomposition.
\end{rem}
\subsection{Smoothing on the $(n-1)$-skeleton of the cubic decomposition}
\begin{lem}\label{Lemma: smoothing in the (n-1)-skeleton}
Let $Q\subset\r^n$ be an $n$-dimensional cube with side length $R$. Let $V\in L^p(Q,\mathbb{R}^m)$. Let $\varepsilon>0$. There exists $V_\varepsilon\in C^\infty_c(Q, \mathbb{R}^m)$ such that
\begin{align*}
    \int_Q V_\varepsilon\, d\L^n=\int_Q V\, d\L^n
\end{align*}
and
\begin{align*}
    \lVert V_\varepsilon-V\rVert_{L^p(Q)}<\varepsilon.
\end{align*}
\begin{proof}
    Without loss of generality, we will assume that $Q$ is centered in the origin of $\r^n$. Let $\psi\in C^\infty_c(\frac{1}{2}Q)$ and $r_0\in (1/2,1)$ such that
    \begin{align*}
        \int_Q \psi\, d\L^n=1\quad\text{ and } R^n(1-r_0^n)<\frac{\varepsilon}{\lVert \psi\rVert_{L^p(Q)}}.
    \end{align*}
    Let $r\in (r_0,1)$ be such that
    \begin{align*}
        \lVert V\rVert_{L^p(Q\smallsetminus rQ)}\leq\min\left\{\left(\frac{\varepsilon}{\lVert \psi\rVert_{L^p(Q)}}\right)^\frac{1}{p},\varepsilon\right\}.
    \end{align*}
    Set
    \begin{align*}
        s:=\int_{Q\smallsetminus rQ}V\, d\L^n, \quad \tilde{V}:=\rchi_{rQ}V+s\psi\in L^p(Q,\r^m).
    \end{align*}
    Then
    \begin{align*}
        \int_Q \tilde{V}\, d\L^n=\int_{rQ}V\, d\L^n+\left(\int_{Q\smallsetminus rQ}V\, d\L^n\right)\int_Q \psi\, d\L^n=\int_Q V\, d\L^n.
    \end{align*}
    Moreover
    \begin{align*}
        \lvert s\rvert&=\left\lvert \int_{Q\smallsetminus rQ}V\, d\L^n\right\rvert\leq\lvert Q\smallsetminus rQ\rvert^\frac{1}{p'}\lVert V\rVert_{L^p(Q\smallsetminus rQ)}\\
        &\leq (R^n(1-r^n))^\frac{1}{p'}\left(\frac{\varepsilon}{\lVert\psi\rVert_{L^p(Q)}}\right)^\frac{1}{p}\le \frac{\varepsilon}{\lVert\psi\rVert_{L^p(Q)}}.
    \end{align*}
    Therefore 
    \begin{align*}
       \lVert s\psi\rVert_{L^p(Q)}=\lVert \psi\rVert_{L^p(Q)}\lvert s\rvert\leq \varepsilon 
    \end{align*}
        and, by choice of $\tilde{V}$,
    \begin{align*}
        \lVert V-\tilde{V}\rVert_{L^p(Q)}\leq \lVert s\psi\rVert_{L^p(Q)}+\lVert V\rVert_{L^p(Q\smallsetminus rQ)}\leq 2\varepsilon.
    \end{align*}
    Notice that $\tilde V\big\vert_{Q\smallsetminus rQ}\equiv 0$.\\
    Let $\eta\in C_c^\infty(B_1(0))$ with $\int_{B_1(0)}\eta\, d\L^n=1$.
    For any $\delta>0$ let
    \begin{align*}
        \eta_\delta(x):=\frac{1}{\delta^n}\eta\left(\frac{x}{\delta}\right)\quad\forall\, x\in \mathbb{R}^n.
    \end{align*}
    Choose $\delta_0>0$ such that
    \begin{align*}
        2\delta_0<\operatorname{dist}\big(\partial Q, \partial (rQ)\big)\quad\text{ and }\quad
        \lVert\tilde{V}-\tilde{V}\ast\eta_{\delta_0}\rVert_{L^p(Q)}\leq\varepsilon.
    \end{align*}
    Set $V_\varepsilon:=\tilde{V}\ast\eta_{\delta_0}$.
    Then $V_\varepsilon\in C^\infty_c(Q)$,
    \begin{align*}
        \int_QV_\varepsilon\, d\L^n=\int_{\mathbb{R}^n}\eta_{\delta_0}\, d\L^n\int_Q\tilde{V}\, d\L^n=\int_Q \tilde{V}\, d\L^n=\int_QV\, d\L^n
    \end{align*}
    and
    \begin{align*}
        \lVert V-V_\varepsilon\rVert_{L^p(Q)}\leq \lVert V-\tilde{V}\rVert_{L^p(Q)}+\lVert \tilde{V}-V_\varepsilon\rVert_{L^p(Q)}\leq 3\varepsilon.
    \end{align*}
\end{proof}
\end{lem}
\subsection{Extensions on good and bad cubes}
\begin{lem}[Extension on the good cubes]\label{appendix: extension on the good cubes}
    Let $\Omega\subset\r^n$ be a bounded, connected Lipschitz domain and $p\in [1,\infty)$. Let $f\in L^p(\partial\Omega)$ with
    \begin{align}\label{equation: Assumption for extension on good cubes}
        \int_{\partial\Omega} f\, d\H^{n-1}=0 
    \end{align}
    There exists a vector field $V\in L^p(\Omega)$ such that
    \begin{align}\label{Equation: Integration by parts for new VF V}
        \int_\Omega V\cdot \nabla\varphi\, d\L^n=\int_{\partial \Omega}f\varphi\, d\H^{n-1}\qquad\forall\,\varphi\in C^\infty(\mathbb{R}^n)
    \end{align}
    and
    \begin{align}\label{appendix estimate extension good cubes}
        \int_\Omega \lvert V\rvert^p\, d\L^n\leq C(p,\Omega)\int_{\partial \Omega}\lvert f\rvert^p\, d\H^{n-1}
    \end{align}
    for some constant $C(p,\Omega)$ depending only on $p$ and $\Omega$.\\
    Moreover if $p=1$, then $V\in L^q(\Omega)$ for any $q\in\big[1,\frac{n}{n-1}\big)$.
\end{lem}
\begin{rem}
Observe that (\ref{Equation: Integration by parts for new VF V}) implies that $V$ is a distributional solution of the following Neumann problem
\begin{align*}
    \begin{cases}
        \operatorname{div}(V)=0 &\text{ in }\Omega\\
        V\cdot \nu_{\partial Q}=f&\text{ on }\partial\Omega.
    \end{cases}
\end{align*}
\end{rem}
\begin{proof}
\,\\
\textbf{Step 1:} First we consider the case $p\in (1,\infty)$.\\
    Let $p':=\displaystyle{\frac{p}{p-1}}$.
    For any $u\in W^{1,p'}(\Omega)$ let
    \begin{align*}
        E_p(u)=\frac{1}{p'}\int_\Omega \lvert \nabla u\rvert^{p'}\,d\L^n-\int_{\partial \Omega} fu\,d\H^{n-1}.
    \end{align*}
    Recall that any function $u\in W^{1, p'}(\Omega)$ has a trace in $L^{p'}(\partial \Omega)$, and that the trace operator is continuous.
    Thus for any $v\in W^{1,p'}(\Omega)$ with $\displaystyle{\int_\Omega v=0}$ by Poincar\'e Lemma there holds
    \begin{align*}
        \left\lvert\int_{\partial\Omega} fv\,d\H^{n-1}\right\rvert\leq \lVert f\rVert_{L^p(\partial \Omega)}\lVert v\rVert_{L^{p'}(\partial \Omega)}\leq C(p,\Omega)\lVert f\rVert_{L^p(\partial \Omega)}\lVert \nabla v\rVert_{L^{p'}(\Omega)}
    \end{align*}
    for some constant $C(p,\Omega)$ depending only on $p$ and $\Omega$. In particular the energy $E_p$ is well defined on $W^{1,p'}(\Omega)$.\\
    Let
    \begin{align*}
        \dot W^{1,p'}(\Omega):=\left\{v\in W^{1,p'}(\Omega), \int_\Omega v\,d\L^n=0\right\}
    \end{align*}
    and observe that $E_p$ is strictly convex on $\dot W^{1,p'}(\Omega)$.
    Let $u$ be the unique minimizer of $E_p$ in $\dot W^{1,p'}(\Omega)$. Then\footnote{The argument above shows that \eqref{Equation: V satisfies equation weakly} holds for any $\varphi\in C^\infty(\mathbb{R}^n)$ with $\int_\Omega \varphi d\L^n=0$, but assumption \eqref{equation: Assumption for extension on good cubes} implies that \eqref{Equation: V satisfies equation weakly} remains valid for any $\varphi\in C^\infty(\mathbb{R}^n)$.}
    \begin{align}
    \label{Equation: V satisfies equation weakly}
        \int_\Omega \lvert \nabla u\rvert^{p'-2}\nabla u\cdot \nabla\varphi\,d\L^n=\int_{\partial \Omega}f\varphi\,d\H^{n-1},\qquad\forall \varphi\in C^\infty(\mathbb{R}^n).
    \end{align}
    Moreover, as $u$ is a minimizer of $E_p$, $E_p(u)\leq E_p(0)=0$. It follows that
    \begin{align*}
        \frac{1}{p'}\int_\Omega\lvert\nabla u\rvert^{p'}\,d\L^n\leq \int_{\partial \Omega}fu\,d\H^{n-1}\leq\lVert f\rVert_{L^p(\partial \Omega)}\lVert u\rVert_{L^{p'}(\partial \Omega)}\leq C(p,\Omega)\lVert f\rVert_{L^p(\partial \Omega)}\lVert \nabla u\rVert_{L^{p'}(\Omega)}.
    \end{align*}
    Thus
    \begin{align*}
        \int_\Omega \lvert \nabla u\rvert^{p'}\,d\L^n\leq \left(p' C(p,\Omega)\right)^p\int_{\partial \Omega}\lvert f\rvert^p\,d\H^{n-1}.
    \end{align*}
    Set $V:=\lvert\nabla u\rvert^{p'-2}\nabla u\quad\text{in }\Omega$.
    Then by (\ref{Equation: V satisfies equation weakly})
    \begin{align*}
        \int_\Omega V\cdot\nabla \varphi\,d\L^n=\int_{\partial\Omega} f\varphi\,d\H^{n-1}\quad\forall\varphi\in C^\infty(\mathbb{R}^n)
    \end{align*}
    and
    \begin{align*}
        \int_{\Omega}\lvert V\rvert^p\,d\L^n=\int_{\Omega}\lvert \nabla u\rvert^{p'}\,d\L^n\leq \left(p' C(p,\Omega)\right)^p\int_{\partial \Omega}\lvert f\rvert^p\,d\H^{n-1}.
    \end{align*}
    \textbf{Step 2:} Next we consider the case $p=1$.\\
    Let $s>n$. For any $u\in W^{1,s}(\Omega)$ let
    \begin{align*}
        E_s(u)=\frac{1}{s}\int_\Omega\lvert \nabla u\rvert^s\,d\L^n-\int_{\partial \Omega}f u\,d\H^{n-1}.
    \end{align*}
    Notice that $E_s$ is well defined and strictly convex in $\dot{W}^{1,s}(\Omega)$.\\
    Recall the Sobolev embedding
    \begin{align*}
        W^{1,s}(\Omega)\hookrightarrow C^{0,\alpha}(\overline{\Omega})
    \end{align*}
    for $\alpha=1-\frac{n}{s}$.
    Then for any $u\in W^{1,s}(\Omega)$ the trace of $u$ on $\partial \Omega$ lies in $C^{0,\alpha}(\partial \Omega)$ and if $\int_\Omega u\,d\L^n=0$.
    Poincar\'e inequality implies
    \begin{align*}
        \lVert u\rVert_{L^{\infty}(\partial \Omega)}\leq C(s,\Omega)\lVert \nabla u\rVert_{L^s(\Omega)}
    \end{align*}
    for some constant $C(s,\Omega)$ depending only on $s$ and $\Omega$.\\
    Let $u$ be the unique minimizer of $E_s$ in $\dot{W}^{1,s}(\Omega)$.
    Then since $E_s(u)\leq E_s(0)=0$ there holds
    \begin{align*}
        \frac{1}{s}\int_\Omega\lvert \nabla u\rvert^s\,d\L^n&\leq \int_{\partial\Omega}fu\,d\H^{n-1}\\
        &\leq\lVert f\rVert_{L^1(\partial \Omega)}\lVert u\rVert_{L^\infty(\partial\Omega)}\leq C(s,\Omega)\lVert f\rVert_{L^1(\partial\Omega)}\lVert \nabla u\rVert_{L^s(\Omega)}.
    \end{align*}
    Therefore
    \begin{align*}
        \left(\int_\Omega \lvert\nabla u\rvert^s\,d\L^n\right)^{\frac{s-1}{s}}\leq s C(s,\Omega) \int_{\partial \Omega}\lvert f\rvert\,d\H^{n-1}.
    \end{align*}
    Moreover, since $u$ is a minimizer of $E_s$,
    \begin{align*}
        \int_\Omega\lvert\nabla u\rvert^{s-2}\nabla u\cdot\nabla\varphi\,d\L^n=\int_{\partial \Omega}f\varphi\,d\H^{n-1}\quad\forall \varphi\in C^\infty(\mathbb{R}^n).
        \end{align*}
    Similarly as before set $V:=\lvert \nabla u\rvert^{s-2}\nabla u\quad\text{in }\Omega$.
    Then
    \begin{align*}
        \int_{\Omega}\lvert V\rvert\,d\L^n=\int_\Omega \lvert\nabla u\rvert^{s-1}\,d\L^n\leq \L^{n}(\Omega)^\frac{1}{s}\left(\int_\Omega\lvert\nabla u\rvert^s\,d\L^n\right)^\frac{s-1}{s}\leq sC(s,\Omega)\L^{n}(\Omega)^\frac{1}{s}\int_{\partial\Omega}\lvert f\rvert\,d\H^{n-1}.
    \end{align*}
    Moreover
    \begin{align*}
        \int_\Omega \lvert V\rvert^{\frac{s}{s-1}}\,d\L^n=\int_\Omega \lvert \nabla u\rvert^s\,d\L^n<\infty
    \end{align*}
\end{proof}
\begin{rem}\label{appendix scaling argument}
    Let $Q\subset \mathbb{R}^n$ be the unit cube and let $C(p,Q)$ be the corresponding constant in (\ref{appendix estimate extension good cubes}). By an easy scaling argument one sees that for any $\varepsilon>0$ one can choose $C(p,\varepsilon Q)=\varepsilon C(p,Q)$.
\end{rem}
\begin{lem}[Extension on the bad cubes]\label{Lemma: extension on the bad cubes}
Let $Q:=Q_{\eps}(c_Q)\subset\r^n$ and $r(x):=\lVert x-c_Q\rVert_{\infty}$, for every $x\in\r^n$. 

Consider any vector field $V:Q\rightarrow\r^n$ having the form
\begin{align*}
    V(x):=\frac{1}{2^{n-1}}f\bigg(\frac{\eps}{2}\frac{x-c_Q}{r(x)}+c_Q\bigg)\frac{x-c_Q}{r(x)^n} \qquad\forall\,x\in Q,
\end{align*}
for some $f\in L^{\infty}(\partial Q)$. Then, the following facts hold:
\begin{enumerate}
    \item $V\in L^{p}(Q)$ for every $p\in\big[1,n/(n-1)\big)$;
    \item for some constant $C(n,p)>0$ depending only on $n$ and $p$ we have
        \begin{align}\label{estimate bad cubes}
            \int_Q\lvert V\rvert^p\, d\L^n\le\eps C(n,p)\int_{\partial Q}\lvert f\rvert^p\, d\H^{n-1},
        \end{align}
    \item for every $\varphi\in C^{\infty}(\r^n)$ we have
        \begin{align}\label{divergence bad cubes}
            \int_{Q}V\cdot\nabla\varphi\, d\L^n=\int_{\partial Q}f\varphi\, d\H^{n-1}-\bigg(\int_{\partial Q}f\, d\H^{n-1}\bigg)\varphi(c_Q).
        \end{align}
\end{enumerate}

\begin{proof}
Without losing generality, we assume that $\eps=1$ and $c_Q=0$. First, notice that $r:\r^n\rightarrow\r$ is a Lipschitz map such that $\lvert\nabla r(x)\rvert=1$, for a.e. $x\in\r^n$. Moreover, since all the norms are equivalent on $\r^n$ there exists a constant $\tilde C(n)>0$ depending only on $n$ such that $\lvert x\rvert\le\tilde Cr(x)$, for a.e. $x\in\r^n$. Now choose any $p\in\big[1,n/(n-1)\big)$. By coarea formula we have
\begin{align*}
    \int_Q\lvert V\rvert^p\, d\L^n&\le\frac{\tilde C^p}{2^{(n-1)p}}\int_{0}^\frac{1}{2}\frac{1}{\rho^{(n-1)p}}\int_{\partial Q_{2\rho}(0)}\bigg|f\bigg(\frac{x}{2\rho}\bigg)\bigg|^p\, d\H^{n-1}(x)\,d\rho\\
    &=\frac{\tilde C^p}{2^{(n-1)(p-1)}}\bigg(\int_{0}^\frac{1}{2}\frac{1}{\rho^{(n-1)(p-1)}}\,d\rho\bigg)\bigg(\int_{\partial Q}\lvert f(y)\rvert^p\, d\H^{n-1}(y)\bigg)\\ 
    &=C\int_{\partial Q}\lvert f\rvert^p\, d\H^{n-1},
\end{align*}
with
\begin{align*}
    C=C(n,p):=\frac{\tilde C^p}{2^{(n-1)(p-1)}} \int_{0}^\frac{1}{2}\frac{1}{\rho^{(n-1)(p-1)}}\,d\rho<+\infty.
\end{align*}
Hence, $1.$ and $2.$ follow in once. We remark that the condition $p\in\big[1,n/(n-1)\big)$ is needed in order to guarantee the convergence of the integral in $\rho$. 

We still need to prove 3. Pick any $\varphi\in C^{\infty}(\r^n)$. By the coarea formula we have 
\begin{align*}
    \int_{Q}V\cdot\nabla\varphi\, d\L^n&=\frac{1}{2^{n-1}}\int_{0}^\frac{1}{2}\frac{1}{\rho^{n}}\int_{\partial Q_{2\rho}(0)}f\bigg(\frac{x}{2\rho}\bigg)\big(x\cdot\nabla\varphi(x)\big)\, d\H^{n-1}(x)\,d\rho\\ 
    &=2\int_{0}^\frac{1}{2}\int_{\partial Q}f(y)\big(y\cdot\nabla\varphi(2\rho y)\big)\, d\H^{n-1}(y)\,d\rho\\ 
    &=\int_{\partial Q}f(y)\int_{0}^\frac{1}{2}\frac{d}{d\rho}\big(\varphi(2\rho y)\big)\,d\rho\, d\H^{n-1}(y)\\
    &=\int_{\partial Q}f\varphi\, d\H^{n-1}-\bigg(\int_{\partial Q}f\, d\H^{n-1}\bigg)\varphi(0)
\end{align*}
and $3.$ follows.
\end{proof}
\end{lem}

\subsection{Proof of Theorem 1.1}
We are finally ready to prove Theorem \ref{Theorem: strong approximation for vector fields, now really for vector fields}. 

\begin{proof}
Let $V\in L_{\z}^p(Q_1(0))$ and let $\eps\in E_V$ (constructed in Lemma \ref{Lemma: cubic decomposition}).
First, we notice that by using Lemma \ref{Lemma: smoothing in the (n-1)-skeleton} separately on every face $F\in\mathscr{F}_{\eps}$ we can build a vector field $V_{\eps}\in C^{\infty}(S_{\eps})$ such that 
\begin{align*}
    \int_{\partial Q}V_{\eps}\cdot\nu_{\partial Q}\, d\H^{n-1}=\int_{\partial Q}V\cdot\nu_{\partial Q}\, d\H^{n-1}\qquad\forall\,Q\in\mathscr{C}_{\eps}
\end{align*}
and
\begin{align*}
    \sum_{Q\in \mathscr{C}_\varepsilon}\int_{\partial Q}\lvert V_\varepsilon-V\rvert^p f(c_Q)d\mathscr{H}^{n-1}<\varepsilon.
\end{align*}
Let $\tilde V_{\eps}$ be the vector field defined $\L^n$-a.e. on $\Omega_{\eps}$ as follows:
\begin{enumerate}
    \item if $Q\in \mathscr{C}_{\eps}$ is a good cube, then we let $\tilde V_{\eps}:=W_{\eps}+(V)_Q$ on $Q$, where $W_{\eps}$ is the extension of the datum $f:=\big(V_{\eps}-(V)_Q\big)\cdot\nu_{\partial Q}$ given by Lemma \ref{appendix: extension on the good cubes} (notice that for any good cube condition \eqref{equation: Assumption for extension on good cubes} is satisfied by our choice of $f$);
\item if $Q\in \mathscr{C}_{\eps}$ is a bad cube, then we let 
    \begin{align*}
        \tilde V_{\eps}:=\frac{1}{2^{n-1}}f\bigg(\frac{\eps}{2}\frac{x-c_Q}{\lVert x-c_Q\rVert_{\infty}}+c_Q\bigg)\frac{x-c_Q}{\lVert x-c_Q\rVert_{\infty}^n}, \qquad\forall\,x\in Q,
    \end{align*}
    with $f:=V_{\eps}\big|_{\partial Q}\cdot\nu_{\partial Q}\in L^{\infty}(\partial Q)$.
\end{enumerate}
We recall that no bad cubes will appear in the cubic decomposition in case $p\in\big[n/(n-1),+\infty\big)$ (see Remark \ref{Remark: no bad cubes if p is big enough}).

\textbf{Claim 1}. We claim that
\begin{align*}
    \operatorname{div}(\tilde V_{\eps})=\sum_{Q\in\mathscr{C}_{\eps}^b}d_Q\delta_{c_Q} \qquad\mbox{ distributionally on } \Omega_{\eps},
\end{align*}
where 
\begin{align*}
    d_{Q}:=\int_{\partial Q}V_{\eps}\cdot\nu_{\partial Q}\, d\H^{n-1}=\int_{\partial Q}V\cdot\nu_{\partial Q}\, d\H^{n-1}\in\z\smallsetminus\{0\},\qquad\forall\,Q\in\mathscr{C}_{\eps}^b.
\end{align*}
Indeed, pick any $\varphi\in C_c^{\infty}(\Omega_{\eps})$. Let $Q\in\mathscr{C}_{\eps}$ be a good cube. By the properties of the extension given by Lemma \ref{appendix: extension on the good cubes} and the divergence theorem we have 
\begin{align*}
    \int_{Q}\tilde V_{\eps}\cdot\nabla\varphi\, d\L^n&=\int_{\partial Q}\big(V_{\eps}-(V)_Q\big)\cdot\nu_{\partial Q}\varphi\, d\H^{n-1}+\int_{\partial Q}\big((V)_Q\cdot \nu_{\partial Q}\big)\varphi\, d\H^{n-1}\\
    &=\int_{\partial Q}(V_{\eps}\cdot\nu_{\partial Q})\varphi\, d\H^{n-1}.
\end{align*}
On the other hand, let $Q\in \mathscr{C}_{\eps}$ be a bad cube. By \eqref{divergence bad cubes}, we have 
\begin{align*}
    \int_{Q}\tilde V_{\eps}\cdot\nabla\varphi\, d\L^n=\int_{\partial Q}(V_{\eps}\cdot\nu_{\partial Q})\varphi\, d\H^{n-1}-d_Q\langle\delta_{c_Q},\varphi\rangle.
\end{align*}
Hence we conclude that
\begin{align*}
    \int_{\Omega_{\eps}}\tilde V_{\eps}\cdot\nabla\varphi\, d\L^n&=\sum_{Q\in\mathscr{C}_{\eps}}\int_{\partial Q}(V_{\eps}\cdot\nu_{\partial Q})\varphi\, d\H^{n-1}-\sum_{Q\in\mathscr{C}_{\eps}^b}d_Q\langle\delta_{c_Q},\varphi\rangle\\
    &=-\sum_{Q\in\mathscr{C}_{\eps}^b}d_Q\langle\delta_{c_Q},\varphi\rangle.
\end{align*}
The claim follows. 

\textbf{Claim 2}. We claim that $\lVert\tilde V_{\eps}-V\rVert_{L^p(\Omega_{\eps},\mu)}\rightarrow 0$ as $\eps\rightarrow 0^+$ in $E_V$.\\
Recall estimate \eqref{estimate on f} and notice that
\begin{align*}
    \lVert \tilde V_{\eps}-V\rVert_{L^p(\Omega_{\eps},\mu)}^p\le(1+C)(A_{\eps}+B_{\eps}),
\end{align*}
with
\begin{align*}
    A_{\eps}&:=\sum_{Q\in\mathscr{C}_{\eps}^g}\int_{Q}\lvert\tilde V_{\eps}-V\rvert^p f(c_Q)\, d\L^n,\\
    B_{\eps}&:=\sum_{Q\in\mathscr{C}_{\eps}^b}\int_{Q}\lvert\tilde V_{\eps}-V\rvert^p f(c_Q)\, d\L^n.
\end{align*}
By triangle inequality and by the estimate in Lemma \ref{appendix: extension on the good cubes}, we have that
\begin{align*}
    A_{\eps}&\le 2^{p-1}\Bigg(\sum_{Q\in\mathscr{C}_{\eps}^g}\int_{Q}\lvert\tilde V_{\eps}-(V)_Q\rvert^p f(c_Q)\, d\L^n+\sum_{Q\in\mathscr{C}_{\eps}^g}\int_{Q}\lvert V-(V)_Q\rvert^p f(c_Q)\, d\L^n\Bigg)\\
    &\le2^{p-1}\Bigg(\sum_{Q\in\mathscr{C}_{\eps}^g}\int_{Q}\lvert W_{\eps}\rvert^p f(c_Q)\, d\L^n+\sum_{Q\in\mathscr{C}_{\eps}^g}\int_{Q}\lvert V-(V)_Q\rvert^p f(c_Q)\, d\L^n\Bigg)\\
    &\le 2^{p-1}\Bigg(\eps C_p\sum_{Q\in\mathscr{C}_{\eps}^g}\int_{\partial Q}\lvert V_{\eps}-(V)_Q\rvert^p\, f(c_Q)d\H^{n-1}+\sum_{Q\in\mathscr{C}_{\eps}^g}\int_{Q}\lvert V-(V)_Q\rvert^pf(c_Q)\, d\L^n\Bigg),
\end{align*}
where $C_p:=C(p,Q)$ (see Remark \ref{appendix scaling argument}). Again by triangle inequality and because of our choice of $V_{\eps}$, we have
\begin{align*}
    \eps\sum_{Q\in\mathscr{C}_{\eps}^g}\int_{\partial Q}\lvert V_{\eps}-(V)_Q\rvert^pf(c_Q)\, d\H^{n-1}&\le 2^{p-1}\Bigg(\eps\sum_{Q\in\mathscr{C}_{\eps}^g}\int_{\partial Q}\lvert V_{\eps}-V\rvert^pf(c_Q)\, d\H^{n-1}\\
    &\quad+\eps\sum_{Q\in\mathscr{C}_{\eps}^g}\int_{\partial Q}\lvert V-(V)_Q\rvert^pf(c_Q)\, d\H^{n-1}\Bigg)\\
    &\le 2^{p-1}\Bigg(2n\eps^2+\eps\sum_{Q\in\mathscr{C}_{\eps}^g}\int_{\partial Q}\lvert V-(V)_Q\rvert^pf(c_Q)\, d\H^{n-1}\Bigg).
\end{align*}
Thus by Lemma \ref{Lemma: cubic decomposition} it follows that 
\begin{align*}
    \eps\sum_{Q\in\mathscr{C}_{\eps}}\int_{\partial Q}\lvert V_{\eps}-(V)_Q\rvert^pf(c_Q)\, d\H^{n-1}\rightarrow 0\qquad\mbox{ as } \eps\rightarrow 0^+ \mbox{ in } E_V.
\end{align*}
Moreover, by \eqref{estimate on f}
we have
\begin{align*}
    \sum_{Q\in\mathscr{C}_{\eps}^g}\int_{Q}\lvert V-(V)_Q\rvert^pf(c_Q)\, d\L^n
    &\leq 2^p(1+C)\sum_{Q\in\mathscr{C}_{\eps}^g}\fint_{Q_\varepsilon(0)}\int_Q\lvert V(x+y)-V(y)\rvert^pf(y)\, d\L^n\\ 
    &\leq 2^p(1+C)\fint_{Q_\varepsilon(0)}\lVert V(x+\cdot\,)-V\rVert_{L^p(\Omega_\varepsilon,\mu)}^p\to 0
\end{align*}
as $\eps\rightarrow 0^+$ in $E_V$. 
Hence, $A_{\eps}\rightarrow 0$ as $\eps\rightarrow 0^+$ in $E_V$. 

On the other hand, by \eqref{estimate bad cubes} we have
\begin{align*}
    B_{\eps}&\le2^{p-1}\Bigg(\sum_{Q\in\mathscr{C}_{\eps}^b}\int_{Q}\lvert\tilde V_{\eps}\rvert^pf(c_Q)\, d\L^n+\sum_{Q\in\mathscr{C}_{\eps}^b}\int_{Q}\lvert V\rvert^pf(c_Q)\, d\L^n\Bigg)\\
    &\le2^{p-1}\Bigg(C\eps\sum_{Q\in\mathscr{C}_{\eps}^b}\int_{\partial Q}\lvert V_{\eps}\rvert^pf(c_Q)\, d\H^{n-1}+\int_{\Omega_{\eps}^b}\lvert V\rvert^pf(c_Q)\, d\L^n\Bigg).
\end{align*}
We notice that
\begin{align*}
    \eps\sum_{Q\in\mathscr{C}_{\eps}^b}\int_{\partial Q}\lvert V_{\eps}\rvert^pf(c_Q)\, d\H^{n-1}&\le 2^{p-1}\Bigg(\eps\sum_{Q\in\mathscr{C}_{\eps}^b}\int_{\partial Q}\lvert V_{\eps}-V\rvert^pf(c_Q)\, d\H^{n-1}\\
    &\quad+\eps\sum_{Q\in\mathscr{C}_{\eps}^b}\int_{\partial Q}\lvert V\rvert^pf(c_Q)\, d\H^{n-1}\Bigg)\\
    &\le 4^{p-1}\Bigg(2n\eps^2+\eps\sum_{Q\in\mathscr{C}_{\eps}^b}\int_{\partial Q}\lvert V-(V)_Q\rvert^pf(c_Q)\, d\H^{n-1}\\
    &\quad+\eps\sum_{Q\in\mathscr{C}_{\eps}^b}\int_{\partial Q}\lvert (V)_Q\rvert^pf(c_Q)\, d\H^{n-1}\Bigg).
\end{align*}
Moreover, by \eqref{estimate on f} we have
\begin{align*}
    \eps\sum_{Q\in\mathscr{C}_{\eps}^b}\int_{\partial Q}\lvert (V)_Q\rvert^pf(c_Q)\, d\H^{n-1}&\le \eps\sum_{Q\in\mathscr{C}_{\eps}^b}\int_{\partial Q}\bigg(\fint_{Q}\lvert V\rvert^pf(c_Q)\, d\L^n\bigg)\, d\H^{n-1}\\
    &\le\sum_{Q\in\mathscr{C}_{\eps}^b}\eps\bigg(\fint_{Q}\lvert V\rvert^pf(c_Q)\, d\L^n\bigg)\int_{\partial Q}\, d\H^{n-1}\\
    &\le 2n\sum_{Q\in\mathscr{C}_{\eps}^b}\int_{Q}\lvert V\rvert^pf(c_Q)\, d\L^n=2n\int_{\Omega_{\eps}^b}\lvert V\rvert^pf(c_Q)\, d\L^n\\
    &\le (1+C)2n\int_{\Omega_{\eps}^b}\lvert V\rvert^pf\, d\L^n.
\end{align*}
Thus, we have obtained
\begin{align*}
    B_\varepsilon&\le C\Bigg(\eps^2+\eps\sum_{Q\in\mathscr{C}_{\eps}^b}\int_{\partial Q}\lvert V-(V)_Q\rvert^pf(c_Q)\, d\H^{n-1}+\int_{\Omega_{\eps}^b}\lvert V\rvert^p f\, d\L^n\Bigg),
\end{align*}
for some constant $C>0$ which does not depend on $\eps\in E_V$. By Lemma \ref{Lemma: cubic decomposition} and Remark \ref{Remark: no bad cubes if p is big enough} we get that $B_{\eps}\rightarrow 0$ as $\eps\rightarrow 0^+$ in $E_V$. Hence, the claim follows. 

Next we show that by rescaling $V_\varepsilon$ we obtain a vector field with similar properties defined on the whole $Q_1(0)$.
Let 
\begin{align*}
    \alpha_{\eps}:=\sup\{\alpha\in[1/2,1) \mbox{ s.t. } Q_1(0)\subset\alpha^{-1}\Omega_{\eps}\}, \qquad\forall\,\eps\in E_V. 
\end{align*}
Notice that $\alpha_{\eps}\rightarrow 1^-$ as $\eps\rightarrow 0^+$ in $E_V$. Define the vector field $\overline{V}_{\eps}:=\alpha_{\eps}^{n-1}\tilde V_{\eps}(\alpha_{\eps}\,\cdot):Q_1(0)\rightarrow\r^n$. It's straightforward that $\overline{V}_{\eps}\in L^p(Q_1(0),\mu)$ in case $p>1$ and $\overline{V}_{\eps}\in L^s(Q_1(0),\mu)$ for some $s>1$ in case $p=1$, for every given $\eps\in E_V$. A direct computation also shows that the distributional divergence of $\overline{V}_{\eps}$ on $Q_1(0)$ is given by
\begin{align*}
    \operatorname{div}(\overline{V}_{\eps})=\sum_{Q\in\mathscr{C}_{\eps}^b}\overline{d}_Q\delta_{\alpha_{\eps}^{-1}c_Q},
\end{align*}
with 
\begin{align*}
    \overline{d}_{Q'}=\begin{cases}d_{Q'} &\mbox{ if } \alpha_{\eps}^{-1}c_{Q'}\in Q_1(0),\\0 &\mbox{ otherwise.}\end{cases}
\end{align*}

We claim that $\overline{V}_{\eps}\rightarrow V$ in $L^p(Q_1(0),\mu)$. Indeed, we have 
\begin{align*}
    \int_{Q_1(0)}\lvert\overline{V}_{\eps}-V\rvert^pf\, d\L^n&=\int_{Q_1(0)}\lvert\alpha_{\eps}^{n-1}\tilde V_{\eps}(\alpha_{\eps}x)-V(x)\rvert^p\, f(x)d\L^n(x)\\ 
    &=\alpha_{\eps}^{p(n-1)-n}\int_{\alpha_{\eps}Q_1(0)}\lvert\tilde V_{\eps}(y)-\alpha_{\eps}^{-(n-1)}V(\alpha_{\eps}^{-1}y)\rvert^p f(\alpha_\varepsilon^{-1}y)\, d\L^n(y)\\
    &=\alpha_{\eps}^{p(n-1)-n}\bigg(\int_{\Omega_{\eps}}\lvert\tilde V_{\eps}(y)-\alpha_{\eps}^{-(n-1)}V(\alpha_{\eps}^{-1}y)\rvert^p f(y)\, d\L^n(y)\\
    &\quad+\int_{\Omega_{\eps}}\lvert\tilde V_{\eps}(y)-\alpha_{\eps}^{-(n-1)}V(\alpha_{\eps}^{-1}y)\rvert^p \frac{f(\alpha_\varepsilon^{-1}y)-f(y)}{f(y)}f(y)\, d\L^n(y)\bigg)\\
    &\le C_{n,p}\bigg(\int_{\Omega_{\eps}}\lvert\tilde V_{\eps}-V\rvert^pf\, d\L^n+\int_{\Omega_{\eps}}\lvert V-P_{\alpha_{\eps}^{-1}}V\rvert^pf\, d\L^n\bigg),
\end{align*}
(see Lemma \ref{lemma: continuity of the dilation operator} for the definition of $P_{\alpha_{\eps}^{-1}}V$). By Lemma \ref{lemma: continuity of the dilation operator} and since $\tilde V_{\eps}\rightarrow V$ in $L^p(\Omega_{\eps},\mu)$, our claim follows. 

Thus, we have built a vector field $\overline{V}_{\eps}$ such that:
\begin{enumerate}
    \item $\overline{V}_{\eps}\in L^p(Q_1(0),\mu)$ and $\overline{V}_{\eps}\in L^s(Q_1(0),\L^n)$ for $s=p$ if $p>1$ and $s>1$ if $p=1$\footnote{In fact even when $\mu$ is different from $\L^n$, $\tilde{V}_\varepsilon$ is constructed through Lemmata \ref{appendix: extension on the good cubes} and \ref{Lemma: extension on the bad cubes} as extension of a smooth boundary datum, thus $\tilde{V}_\varepsilon$ lies in $L^{r}(\Omega_{\varepsilon})$ for any $r\in [1,\frac{n}{n-1})$ if $p<n$ and in $L^r(\Omega_\varepsilon)$ for any $r\in [1,\infty)$ if $p\geq \frac{n}{n-1}$. It follows that $\overline V_\varepsilon \in L^{r}(\Omega_{\varepsilon})$ for any $r\in [1,\frac{n}{n-1})$ if $p<n$ and $\overline V_\varepsilon \in L^r(\Omega_\varepsilon)$ for any $r\in [1,\infty)$ if $p\geq \frac{n}{n-1}$.};
    \item the distributional divergence of $\overline{V}_{\eps}$ on $Q_1(0)$ is given by a finite sum of delta distributions supported on a finite set $X_{\eps}\subset Q_1(0)$ with integer weights $\{d_x \mbox{ s.t. } x\in X_{\eps}\}$;
    \item $\lVert\overline{V}_{\eps}-V\rVert_{L^p(Q_1(0),\mu)}\rightarrow 0$ as $\eps\rightarrow 0^+$ in $E_V$.
\end{enumerate}

Now we are ready to reach the conclusions 1 and 2 of Theorem \ref{Theorem: strong approximation for vector fields, now really for vector fields}. 

\begin{enumerate}
    \item If $q\in [0,1]$ and $p\in\big[1,\frac{n}{n-1}\big)$ we possibly have $X_{\eps}\neq\emptyset$, since bad cubes can appear in the cubic decompositions. Since $\overline{V}_{\eps}$ always belongs to $L^s(Q_1(0))$ for some $s>1$ (with $s=p$ if $p$ itself is already greater than 1),
    we can Hodge-decompose $\overline{V}_{\eps}^\flat$ as $\overline{V}_{\eps}^\flat=d\varphi+d^\ast A$ 
    for some $A\in \Omega_{W^{1,s}}^2(Q_1(0))$ and some $\varphi\in W^{1,s}(Q_1(0))$. 
    Applying $d^\ast$ to the previous decomposition we obtain
    \begin{align*}
         \Delta\varphi=d^\ast (\overline{V}_\varepsilon^\flat)=\operatorname{div}(\overline{V}_{\eps})=\sum_{x\in X_{\eps}}d_x\delta_x.
    \end{align*}
    By standard elliptic regularity, $\varphi\in C^{\infty}(Q_1(0)\smallsetminus X_{\eps})$. Choose $A_{\eps}\in \Omega^2(Q_1(0))$ such that $\lVert A_{\eps}-A\rVert_{\Omega^2_{W^{1,s}}(Q_1(0))}\le\eps$. Then 
    $\lVert d^\ast A_\varepsilon-d^\ast A\rVert_{\Omega^1_{L^p(\mu)}(Q_1(0))}\leq \varepsilon$.
    Let $\upsilon_{\eps}:=d\varphi+d^\ast(A_{\eps})$ and let $U_\varepsilon:=\upsilon_\eps^\#$. Then $U_{\eps}\in L_R^p(Q_1(0),\mu)$ for every $\eps\in E_V$ and $U_{\eps}\rightarrow V$ in $L^p(Q_1(0),\mu)$. 
    \item If $q\in (-\infty,0]$ and $p\in\big[\frac{n}{n-1},+\infty\big)$ no bad cubes are allowed in the cubic decomposition, thus $(\overline{V}_{\eps})_{\eps\in E_V}$ is a sequence of divergence-free vector fields converging to $V$ in $L^p(Q_1(0),\mu)$ as $\varepsilon\to 0^+$ in $E_V$. Hence $V$ itself is divergence-free.
\end{enumerate}
\end{proof}
\begin{rem}
\label{Remark: F_k in Lq}
Notice that if $p=1$, $V_k\in L^s(Q^n_1(0))$ for any $k\in \mathbb{N}$ and for any $s\in\big[1,\frac{n}{n-1}\big)$.
\end{rem}
\begin{rem}
    Observe that this proof can be used to show that the analogous approximation result holds if we assume that $V$ satisfies the first three conditions of Definition \ref{Definition: integer valued fluxes on unit cube} and in addition we require that for every $\rho\in R_{F, x_0}$ we have that
    \begin{align*}
        \int_{\partial Q_\rho(x_0)}i^\ast_{\partial Q_\rho(x_0)}F\in S
    \end{align*}
    for a set $S\subset\mathbb{R}$ such that $0\in S$ and $0$ is an isolated point in $S$. In this case the vector field $V$ can be approximated in $L^p$ by a sequence of vector fields $(V_n)_{n\in \mathbb{N}}$ smooth outside a finite set of points and such that for any $n\in \mathbb{N}$, $\operatorname{div}(V_n)$ is a finite sum of deltas with coefficients in $S$.
\end{rem}
\begin{lem}\label{Lemma: caso n=1}
Let $I\subset\r$ be a connected interval and $p\in[1,+\infty)$. Then
\begin{align*}
    \overline{L^p_R(I)}^{L^p}= L^p(I,\mathbb{Z})+ \mathbb{R}=L^p_\mathbb{Z}(I).
\end{align*}
\begin{proof}
    We start by showing the first equality.
    Notice that
    \begin{align*}
        L^p_{R}(I)=\left\{V=c+\sum_{j\in J}a_j\rchi_{I_j} \mbox{ : }(I_j)_{j\in J} \text{ is a finite partition of }I, a_j\in \mathbb{Z}\,\forall\, j\in J, \, c\in [0,1)\right\}.
    \end{align*}
    In other words, $L^p_R(I)$ consists of all integer-valued step functions and their translations by a constant.\\
    First we show the inclusion $"\supset"$: let $f=g+a$ with $g\in L^p(I,\mathbb{Z})$ and $a\in [0,1)$ and let $\varepsilon>0$. For any $k\in \mathbb{N}$ set $g_k:=\mathbbm{1}_{\lvert g\rvert\leq k}g$. Then there exists $K\in \mathbb{N}$ such that $\lVert g_K-g\rVert_{L^p(I)}<\frac{\varepsilon}{2}$.\\
    Now for any $j\in\{-K,..., K\}$ $g_K^{-1}(j)=g^{-1}(j)$ is a measurable set, therefore there exists a finite collection of disjoint intervals $(I_i^j)_{i\in J^j}$ such that $\L\left(\bigcup_{i\in J^j}I_i^j\triangle g^{-1}(j)\right)\leq \frac{\varepsilon}{2(2K+1)^2}$.\\
    For any $j\in\{-K,..., K\}$ set $A_j:=\bigcup_{i\in J^j}I_i^j\smallsetminus\left(\bigcup_{j'\neq j}\bigcup_{i\in J^{j'}}I_i^{j'}\right)$. Then $A_j$ is a finite union of intervals and
    \begin{align*}
        \L(A_j\triangle g^{-1}(j))\leq \L(A_j\smallsetminus g^{-1}(j))+\sum_{j'\neq j}\L\left(\bigcup_{i\in J^j}I_i^{j'}\cap g^{-1}(j)\right)\leq \frac{\varepsilon}{2(2K+1)}.
    \end{align*}
    Set \begin{align*}
        \tilde{g}_K=\begin{cases}
            j& \text{if }x\in A_j\\
            0&\text{otherwise}
        \end{cases}
    \end{align*}
    Then by construction $\tilde{g}_K\in L^p_R(I)$ and $\lVert \tilde{g}_K-g\rVert_{L^p(I)}\leq \varepsilon$.
    We conclude that any $g\in L^p(I, \mathbb{Z})$ lies in the closure of $L^p_R(I)$. This shows $"\supset"$.
    As $L^p(I,\mathbb{Z})+ \mathbb{R}=L^p(I,\mathbb{Z})+ [0,1)$ is closed in $L^p(I)$, the inclusion "$\subset$" holds as well.\\
    Next we show the second equality.
    Let's start with "$\subset$". Let $g\in L^p(I,\mathbb{Z})$, $a\in \mathbb{R}$ and $f=g+a$. Let $x_0\in I$, then for a.e. $r\in (0, \operatorname{dist}(x_0, \partial I))$ we have $f(x_0+r)-f(x_0-r)=g(x_0+r)-g(x_0-r)\in \mathbb{Z}$.
    Let $\tilde{R}_{F,x_0}$ denote the set of all such $r$. Set $R_{R,x_0}$ to be the intersection of $\tilde{R}_{F,x_0}$ with the set of Lebesgue points of $f$. Then $f$ satisfies Definition \ref{Definition: integer valued fluxes on unit cube} and thus $f\in L^p_\mathbb{Z}(I)$.\\
    To show "$\supset$" let $f\in L^p_\mathbb{Z}(I)$. Set $F: I\to\mathbb{S}^1,$ $x\mapsto e^{i2\pi f(x)}$. Then $F$ is a measurable bounded function. Notice that for any $x_0\in I$, for a.e. $r\in (0,\operatorname{dist}(x_0,\partial I))$ there holds $F(x_0-r)=F(x_0+r)$. This implies that $F$ is constant: this can be seen for instance approximating $F$ by smooth functions with the same symmetry properties away from $\partial I$ (convolving with a symmetric mollifier with small support), which then have to be constant. Choose $a\in \mathbb{R}$ such that $F\equiv e^{i2\pi a}$, then $f-a\in L^p(I,\mathbb{Z})$. This completes the proof.
\end{proof}
\end{lem}
\begin{rem}
\label{Remark: Closure of L_R and closure in case n=1}
    From Theorem \ref{Theorem: strong approximation for vector fields, now really for vector fields} it follows immediately that
    \begin{align*}
        \overline{L^p_R(Q_1(0))}^{L^p}=L^p_\mathbb{Z}(Q_1(0)).
    \end{align*}
    To see this it is enough to check that $L^p_\mathbb{Z}(Q_1(0))$ is closed in $L^p$, which can be shown by simple application of the coarea formula. 

    Notice that in case $p\in[n/(n-1),+\infty)$ we can approximate $V$ strongly in $L^p$ with smooth and divergence free vector fields. This is a straightforward consequence of Hodge decomposition. 
\end{rem}
\subsection{A characterization of $\Omega_{p,\mathbb{Z}}^{n-1}$}
First of all we apply the Strong approximation Theorem to obtain a useful characterization of the class $F\in\Omega_{p,\z}^{n-1}(Q_1(0))$.
\begin{thm}\label{Theorem: characterization of the integer valued fluxes class}
	Let $n\in \mathbb{N}_{>0}$, $p\in [1,+\infty)$. Let $F\in\Omega_p^{n-1}(Q^n_1(0))$. Then, the following are equivalent:
	\begin{enumerate}
		\item there exists $L\in\R_1(Q_1(0))$ such that $\partial L=*dF$ in $(W^{1,\infty}_0(Q_1(0)))^\ast$ and 
		\begin{align*}
			\m(L)=\sup_{\substack{\varphi\in\mathcal{D}(Q_1(0)), \\ ||d\varphi||_{L^{\infty}}\le 1}}\int_{Q_1(0)}F\wedge d\varphi.
		\end{align*}
		\item for every Lipschitz function $f:\overline{Q_1(0)}\rightarrow [a,b]\subset\r$ such that $f\vert_{\partial Q_1(0)}\equiv b$, we have 
		\begin{align*}
		\int_{f^{-1}(t)}i_{f^{-1}(t)}^*F\in\z, \qquad\mbox{ for $\L^1$-a.e. } t\in [a,b];
		\end{align*}
		\item $F\in\Omega_{p,\z}^{n-1}(Q_1(0))$.
	\end{enumerate}
\begin{proof}
We just need to show that $1\Rightarrow 2$, $2\Rightarrow 3$ and $3\Rightarrow 1$. We prove these implications separately.

$1\Rightarrow 2$. Assume 1. Let $L\in\R_1(Q_1(0))$ be given by 
\begin{align*}
    \langle L,\omega\rangle=\int_{\Gamma}\theta\langle\omega,\vec L\rangle\, d\H^1, \qquad\forall\,\omega\in\D^1(Q_1(0)),
\end{align*}
where $\Gamma\subset Q_1(0)$ is a locally $1$-rectifiable set, $\vec L$ is a Borel measurable unitary vector field on $\Gamma$ and $\theta\in L^1(\Gamma,\H^1)$ is a $\z$-valued function. Pick any Lipschitz function $f:Q_1(0)\rightarrow\r$ and let $\varphi\in C_c^{\infty}((-\infty,b))$ be such that $\int_{\r}\varphi\,d\L^1=0$. By the coarea formula we have 
\begin{align*}
    \int_{Q_1(0)}F\wedge f^*(\varphi\vol_{\r})=\int_{-\infty}^{+\infty}\varphi(t)\bigg(\int_{f^{-1}(t)}i_{f^{-1}(t)}^*F\bigg)\, dt.
\end{align*}
At the same time, by the coarea formula for countably $1$-rectifiable sets, we have 
\begin{align*}
    \langle L, f^*(\varphi\vol_{\r})\rangle=\int_{\Gamma}\theta\langle f^*(\varphi\vol_{\r}),\vec L\rangle\, d\H^1=\int_{-\infty}^{+\infty}\varphi(t)\bigg(\int_{\Gamma\cap f^{-1}(t)}\tilde\theta\bigg)\, dt,
\end{align*}
where $\tilde\theta:\Gamma\rightarrow\z$ is given by $\tilde\theta:=\text{sgn}(\langle f^*\vol_{\r},\vec L\rangle)\theta$. Let $\Phi\in C_c^{\infty}((-\infty,b))$ satisfy $d\Phi=\varphi\vol_{\r}$. Notice that since $f$ is proper $f^\ast\Phi\in W_0^{1,\infty}(Q_1(0))$. Then, by hypothesis, we have
\begin{align*}
    \int_{Q_1(0)}F\wedge f^*(\varphi\vol_{\r})&=\int_{Q_1(0)}F\wedge d f^\ast\Phi=\langle*dF,f^*\Phi\rangle=\langle\partial L,f^*\Phi\rangle\\
    &=\langle L,d(f^*\Phi)\rangle=\langle L,f^*(d\Phi)\rangle=\langle L,f^*(\varphi\vol_{\r})\rangle.
\end{align*}
Therefore
\begin{align*}
    \int_{-\infty}^{\infty}\varphi(t)\bigg(\int_{f^{-1}(t)}i_{f^{-1}(t)}^*F-\int_{\Gamma\cap f^{-1}(t)}\tilde\theta\bigg)\, dt=0, \qquad\forall\,\varphi\in C_c^{\infty}((-\infty,b)) \mbox{ s.t. } \int_{\r}\varphi=0.
\end{align*} 
We conclude that
\begin{align*}
    \int_{f^{-1}(t)}i_{f^{-1}(t)}^*F-\int_{\Gamma\cap f^{-1}(t)}\tilde\theta=c, \qquad\mbox{ for $\L^1$-a.e. } t\in[a,b],
\end{align*}
for some constant $c\in\r$. We claim that $c=0$. In fact let $m\in\n\smallsetminus\{0\}$.
Integrating both sides on $(-m,m)$ we get 
\begin{align}
\label{equation normalization}
    \int_{\{\lvert f\rvert< m\}}F\wedge f^*\vol_{\r}-\int_{\Gamma\cap\{\lvert f\rvert< m\}}\theta\langle f^*\vol_{\r},\vec L\rangle=2mc.
\end{align}
Since $f^*\vol_{\r}=df$, we have
\begin{align*}
    \int_{Q_1(0)}F\wedge f^*\vol_{\r}-\int_{\Gamma}\theta\langle f^*\vol_{\r},\vec L\rangle=\int_{Q_1(0)}F\wedge df-\langle L,df\rangle=\langle *dF-\partial L,f\rangle=0.
\end{align*}
Thus, by letting $m\rightarrow+\infty$ in \eqref{equation normalization}, we get that the left-hand-side converges to $0$ whilst the right-hand-side diverges to $+\infty$, unless $c=0$. Hence we conclude that $c=0$, i.e. 
\begin{align*}
    \int_{f^{-1}(t)}i_{f^{-1}(t)}^*F=\int_{\Gamma\cap f^{-1}(t)}\tilde\theta\in\z, \qquad\mbox{ for $\L^1$-a.e. } t\in [a,b],
\end{align*}
since $\Gamma\cap f^{-1}(t)$ consists of finitely many points for $\L^1$-a.e. $t\in \r$.\\
$2\Rightarrow 3$. Assume 2. Given $x_0\in Q_1(0)$, let $f_{x_0}:=\min\big\{\lVert\cdot-x_0\rVert_{\infty},\frac{r_{x_0}}{2}\big\}$. We claim that we can find $R_{F,x_0}$ as in Definition \ref{Definition: integer valued fluxes on unit cube}. Indeed, let $L\subset Q_1(0)$ be the set of the Lebesgue points of $F$. Let $r_{x_0}:=2\dist_\infty(x_0,\partial Q_1(0))$. Then, by the coarea formula, we have
\begin{align*}
    0=\L^n\big(Q_{r_{x_0}}(x_0)\smallsetminus L\big)=\frac{1}{2^n}\int_{0}^{r_{x_0}}\H^{n-1}\big((Q_1(0)\smallsetminus L)\cap\partial Q_{\rho}(x_0)\big)\,d\rho,
\end{align*}
which implies that there exists a set $E_{x_0}\subset(0,r_{x_0})$ such that
\begin{enumerate}
    \item $\L^1\big((0,r_{x_0})\smallsetminus E_{x_0}\big)=0$;
    \item for every $\rho\in E_{x_0}$, $\H^{n-1}$-a.e. $x\in\partial Q_{\rho}(x_0)$ is a Lebesgue point for $F$.
\end{enumerate} 
Hence, for every $\rho\in E_{x_0}$ it makes sense to consider the pointwise restriction of $F$ to $\partial Q_\rho(x_0)$. Notice that, by the coarea formula, we have 
\begin{align*}
    \int_{E_{x_0}}\bigg(\int_{\partial Q_{\rho}(x_0)}\big|F\big|^p\, d\H^{n-1}\bigg)\,d\rho=2^n\int_{Q_{r_{x_0}}(x_0)}|F|^p\,d\L^n<+\infty,
\end{align*}
which implies 
\begin{align*}
    \int_{\partial Q_{\rho}(x_0)}\big|F\big|^p\, d\H^{n-1}<+\infty, \qquad\mbox{ for $\L^1$-a.e. }\rho\in (0,E_{x_0}).
\end{align*}
Moreover, by Statement $2.$, we have 
\begin{align*}
    \int_{f_{x_0}^{-1}(\rho)}i_{f_{x_0}^{-1}(\rho)}^*F=\int_{\partial Q_{\rho}(x_0)}i_{\partial Q_{\rho}(x_0)}^*F\in\z, \qquad\mbox{ for $\L^1$-a.e. }\rho\in (0,E_{x_0}).
\end{align*}
Our claim follows immediately.

$3\Rightarrow 1$. Assume 3. By Theorem \ref{Theorem: strong approximation for vector fields}, we can find a sequence $\{F_k\}_{k\in\n}\subset\Omega_{p,R}^{n-1}(Q_1(0))$ such that $F_k\rightarrow F$ strongly in $L^p$. Since the estimate
\begin{align*}
	\lvert\langle T_{F_k}-T_{F},\omega\rangle\rvert\le C\lvert\lvert F_k-F\rvert\rvert_{L^p}
\end{align*}
holds for every $\omega\in\mathcal{D}^1(Q_1(0))$ such that $||\omega||_{L^{\infty}}\le 1$ and for every $k\in\n$, we conclude that
\begin{align}\label{useful equation bla}
	\sup_{\substack{\varphi\in W^{1,\infty}_0(Q_1(0)),\\\lVert d\varphi\rVert_{L^\infty}\leq 1}}\langle*dF_k-*dF,\varphi\rangle\le C\lvert\lvert F_k-F\rvert\rvert_{L^p}\rightarrow 0\quad \text{ as }k\to\infty.
\end{align}
Fix any $\varepsilon\in(0,1)$. By \eqref{useful equation bla}  we can find a subsequence $\big\{F_{k_j(\varepsilon)}\big\}_{j\in\n}\subset\Omega_{p,R}^{n-1}(Q_1(0))$ such that
\begin{align*}
	\lVert *dF_{k_j(\varepsilon)}-*dF_{k_{j+1}(\varepsilon)}\rVert_{(W_0^{1,\infty}(Q_1(0)))^\ast}\le\frac{\varepsilon}{2^j}, \qquad \mbox{ for every } j\in\n.
\end{align*}
For every $h\in\n$, let $L_h^{\varepsilon}$ be a minimal connection for the singular set of $F_{k_h(\varepsilon)}$ (the existence of such a minimal connection is proved in Proposition \ref{appendix proposition existence minimal connection finitely many singularities}. Analogously, for every $j\in\n$, let $L_{j,j+1}^{\varepsilon}$ be the minimal connection for the singular set of $F_{k_j(\varepsilon)}-F_{k_{j+1}(\varepsilon)}$.
Define the following sequence of integer $1$-currents on $Q_1(0)$:
\begin{align*}
	\tilde L_h^{\varepsilon}:=\begin{cases}
	L_0^{\varepsilon} & \mbox{ if } h=0,\\ \displaystyle{L_0^{\varepsilon}-\sum_{j=0}^{h-1}L_{j,j+1}^{\varepsilon}}& \mbox{ if } h>0,
	\end{cases} \qquad \mbox{ for every } h\in\n.
\end{align*}
Clearly
\begin{align*}
\partial\tilde L_h^{\varepsilon}=\partial L_0^{\varepsilon}-\sum_{j=0}^{h-1}\partial L_{j,j+1}^{\varepsilon}=\partial L_0^{\varepsilon}-\sum_{j=0}^{h-1}(\partial L_j^{\varepsilon}-\partial L_{j+1}^{\varepsilon})=\partial L_h^{\varepsilon}=*dF_{k_h(\varepsilon)}.
\end{align*}
Moreover, since $L_{j,j+1}^{\varepsilon}$ is a minimal connection, it holds that
\begin{align*}
	\mathbb{M}(L_{j, j+1}^\varepsilon)=\lVert *dF_{k_j(\varepsilon)}-*dF_{k_{j+1}(\varepsilon)}\rVert_{(W_0^{1,\infty}(Q_1(0)))^\ast}\le\frac{\varepsilon}{2^j}, \qquad \mbox{ for every } j\in\n.
\end{align*}
Thus,
\begin{align*}
	\mathbb{M}(\tilde L_{h+1}^{\varepsilon}-\tilde L_h^{\varepsilon})=\mathbb{M} (L_{h,h+1}^{\varepsilon})\le\frac{\varepsilon}{2^h}, \qquad \mbox{ for every } h\in\n,
\end{align*}
which amounts to saying that the sequence $\{\tilde L_h^{\varepsilon}\}_{h\in\n}$ is a Cauchy sequence in mass. Hence, by the closure of integer currents under mass convergence (see Lemma \ref{appendix Lemma closure integer currents}), there exists an integer $1$-current $\tilde L^{\varepsilon}\in\mathcal{R}_1(Q_1(0))$ such that 
\begin{align*}
	\mathbb{M} (\tilde L_h^{\varepsilon}-\tilde L^{\varepsilon})\rightarrow 0\quad\text{ as }h\to\infty,
\end{align*}
Notice that
\begin{align*}
	\partial \tilde L^{\varepsilon}=\lim_{h\rightarrow\infty}\partial \tilde L_h^{\varepsilon}=\lim_{h\rightarrow\infty}*dF_{k_h(\varepsilon)}=*dF\quad \text{ in }(W_0^{1,\infty}(Q_1(0)))^\ast.
	\end{align*}
The family of integer $1$-cycles $\{\tilde L^{\varepsilon}-\tilde L^{1/2}\}_{0<\varepsilon<1}\subset\mathcal{R}_1(Q_1(0))$ is uniformly bounded in mass. Indeed, first we notice that by (2.10) it holds that
\begin{align*}
    \m(L_{h}^{\eps})&=\lVert *dF_{k_h(\varepsilon)}\rVert_{(W_0^{1,\infty}(Q_1(0)))^\ast}\le C, \quad\forall\,h\in\mathbb{N},\,\forall\,\eps\in(0,1),
\end{align*} 
where $C>0$ is a constant independent on $h$ and $\eps$. Thus, we have
\begin{align*}
    \m(\tilde L_h^{\eps})&\le\m(L_0^{\eps})+\sum_{j=0}^{h-1}\m(L_{j,j+1}^{\eps})\le C+\sum_{j=0}^{+\infty}\frac{\eps}{2^j}\le C+\eps\le C+1\quad\forall\,h\in\mathbb{N},\,\forall\,\eps\in(0,1).
\end{align*}
Since $\tilde L_{h}^{\eps}\to\tilde L^{\eps}$ in mass as $h\to+\infty$ for every $\eps\in(0,1)$, we have
\begin{align*}
\m(\tilde L^{\eps}-\tilde L^{1/2})\le \m(\tilde L^{\eps})+\m(\tilde L^{1/2})\le 2(C+1).
\end{align*}

Hence, by standard compactness arguments for currents (see for instance \cite{krantz_parks-geometric_integration_theory}, Theorem 7.5.2), we can find a sequence $\varepsilon_k\rightarrow 0$ and an integer 1-cycle $\tilde L\in\mathcal{R}_1(Q_1(0))$ with finite mass such that $\tilde L^{\varepsilon_k}-\tilde L^{1/2}\rightarrow  \tilde L$ weakly in $\mathcal{D}_1(Q_1(0))$ as $k\rightarrow+\infty$. If we let $L:=\tilde L^{1/2}+\tilde L$ then we get $\tilde L^{\varepsilon_k}\rightarrow  L$ weakly in $\mathcal{D}_1(Q_1(0))$.
By construction, $L$ is again an integer $1$-current with finite mass such that $\partial L=\partial\tilde L^{1/2}=*dF$ in $(W^{1,\infty}(Q_1(0)))^\ast$. We claim that
\begin{align*}
	\m(L)=\inf_{\substack{T\in\mathcal{M}_1(Q_1(0)), \\ \partial T=*dF}}\m(T).
\end{align*}
By contradiction, assume that we can find $T\in\mathcal{M}_1(Q_1(0))$ such that $\partial T=*dF$ and 
\begin{align*}
	\m(T)<\m(L)\le\liminf_{k\rightarrow\infty}\m(\tilde L^{\varepsilon_k}),
\end{align*}
where the last inequality follows by weak convergence and lower semicontinuity of the mass. Then, we can find some $h\in\n$ such that
\begin{align*}
	\m(T)<\m(\tilde L^{\varepsilon_{h}})-2\varepsilon_{h}.
	\end{align*}
Moreover, since $\m(L_0^{\varepsilon}-\tilde L^{\varepsilon})\le 2\varepsilon$ for every $0<\varepsilon<1$, it holds that
\begin{align*}
	\m(L_0^{\varepsilon_h}-\tilde L^{\varepsilon_h})\le 2\varepsilon_{h}.
\end{align*}
We define $\tilde T:=T+L_0^{\varepsilon_h}-\tilde L^{\varepsilon_h}$ and we notice that $\partial\tilde T=*dF_{k_0(\varepsilon_h)}$. Moreover, by the minimality of $L_0^{\varepsilon_h}$, we conclude that
\begin{align*}
	\m(L_0^{\varepsilon_h})\le\m(\tilde T)\le\m(T)+\m(L_0^{\varepsilon_h}-\tilde L^{\varepsilon_h})<\m(L_0^{\varepsilon_h}),
\end{align*}
which is a contradiction. Thus, our claim follows.

Since $L\in\mathcal{R}_1(Q)$, we get that
\begin{align*}
	\m(L)=\inf_{\substack{T\in\mathcal{R}_1(Q_1(0)), \\ \partial T=*dF}}\m(T)=\inf_{\substack{T\in\mathcal{M}_1(Q_1(0)), \\ \partial T=*dF}}\m(T)
\end{align*}
and, by Lemma \ref{appendix lemma duality}, we have
\begin{align*}
	\inf_{\substack{T\in\mathcal{M}_1(Q_1(0)), \\ \partial T=*dF}}\m(T)=\sup_{\substack{\varphi\in\mathcal{D}(Q_1(0)), \\ \lvert\lvert d\varphi\rvert\rvert_{L^{\infty}}\le 1}}\int_MF\wedge d\varphi.
\end{align*}
Hence, 1. follows.
\end{proof}
\end{thm}
\begin{rem}\label{Remark: existence of connection is enough for strong approximability on cube}
We notice that in the proof of Theorem \ref{Theorem: characterization of the integer valued fluxes class} we have never used the minimality property of $L$ while showing that $1\Rightarrow 2$. Hence whenever $F\in\Omega_p^{n-1}(Q^n_1(0))$ admits a connection we have $F\in\Omega_{p,\z}^{n-1}(Q^n_1(0))$ in the sense of Definition \ref{Definition: integer valued fluxes on unit cube} and the conclusions of Theorem \ref{Theorem: strong approximation for vector fields} hold for $F$.
\end{rem}
By the previous remark we can deduce the following result from the proof of Theorem \ref{Theorem: characterization of the integer valued fluxes class}.
\begin{cor}
\label{Corollary: existence of connections implies existence of minimal connection}
Let $F\in \Omega_p^{n-1}(Q_1^n(0))$ and assume that there exists an integer $1$-current of finite mass $I\in \mathcal{R}_1(Q_1(0))$ such that $\partial I=\ast dF$. Then there exists an integer $1$-current $L\in \mathcal{R}_1(Q_1(0))$ of finite mass such that $\partial L=\ast dF$ and
\begin{align*}
    \m(L)=\inf_{\substack{T\in\mathcal{R}_1(Q_1(0))\\\partial T=\ast d F}}\m(T).
\end{align*}
\end{cor}
In other words, whenever there exists a connection for $F$, then there exists a minimal connection for $F$.
\subsection{The case of $\partial Q_1^{n+1}(0)$}
In order to extend the previous result to more general manifolds we introduce the following definition.
\begin{dfn}\label{definition: Lipschitz forms}
Let $M$ be a Lipschitz $m$-manifold embedded in $\mathbb{R}^n$. Let $p\in [1,\infty)$. Set
\begin{align*}
    \Omega_{p,R,\infty}^1(M):=\left\{\alpha\in \Omega_{p}^1(M)\cap\Omega_{L^\infty_{loc}}^1(M\smallsetminus S) \mbox{ : }\hspace{-1mm}\ast\hspace{-0.5mm}d\alpha=\sum_{i\in I}d_i\delta_{p_i}\right\},
\end{align*}
where $I$ is a finite index set, $d_i\in \mathbb{Z}$, $p_i\in S$ for any $i\in I$ and $F:=\{p_i\}_{i\in I}$.
\end{dfn}
The previous definition is motivated by the following observation: let $M$ be a Lipschitz $m$-manifold, $N$ a smooth $m$-manifold, $\varphi: M\to N$ a bi-Lipschitz map. Let $F\in \Omega_{p, R}^1(N)$. Then $\varphi^\ast F\in \Omega_{p,R,\infty}^\infty$ (see Lemma \ref{Lemma: bi-Lipschitz maps preserve approximability properties}).
\begin{cor}\label{Corollary: strong approximation on boundary of cube}
Let $n\in \mathbb{N}_{\geq 2}$. Assume that $F\in\Omega_{p,\z}^{n-2}(\partial Q^n_1(0))$ admits a connection. Then, the following facts hold:
\begin{enumerate}
    \item if $p\in\big[1,(n-1)/(n-2)\big)$, then there exists a sequence $\{ F_k\}_{k\in\n}\subset \Omega_{p,R,\infty}^{n-2}(\partial Q^n_1(0))$ such that $F_k\rightarrow F$ strongly in $L^p$;
    \item if $p\in\big[(n-1)/(n-2),+\infty\big)$, then $*dF=0$ distributionally on $\partial Q^n$. 
\end{enumerate}
\begin{proof}


Let $N=(0,...,0,\frac{1}{2})$ be the north pole in $\partial Q^n_1(0)\subset\r^m$ and let
\begin{align*}
    U:=\big\{(x_1,...,x_{n-1},x_n)\in\partial Q^n_1(0) \mbox{ s.t. } x_n=\frac{1}{2}\big\}
\end{align*}
be the upper face of $\partial Q^n_1(0)$. Let $q:=(n-1)-(n-2)p$. For every $x=(x_1,...,x_n)\in\r^n$, we let $x':=(x_1,...,x_{n-1})\in\r^{n-1}$. Define $\Phi:\partial Q^n_1(0)\smallsetminus N\subset\r^n\to Q^{n-1}_1(0)$ by
\begin{align*}
    \Phi(x):=\begin{cases}\displaystyle{\bigg(\frac{1}{2}-\frac{\sqrt{2}}{4}\lVert x'\rVert_{\infty}^\frac{1}{2}\bigg)\frac{x'}{\lVert x'\rVert_{\infty}}} & \mbox{ on }\,U\smallsetminus N,\vspace{2mm}\\\,g(x) & \mbox{ on }\,\partial Q^n_1(0)\smallsetminus U,\end{cases}
\end{align*} 
where the map $g:\overline{\partial Q_1^n(0)\smallsetminus U}\to \overline{Q_{1/2}^{n-1}(0)}$ is any bi-Lipschitz homeomorphism such that $g\equiv\bigg(\frac{1}{2}-\frac{\sqrt{2}}{4}\lVert x'\rVert_{\infty}^\frac{1}{2}\bigg)\frac{x'}{\lVert x'\rVert_{\infty}}$ on $\partial U$. Notice that $\Phi$ is an homeomorphism from $\partial Q^n_1(0)\smallsetminus N$ to $Q^{n-1}_1(0)$. We denote its inverse map by $\Psi$.\\
We have that $\Psi$ is Lipschitz on $Q^{n-1}_1(0)$ and $\Phi$ is Lipschitz on every compact set $K\subset\partial Q^n_1(0)\smallsetminus N$, since there exists $C>0$ such that
\begin{align*}
    \lvert d\Phi(x)\rvert\le\frac{C}{\lVert x-N\rVert_{\infty}^\frac{1}{2}}, \qquad\forall\, x\in\partial Q^n_1(0)\smallsetminus N,\\
     \lvert d\Psi(y)\rvert\le C\bigg(\frac{1}{2}-\lVert y\rVert_{\infty}\bigg), \qquad\forall\, y\in Q_1^{n-1}(0).
\end{align*}

Define $\tilde F:=\Psi^*F$ and fix any $\eps\in(0,1/4)$. Notice that
\begin{align*}
    \int_{Q^{n-1}_1(0)}\bigg(\frac{1}{2}-\lVert\,\cdot\,\rVert_{\infty}\bigg)^{q}\lvert\tilde F\rvert^p\, d\H^{n-1}&\le C\int_{Q^{n-1}_1(0)}\bigg(\frac{1}{2}-\lVert\,\cdot\,\rVert_{\infty}\bigg)^{n-1}\lvert F\circ\Psi\rvert^p\, d\H^{n-1}\\
    &\le C\int_{\partial Q^n_1(0)}\lvert F\lvert^p\, d\H^{n-1}<+\infty.
\end{align*}
Moreover if $I\in \mathcal{R}_1(\partial Q_1^n(0))$ and $\ast dF =\partial I$, then $\Phi_\ast I\in \mathcal{R}_1 (\partial Q_1^{n-1}(0))$ and $\ast d \tilde F=\partial \Phi_\ast I$ (see Lemma \ref{Lemma: bi-Lipschitz maps preserve approximability properties}).
This implies that $\tilde F\in\Omega_{p,\z}^{n-2}(Q^{n-1}_1(0),\mu)$ with $\mu:=(\frac{1}{2}-\lVert\cdot\rVert_{\infty})^q\, \L^{n-1}$ in the sense of Definition \ref{Definition: integer valued fluxes on unit cube}.\\
Let's consider first the case $p\in\big[1,(n-1)/(n-2)\big)$. Notice that in this case $q>0$.

By Theorem \ref{Theorem: strong approximation for vector fields, now really for vector fields} with $f:=(\frac{1}{2}-\lVert\cdot\rVert_{\infty})^q$ there exists a $(n-2)$-form $\tilde F_{\eps}\in\Omega_{p,R}^{n-2}(\Omega_{\eps,a_{\eps}})$ such that $\lVert \tilde F_\varepsilon-\tilde F\rVert_{L^p(S_{\varepsilon, a_\varepsilon})}\leq \varepsilon$ and
\begin{align}\label{Equation: convergence of approximation of partial Q}
    \lVert \tilde F_{\eps}-\tilde F\rVert^p_{L^p(\Omega_{\eps,a_{\eps}},\mu)}\rightarrow 0
\end{align}
as $\eps\rightarrow 0^+$ in $E_{\tilde F}$.

Define $F_{\eps}:=\Phi_{a_{\eps}}^*\tilde F_{\eps}$ on $\partial Q_1^n(0)\smallsetminus U_{\eps}$, with $\Phi_{a_{\eps}}:=\Phi+a_{\eps}$ and $U_{\eps}:=\Phi_{a_{\eps}}^{-1}(Q_1^{n-1}(0)\smallsetminus\Omega_{\eps,a_{\eps}})$.

Notice that
\begin{align*}
    \lVert F_{\eps}-F\rVert_{L^p(\partial Q^n_1(0)\smallsetminus U_{\eps})}^p&\le C\int_{\partial Q^n_1(0)\smallsetminus U_{\eps}}\frac{1}{\lVert\cdot - N\rVert_{\infty}^{(n-2)p}}\lvert\tilde F_{\eps}\circ\Phi_{a_{\eps}}-\tilde F\circ\Phi_{a_{\eps}}\rvert^p\, d\H^{n-1}\\
    &\quad+C\int_{\partial Q^n_1(0)\smallsetminus U_{\eps}}\frac{1}{\lVert\cdot - N\rVert_{\infty}^{(n-2)p}}\lvert\tilde F\circ\Phi_{a_{\eps}}-\tilde F\circ\Phi\rvert^p\, d\H^{n-1}\\
    &\le C\bigg(\int_{\Omega_{\eps,a_{\eps}}}\lvert\tilde F_{\eps}-\tilde F\rvert^p f\, d\H^{n-1}+\int_{Q_{1-\eps}^{n-1}(0)}\lvert\tilde F(\cdot -a_{\eps})-\tilde F\rvert^p f\, d\H^{n-1}\bigg).
\end{align*}
The first term tends to zero as $\varepsilon\to 0^+$ in $E_{\tilde F}$ by \eqref{Equation: convergence of approximation of partial Q}, while the second tends to zero as $\varepsilon\to 0^+$ by \eqref{Equation: Continuity wrt argument for weighted measure}.
Therefore we have
\begin{align*}
    \lVert F_{\eps}-F\rVert_{L^p(\partial Q_1^n(0)\smallsetminus U_{\eps})}^p\rightarrow 0
\end{align*}
as $\eps\rightarrow 0^+$ in $E_{\tilde F}$.

Now we notice that
\begin{align*}
    \int_{\partial U_{\eps}}i_{\partial U_{\eps}}^*F_{\eps}&=-\int_{\partial(\partial Q_1^n(0)\smallsetminus U_{\eps})}i_{\partial(\partial Q_1^n(0)\smallsetminus U_{\eps})}^*F_{\eps}\\
    &= \int_{\partial\Omega_{\eps,a_{\eps}}}i_{\partial\Omega_{\eps,a_{\eps}}}^*\tilde F_{\eps}=\int_{\partial\Omega_{\eps,a_{\eps}}}i_{\partial\Omega_{\eps,a_{\eps}}}^*\tilde F=:b_{\eps}\in\z.
\end{align*}
If $b_{\eps}=0$, then we use Lemma \ref{appendix: extension on the good cubes} to extend $F_{\eps}$ inside $U_{\eps}$. If $b_{\eps}\neq 0$, then we use Lemma \ref{Lemma: extension on the bad cubes} to extend $F_{\eps}$ inside $U_{\eps}$ (notice that $U_{\eps}$ is an $(n-1)$-cube of side-length $8\varepsilon^2$ contained in $U$ and centered at $N$, for $\eps$ sufficiently small). In both cases, the following estimate holds:
\begin{align*}
    \int_{U_{\eps}}\lvert F_{\eps}\rvert^p\, d\H^{n-1}\le C\eps^2\int_{\partial U_{\eps}}\lvert F_\varepsilon\rvert^p\, d\H^{n-2}\leq C\varepsilon^2\left(\int_{\partial U_\varepsilon}\lvert F-F_\varepsilon\rvert^p d\mathscr{H}^{n-2}+\int_{\partial U_\varepsilon}\lvert F\rvert^p d\mathscr{H}^{n-2}\right).
\end{align*}
Notice that the first term on the right hand side tends to zero as $\varepsilon\to 0^+$ in $E_{\tilde F}$, since $\lVert \tilde F_\varepsilon-\tilde F\rVert_{L^p(S_{\varepsilon, a_\varepsilon})}\leq \varepsilon$ for any $\varepsilon\in E_{\tilde{F}}$.
In order to control the second term, pick any $\delta>0$ sufficiently small and notice that, by coarea formula, we have
\begin{align*}
    \fint_0^{\delta}\eps^2\int_{\partial U_{\eps}}\lvert F\rvert^p\, d\H^{n-2}\, d\L^1(\eps)&\le\int_0^{\delta}\varepsilon\int_{Q^{n-1}_{8\varepsilon^2}(N)}\lvert F\rvert^p\, d\H^{n-2}\, d\L^1(\eps)\\
    &\leq \frac{1}{16}\int_0^{8\delta^2}\int_{\partial Q^{n-1}_{\zeta}(N)}\lvert F\rvert^p\, d\H^{n-2}\, d\L^1(\zeta)\\
    &\le \frac{2^{n-1}}{16}\int_{Q^{n-1}_{8\delta^2}(N)}\lvert F\rvert^p\, d\H^{n-1}\rightarrow 0^+
\end{align*}
as $\delta\rightarrow 0^+$. This implies that we can pick a sequence $\{\eps_j\}_{j\in\n}\subset E_{\tilde F}$ such that $\eps_j\rightarrow 0^+$ and 
\begin{align*}
    \eps_j^2\int_{\partial U_{\eps_j}}\lvert F\rvert^p\,d\H^{m-1}\rightarrow 0^+
\end{align*}
as $j\rightarrow+\infty$. Thus 
\begin{align*}
    \lVert F_{\eps_j}-F\rVert_{L^p(U_{\eps_j})}^p\le 2^{p-1}\Big(\lVert F_{\eps_j}\rVert_{L^p(U_{\eps_j})}^p+\lVert F\rVert_{L^p(U_{\eps_j})}^p\Big)\rightarrow 0
\end{align*}
as $j\rightarrow+\infty$, since $\H^{n-1}\big(U_{\eps_j}\big)\rightarrow 0^+$ as $\eps_j\rightarrow 0^+$. Hence, we conclude that
\begin{align*}
    \lVert F_{\eps_j}-F\rVert_{L^p(\partial Q^n_1(0))}^p\rightarrow 0
\end{align*}
as $j\rightarrow+\infty$. Moreover, by construction we have that $\ast d F_{\eps_j}$ is a finite sum of Dirac-deltas with integer coefficients, for any $j\in \mathbb{N}$. Thus arguing as in the final step of the proof of Theorem \ref{Theorem: strong approximation for vector fields, now really for vector fields} (i.e. by Hodge decomposition), for any $j\in \mathbb{N}$ we can find $\hat F_{\varepsilon_j}\in \Omega^1_{p,R,\infty}(\partial Q^n_1(0))$ such that $\lVert F_{\eps_j}-\hat F_{\eps_j}\rVert_{L^p(\partial Q_1^n(0))}<\varepsilon_j$.
The sequence $\{\hat F_j\}_{j\in \mathbb{N}}$ then has the desired properties.
This concludes the proof in the case $p\in\big[1,(n-1)/(n-2)\big)$.\\
If $p\in\big[(n-1)/(n-2),+\infty\big)$ notice that $q\leq 0$. Therefore by Remark \ref{Remark: no bad cubes if p is big enough} we may assume, up to passing to a subsequence, that no bad cube appears in the construction of $\tilde{F}_\varepsilon$. Hence repeating the first part of the proof as in the previous case, we get $b_{\varepsilon_j}=0$ and $\ast d F_{\varepsilon_j} =0$ on $\partial Q^n_1(0)$ for any $j\in\mathbb{N}$. Thus we obtain that $F$ can be approximated in $\Omega_p^{n-1}(\partial Q^n_1(0))$ by a sequence $(F_{\varepsilon_j})_{j\in \mathbb{N}}$ such that for any $j\in \mathbb{N}$ there holds $\ast d F_{\varepsilon_j}=0$, and this property passes to the limit.
\end{proof}
\end{cor}
\begin{rem}
Notice that if $p=1$, for any $k\in \mathbb{N}$ we have that $F_k\in \Omega_{q}^{n-1}(\partial Q^n_1(0))$ for some $q>1$ which does not depend on $k$ (compare with Remark \ref{Remark: F_k in Lq}).
\end{rem}
\begin{lem}
\label{Lemma: bi-Lipschitz maps preserve approximability properties}
Let $M,N\subset \mathbb{R}^n$ be Lipschitz $m$-manifolds in $\mathbb{R}^n$. Let $\varphi: M\to N$ be a bi-Lipschitz map.\\
Let $F\in \Omega_p^{m-1}(M)$ and assume that there exists a 1-rectifiable current $I\in \mathcal{R}_1(M)$ of finite mass such that $\ast d F=\partial I$ in $(W^{1,\infty}_0(M))^\ast$. Then  $\varphi_\ast F\in \Omega_p^{m-1}(N)$ and $\ast d(\varphi_\ast F)=\partial \varphi_\ast I$.\\
If $(F_k)_{k\in \mathbb{N}}$ is a sequence in $\Omega_{R,p,\infty}^{m-1}(M)$ and $F_k\to F$ in $\Omega_p^{m-1}(M)$ as $k\to\infty$, then $(\varphi_\ast F_k)_{k\in \mathbb{N}}$ is a sequence in $\Omega_{p,R,\infty}^{m-1}(N)$ and $\varphi_\ast F_k\to\varphi_\ast F$ in $\Omega_p^{m-1}(N)$ as $k\to\infty$.\\
If in addition we assume that $N$ is smooth and closed or a bounded simply connected Lipschitz domain and for any $k\in \mathbb{N}$ we have $F_k\in \Omega_q^{m-1}(N)$ for some $q>1$ (possibly dependent on $k$), then $\varphi_\ast F$ can be approximated in $\Omega_{p}^{m-1}(N)$ by $(m-1)$-forms in $\Omega_{p,R}^{m-1}(N)$.
\end{lem}
\begin{proof}
Assume that $\ast d F=\partial I$ holds in $(W^{1,\infty}_0(M))^\ast$. Then for any $f\in W^{1,\infty}_0(N)$ $\varphi^\ast f\in W^{1,\infty}_0(M)$ and thus
\begin{align*}
    \langle \ast d (\varphi_\ast F), f\rangle=\int_N\varphi_\ast F\wedge df=\int_M F\wedge d\varphi^\ast f=\langle \partial I, \varphi^\ast f\rangle=\langle I, \varphi^\ast df\rangle=\langle \partial \varphi_\ast I, f\rangle,
\end{align*}
therefore $\ast d(\varphi_\ast F)=\partial \varphi_\ast I$ in $(W^{1,\infty}_0(N))^\ast$.\\
Now assume that $(F_k)_{k\in \mathbb{N}}$ is a sequence in $\Omega_R^{m-1}(M)$ such that $F_k\to F$ in $\Omega_p^{m-1}(M)$ as $k\to\infty$. Then for any $k\in \mathbb{N}$ there exists $I_k\in \mathcal{R}_1(M)$ of finite mass and so that $\partial I_k$ supported in a finite subset of $M$ such that $\ast d F_k=\partial I_k$. As we saw above, $\ast d(\varphi_\ast F_k)=\partial \varphi_\ast I_k$, therefore $\varphi_\ast F_k\in \Omega_{p, R,\infty}^{m-1}(N)$. Moreover we have $\varphi_\ast F_k\to \varphi_\ast F$ in $\Omega_p^{m-1}(N)$ as $k\to\infty$.\\
Finally if $N$ is smooth and closed or a bounded Lipschitz domain and for any $k\in \mathbb{N}$ we have that $F_k\in \Omega_q^{m-1}(N)$ for some $q>1$, we can improve the approximating sequence $(\varphi^\ast F_k)_{k\in \mathbb{N}}$ as follows: for any $k\in \mathbb{N}$, let $\alpha_k\in \Omega_{W^{1,q}}^{m-2}(N)$, $\beta_k\in \Omega_{W^{1,q}}^{m}(N)$ and $h_k\in \Omega_h^{m-1}(N)$ (the space of harmonic $(m-1)$-forms on $N$) such that $\varphi_\ast F_k=d\alpha_k+d^\ast \beta_k+h_k$. Then $\ast\Delta \beta=\ast d \varphi_\ast F_k=\ast \varphi_\ast d F_k$.  Since $\varphi_\ast d F_k=\partial \varphi_\ast I_k$, $\varphi_\ast d F_k$ is supported in a finite set of points, thus $\beta$ is smooth in $N$ outside of a finite number of points. Now let $\tilde{\alpha}_k\in \Omega^{m-2}(N)$ such that $\lVert \alpha_k-\tilde{\alpha}_k\rVert_{W^{1,p}}\leq \frac{1}{2^k}$ and set $\tilde{F}_k:= d\tilde{\alpha}_k+d\beta_k+h_k$. Then by construction $\tilde{F}_k\in \Omega_{R}^{m-1}(N)$ and $\tilde{F}_k\to\varphi_\ast F$ in $\Omega_p^{m-1}(N)$ as $k\to\infty$. 
\end{proof}

Theorem \ref{Theorem: strong approximation for vector fields}, Corollary \ref{Corollary: strong approximation on boundary of cube}, Lemma \ref{Lemma: bi-Lipschitz maps preserve approximability properties} and Remark \ref{Remark: existence of connection is enough for strong approximability on cube} can be combined to obtain the following general statement.
\begin{thm}\label{Theorem: strong approximation for general domains}
Let $M\subset\r^n$ be any embedded $m$-dimensional Lipschitz submanifold of $\r^n$ which is bi-Lipschitz equivalent either to $Q^m_1(0)$ or $\partial Q^{m+1}_1(0)$. Then:
\begin{align*}
    \overline{\Omega_{p,R,\infty}^{m-1}(M)}^{L^p}=\begin{cases}
        \Omega_{p,\z}^{m-1}(M) & \mbox{ if } p\in[1,m/(m-1)),\\
        \big\{F\in\Omega_{p}^{m-1}(M) \mbox{ s.t. } *\hspace{-0.5mm}dF=0\big\} &\mbox{ if } p\in[m/(m-1),+\infty).
    \end{cases}
\end{align*}
Moreover, if $M$ is smooth we have $\overline{\Omega_{p,R}^{m-1}(M)}^{L^p}=\overline{\Omega_{p,R,\infty}^{m-1}(M)}^{L^p}$.
\end{thm}
\subsection{Corollaries of Theorem \ref{Theorem: strong approximation for vector fields}}
Finally we present a couple of Corollaries of Theorem \ref{Theorem: strong approximation for vector fields}.\\
First we show that the boundary of any $I\in\mathcal{R}_1(D)$ having finite mass can be approximated strongly in $(W^{1,\infty}_0(D))^\ast$ by finite sums of deltas with integer coefficients.
\begin{cor}\label{Corollary: Approximation of boundaries of 1 rect currents}
Let $D\subset\r^n$ be any open and bounded domain in $\r^n$ which is bi-Lipschitz equivalent to $Q^n_1(0)$. Let $I\in \mathcal{R}_1(D)$ with finite mass. Then, there exists a vector field $\,V\in L^1(D)$ such that
\begin{align*}
    \operatorname{div}(V)=\partial I \quad \text{ in }\,(W^{1,\infty}_0(D))^\ast.
\end{align*}
Thus $\partial I$ can be approximated strongly in $(W^{1,\infty}_0(D))^\ast$ by finite sums of deltas with integer coefficients. More precisely there exist sequences of points $(P_i)_{i\in \mathbb{N}}$ and $(N_i)_{i\in \mathbb{N}}$ in $D$ such that
\begin{align}\label{Equation: result of cor 1.1, approx by finite sing}
    \partial I= \sum_{i\in \mathbb{N}}(\delta_{P_i}-\delta_{N_i})\,\text{ in }(W^{1,\infty}_0(D))^\ast\text{ and }\sum_{i\in \mathbb{N}}\lvert P_i-N_i\rvert<\infty.
    \end{align}
\end{cor}
\begin{proof}
By Lemma \ref{Lemma: bi-Lipschitz maps preserve approximability properties}, it is enough to consider the case $D=Q_1(0)$. Let $I$ be as above. By \cite[Theorem 5.6]{alberti-baldo-orlandi} there exists a map $u\in W^{1,n-1}(Q_1(0),\s^{n-1})$ such that
\begin{align*}
    \ast d\left(\frac{1}{n}\sum_{i=1}^n(-1)^{i-1}u_i\bigwedge_{j\neq i}d u_j\right)=\alpha_{n-1}\partial I,
\end{align*}
where
$\alpha_{n-1}$ denotes the volume of the $(n-1)$-dimensional ball.\\
Set
\begin{align*}
    \omega:=\frac{1}{n\alpha_{n-1}}\sum_{i=1}^n(-1)^{i-1}u_i\bigwedge_{j\neq i}d u_j
\end{align*}
Notice that $\omega\in \Omega^{n-1}_{1}(Q_1(0))$ and $\ast d \omega=\partial I$.
Now let $V:=(\ast \omega)^{\sharp}$. Then $V\in L^1(Q_1(0))$ and
\begin{align*}
    \operatorname{div}(V)=\partial I.
\end{align*}
By Theorem \ref{Theorem: strong approximation for vector fields}, there exists a sequence $(V_k)_{k\in \mathbb{N}}$ in $L^1_R(Q_1(0))$ such that $V_k\to V$ in $L^1.$ Then
\begin{align*}
    \operatorname{div}(V_k)\to\operatorname{div}(V)=\partial I \quad \text{ in }\,(W^{1,\infty}_0(Q_1(0)))^\ast.
\end{align*}
As for any $k\in \mathbb{N}$ we have that $\operatorname{div}(V_k)$ is a finite sum of deltas with integer coefficients, by \cite[Proposition A.1]{Ponce2004}\footnote{Here Proposition A.1 in \cite{Ponce2004} is applied to the following metric space: for any $x,y\in Q_1(0)$ let $\overline{d}(x,y)=\min\{d(x,y),\operatorname{dist}(x,\partial Q_1(0))+\operatorname{dist}(y,\partial Q_1(0)\}$, where $d$ denotes the Euclidean distance in $Q_1(0)$. Let $(\widetilde Q_1(0), \overline{d})$ denote the completion of $Q_1(0)$ with respect to the distance $\overline{d}$. Then Lipschitz functions on $(\widetilde Q_1(0), \overline{d})$ corresponds to functions in $W^{1,\infty}_0(Q_1(0))$ (with same Lipschitz constant) modulo additive constants.} \eqref{Equation: result of cor 1.1, approx by finite sing} holds.
\end{proof}
\begin{cor}\label{Corollary: representation of partial I}
Let $M$ be a complete Lipschitz $m$-manifold, with or without boundary, compactly contained in the open cube $Q_2(0)$. Let $I\in \mathcal{R}_1(Q_2(0))$ be a rectifiable current of finite mass supported on $M$. Then there exist two sequences of points $(p_i)_{i\in \mathbb{N}}$ and $(n_i)_{i\in \mathbb{N}}$ in $M$ such that
\begin{align*}
    \partial I= \sum_{i\in \mathbb{N}}(\delta_{p_i}-\delta_{n_i})\,\text{ in }(W^{1,\infty}(Q_2(0)))^\ast\text{ and }\sum_{i\in \mathbb{N}}\lvert p_i-n_i\rvert<\infty.
\end{align*}
\end{cor}
\begin{proof}
Let $I\in\mathcal{R}_1(M)$ be a rectifiable 1-current of finite mass supported in $M\subset Q_2(0)$. By Corollary \ref{Corollary: Approximation of boundaries of 1 rect currents} there exists a vector field $V\in L^1(Q_2(0))$ such that $\operatorname{div}(V)=\partial I$. Thus we can apply the arguments of the proof of Theorem \ref{Theorem: strong approximation for vector fields, now really for vector fields} to $V$. For any $\varepsilon\in E_V$ this yields a vector field $\tilde V_\varepsilon\in L^1(\Omega_{a_\varepsilon, \varepsilon})$ with the following properties: all bad cubes $Q\in \mathscr{C}_{a_\varepsilon,\varepsilon}$ are such that $\overline{Q}\cap M\neq \emptyset$, therefore the topological singularities of $\tilde V_\varepsilon$ lie at a distance of at most $\sqrt{n}\varepsilon$ from $M$. Moreover notice that if $\varepsilon$ is sufficiently small, $\int_{Q_2(0)}\operatorname{div}(\tilde V_\varepsilon)\,d\L^n=0$ (one can see this by testing $\operatorname{div}(\tilde V_\varepsilon)$ against a function $\varphi\in C^\infty_c(Q_2(0))$ such that $\varphi\equiv 1$ in a neighbourhood of $M$).
Thus $\operatorname{div}(\tilde V_\varepsilon)$ can represented by
\begin{align*}
    \operatorname{div}(V_\varepsilon)=\sum_{i=1}^{Q^\varepsilon}(\delta_{p_i^\varepsilon}-\delta_{n_i^\varepsilon})
\end{align*}
for some $Q^\varepsilon\in\mathbb{N}$ and points $p_i^\varepsilon$ and $n_i^\varepsilon$ (possibly repeated) in a $\sqrt{n}\varepsilon$-neighbourhood of $M$.
By the argument of Lemma \ref{Lemma: volume of the bad cubes vanishes at the limit} (with $Q_2(0)$ in place of $Q_1(0)$) we have $\varepsilon Q^\varepsilon \to 0$ as $\varepsilon\to 0^+$ in $E_V$.
Now for any $i\in \{1,...,Q^\varepsilon\}$ let $\tilde p_i^\varepsilon$ and $\tilde n_i^\varepsilon$ in $M$ such that $\lvert p_i^\varepsilon-\tilde p_i^\varepsilon\rvert<2\sqrt{n}\varepsilon$ and $\lvert n_i^\varepsilon-\tilde n_i^\varepsilon\rvert<2\sqrt{n}\varepsilon$. Let $I_{p_i^\varepsilon}\in \mathcal{R}_1(Q_2(0))$ be the rectifiable current given by integration on the segment joining $p_i^\varepsilon$ and $\tilde p_i^\varepsilon$ oriented from $p_i^\varepsilon$ to $\tilde p_i^\varepsilon$ and let $I_{n_i^\varepsilon}\in \mathcal{R}_1(Q_2(0))$ be the rectifiable current given by integration on the segment joining $n_i^\varepsilon$ and $\tilde n_i^\varepsilon$ oriented from $\tilde n_i^\varepsilon$ to $n_i^\varepsilon$.
Let $I_\varepsilon\in \mathcal{R}_1(Q_2(0))$ be a rectifiable 1-current of finite mass such that $\operatorname{div}(\tilde V_\varepsilon)=\partial I_\varepsilon$.
Set $\tilde I_\varepsilon= I_\varepsilon+\sum_{i=1}^{Q^\varepsilon}(I_{p_i^\varepsilon}+I_{n_i^\varepsilon})$.
Then
\begin{align*}
    \partial \tilde I_\varepsilon=\sum_{i=1}^{Q^\varepsilon}(\delta_{\tilde p_i^\varepsilon}-\delta_{\tilde n_i^\varepsilon})
\end{align*}
is supported in $M$. Moreover we have
\begin{align*}
    \lVert \partial I-\partial \tilde I_\varepsilon \rVert_{(W^{1,\infty}(M))^\ast}\leq \lVert \partial I-\partial I_\varepsilon\rVert_{(W^{1,\infty}(Q_2(0)))^\ast}+\lVert \partial I_\varepsilon-\partial \tilde I_\varepsilon \rVert_{(W^{1,\infty}(Q_2(0)))^\ast}
\end{align*}
(here $M$ is endowed with the euclidean distance in $Q_2(0)$; notice that we are making use of the fact that any Lipschitz function on $M$ can be extended to a Lipschitz function on $Q_2(0)$ with same Lipschitz constant).
Now since $\lVert \tilde V_\varepsilon-V\rVert_{L^p(\Omega_{a_\varepsilon,\varepsilon)}}\to 0$ as $\varepsilon\to 0$ in $E_V$ and $\partial I$, $\partial I_\varepsilon$ are supported in a compact subset of $Q_2(0)$, the first term on the right hand side tends to zero as $\varepsilon\to 0^+$ in $E_V$. Moreover the second term is bounded by $4\sqrt{n}\varepsilon Q^\varepsilon$ (see for instance Lemma 2 in \cite{brezis}) and thus tends to $0$ as $\varepsilon\to 0^+$ in $E_V$.
This shows that $\partial I$ belongs to the (strong) $(W^{1,\infty}(M))^\ast$ closure of the class of $0$-currents $T$ on $M$ such that
\begin{align}\label{Equation: explicit form for T, sum of deltas in M}
    T=\sum_{j\in J}(\delta_{p_j}-\delta_{n_j})\text{ in }(W^{1,\infty}(M))^\ast\text{ and }\sum_{j\in J}\lvert p_j- n_j\rvert<\infty
\end{align}
for a countable set $J$ and points $p_j$, $n_j$ in $M$. By \cite[Proposition A.1]{Ponce2004} applied to the complete metric space $(M, d)$ (where $d$ denotes the Euclidean distance in $Q_2(0)$) this space is closed in $(W^{1,\infty}(M))^\ast$, therefore $\partial I$ is also of this form.
Since any Lipschitz function $\varphi\in W^{1,\infty}(Q_2(0))$ has a Lipschitz trace $\varphi\big\vert_{M}$ on $M$
and $\langle \partial I, \varphi\rangle_{Q_2(0)}=\langle \partial I, \varphi\big\vert_M\rangle_M$, we conclude that $\partial I$ can be represented as in \ref{Equation: explicit form for T, sum of deltas in M} also as an element of $(W^{1,\infty}(Q_2(0))^\ast$.
\end{proof}
Theorem $\ref{Theorem: strong approximation for vector fields}$ could also be useful to obtain approximation results for Sobolev maps with values into manifolds. For instance we can use it to recover the following result, due to R. Schoen and K. Uhlenbeck (for $p=2$, see \cite{schoen-uhlenbeck-2}, Section 4) and F. Bethuel and X. Zheng (for $p>2$, see \cite[Theorem 4]{bethuel-zheng}).
\begin{cor}
Let $u\in W^{1,p}(Q_1(0), \mathbb{S}^1)$ for some $p\in (1,\infty)$.\\
If $p\geq 2$, then
\begin{align*}
    \operatorname{div}(u\wedge \nabla^\perp u)=0\quad\text{ in }\mathcal{D}'(Q_1(0))
\end{align*}
and $u$ can be approximated in $W^{1,p}$ by a sequence of functions in $C^\infty(Q_1(0), \mathbb{S}^1).$\\
If $p<2$, then
\begin{align}
\label{Equation: form of weak Jacobian of maps to S^1}
    \frac{1}{2\pi}\operatorname{div}(u\wedge \nabla^\perp u)=\partial I,
\end{align}
where $I\in \mathcal{R}_1(Q_1(0))$ is a 1-rectifiable current of finite mass, and $u$ can be approximated in $W^{1,p}$ by a sequence of functions in
\begin{align*}
    \mathcal{R}:=\left\{v\in W^{1,p}(Q_1(0),\mathbb{S}^1); v\in C^\infty(Q_1(0)\smallsetminus A, \mathbb{S}^1),\text{ where }A\text{ is some finite set}\right\}.
\end{align*}
\begin{proof}
First we claim that the vector field $u\wedge \nabla^\perp u$ belongs to $L^p_\mathbb{Z}(Q_1(0))$. In fact notice that for any $x_0\in Q_1(0)$, for a.e. $\rho\in (0, 2\dist_\infty(x_0,\partial Q_1(0)))$ $\partial Q_\rho(x_0)$ consists $\H^{n-1}$-a.e. of Lebesgue points of $u\wedge \nabla^\perp u$. Moreover for almost any such $\rho$ we have
\begin{align*}
   \frac{1}{2\pi}\int_{\partial Q_\rho(x_0)}(u\wedge\nabla^\perp u)\cdot \nu_{\partial Q_\rho(x_0)}=\deg \left(u\big\vert_{\partial Q_\rho(x_0)}\right)\in \mathbb{Z}.
\end{align*}
Hence the vector field $\frac{1}{2\pi}u\wedge\nabla^\perp u$ belongs to $L^p_\mathbb{Z}(Q_1(0))$.\\
Thus if $p\geq 2$ by Theorem \ref{Theorem: strong approximation for vector fields} there holds $\operatorname{div}(u\wedge\nabla^\perp u)=0$ in $\mathcal{D}'(Q_1(0))$, while if $p<2$
there exists a sequence of vector fields $(V_n)_{n\in \mathbb{N}}$ in $L^p_R(Q_1(0))$ such that
\begin{align*}
    V_n\to \frac{1}{2\pi}u\wedge \nabla^\perp u\quad\text{ in }L^p(Q_1(0))\text{ as }n\to\infty.
\end{align*}
For any $n\in \mathbb{N}$ by Hodge decomposition there exist $a_n\in W^{1,p}(Q_1(0))$, $b_n\in W^{1,p}_0(Q_1(0))$ such that
\begin{align*}
   2\pi V_n=\nabla^\perp a_n+\nabla b_n.
\end{align*}
For any $n\in \mathbb{N}$ let $\tilde a_n\in C^\infty(Q_1(0))$ be such that $\lVert \tilde a_n-a_n\rVert_{L^p}\leq \frac{1}{n}$. Moreover notice that there exists $d_n\in W^{1,p}(Q_1(0),\mathbb{S}^1)\cap C^\infty(Q_1(0)\smallsetminus A)$, where $A$ is a finite set, such that
\begin{align*}
    \nabla b_n=d_n\wedge \nabla^\perp d_n.
\end{align*}
In fact
\begin{align*}
    \Delta b_n=2\pi\operatorname{div}(V_n)=2\pi\sum_{i=1}^{Q^n}d_i^n\delta_{p_i^n}
\end{align*}
for some $Q^n\in \mathbb{N}$, $p_i^n\in Q_1(0)$ and $d_i^n\in \mathbb{Z}$, thus
$b_n=-\sum_{i=1}^{Q^n}\log\lvert x-p_i^n\rvert^{d_i^n}+h_n$ for an harmonic function $h_n$. Then $d_n$ can be chosen to be
\begin{align*}
   d_n(x)= e^{-i\tilde{h}_n}\prod_{i=1}^{Q^n} \left(\frac{x-p_i^n}{\lvert x-p_i^n\rvert}\right)^{d_n},
\end{align*}
where $\tilde{h}_n$ is the harmonic conjugate of $h_n$ (the product has to be understood as complex multiplication in $\mathbb{C}\simeq\mathbb{R}^2$).\\
For any $n\in \mathbb{N}$ set $u_n:=e^{i\tilde a_n}d_n$.
Then by construction $u_n\in \mathcal{R}$ and
\begin{align*}
    u_n\wedge \nabla^\perp u_n\to u\wedge \nabla^\perp u \quad\text{ in }L^p(Q_1(0))\text{ as }n\to \infty.
\end{align*}
Therefore there is $c\in [0,2\pi)$ so that up to a subsequence
\begin{align*}
    e^{ic}u_n\to u\quad\text{ in }L^p(Q_1(0))\text{ as }n\to\infty.
\end{align*}
\end{proof}
\end{cor}
\begin{rem}
Equation (\ref{Equation: form of weak Jacobian of maps to S^1}) was obtained in \cite[Theorem 3']{brezis-mironescu-ponce} with the help of the approximation result of F. Bethuel and X. Zheng.
\end{rem}
\section{The weak $L^p$-closure of $\Omega_{p,\z}^{n-1}(Q_1^n(0))$}
In the present section we follow the ideas presented in \cite{petrache-riviere-abelian} in order to prove that the space $\Omega_{p,\z}^{n-1}(Q^n_1(0))$ is weakly sequentially closed for every $n\ge 2$ and $p\in(1,+\infty)$. The main reason why such techniques couldn't be used before in this context for $n\neq 3$ was the lack of a strong approximation theorem like Theorem \ref{Theorem: strong approximation for vector fields} for general dimension $n$. Such result is needed in order to define a suitable notion of distance between the cubical slices of a form $F\in\Omega_{p,\z}^{n-1}(Q^n(0))$, given by $(x\mapsto x_0+\rho x)^*i_{\partial Q_{\rho}(x_0)}^*F$ for $\L^1$-a.e. $\rho$ (see subsections 3.1 and 3.2 for the precise definition). Once we have turned the space of the cubical slices of $F$ into a metric space, we will show that the "slice function" associated to $F$, given by $\rho\mapsto(x\mapsto x_0+\rho x)^*i_{\partial Q_{\rho}(x_0)}^*F$, is locally $\frac{1}{p'}-$H\"older continuous (see subsection 3.3). Moreover we will see that if $\{F_k\}_{k\in\n}\subset\Omega_{p,\z}^{n-1}(Q_1(0))$ converges weakly in $L^p$, then the sequence of the slice functions associated to each $F_k$ is locally uniformly $\frac{1}{p'}-$H\"older continuous. Finally, we will use the previous facts together with some technical lemmata to conclude the proof of Theorem \ref{Theorem: weak closure} for $D=Q^n_1(0)$. Notice that by Theorem 
\ref{Theorem: strong approximation for vector fields} the result is clear if $p\in \big[ n/(n-1),\infty)$, here we will focus on the case $p\in\big(1,n/(n-1)\big)$.
\subsection{Slice distance on $\s^{n-1}$}
Throughout the following section, we will assume that $p\in\big(1,n/(n-1)\big)$. Moreover, we will denote by "$\ast$" the Hodge star operator associated with the standard round metric on $\s^{n-1}$. We will denote by $Z$ the linear subspace of $\Omega_{p}^{n-1}(\s^{n-1})$ given by
\begin{align*}
    Z:=\bigg\{h\in\Omega_p^{n-1}(\s^{n-1}) \mbox{ s.t. } \int_{\s^{n-1}}h\in\z\bigg\}.
\end{align*}
\begin{rem}
It's clear that $Z$ is weakly (and thus strongly) $L^p$ closed in $\Omega_{p}^{n-1}(\s^{n-1})$. Indeed, let $\{h_k\}_{k\in\n}\subset Z$ be any sequence such that $h_k\rightharpoonup h$ weakly in $\Omega_p^{n-1}(\s^{n-1})$, i.e.
\begin{align*}
    \int_{\s^{n-1}}\varphi h_k\rightarrow\int_{\s^{n-1}}\varphi h, \qquad\forall\,\varphi\in L^{p'}(\s^{n-1}).
\end{align*}
Then, the statement follows by picking $\varphi\equiv 1$ and noticing that a convergent sequence of integer numbers is definitively constant.
\end{rem}
Fix any arbitrary point $q\in\s^{n-1}$. We define the functions $d,\tilde d:Z\times Z\rightarrow [0,+\infty]$ by
\begin{align*}
    d(h_1,h_2):=\inf\bigg\{||\alpha||_{L^p}\, \mbox{ s.t. } *\hspace{-0.5mm}(h_1-h_2)=d^*\alpha+\partial I+\bigg(\int_{\s^{n-1}}h_1-h_2\bigg)\delta_q\bigg\}, 
\end{align*}
with $\alpha\in\Omega_p^1(\s^{n-1}), I\in\R_1(\s^{n-1})$, and
\begin{align*}
    \tilde d(h_1,h_2):=\inf\bigg\{||\alpha||_{L^p}\, \mbox{ s.t. } *\hspace{-0.5mm}(h_1-h_2)=d^*\alpha+\partial I+\bigg(\int_{\s^{n-1}}h_1-h_2\bigg)\delta_q\bigg\},
\end{align*}
with $\alpha\in\Omega_p^1(\s^{n-1})$, $I\in\R_1(\s^{n-1})\cap\mathcal{N}_1(\s^{n-1})$. 
\begin{rem}[$d$ and $\tilde d$ are always finite on $Z$]
\label{remark the distances are finite}
We claim that $d,\tilde d<+\infty$. Since obviously $d(h_1,h_2)\le\tilde d(h_1,h_2)$, it is enough to show that $\tilde d(h_1,h_2)<+\infty$, for every $h_1,h_2\in Z$. This just amounts to saying that given any $h_1,h_2\in Z$ we can always find $\alpha\in\Omega_p^1(\s^{n-1})$, $I\in\R_1(\s^{n-1})\cap\mathcal{N}_1(\s^{n-1})$ satisfying 
\begin{align*}
    *(h_1-h_2)=d^*\alpha+\partial I+\bigg(\int_{\s^{n-1}}h_1-h_2\bigg)\delta_q.
\end{align*}
Indeed, let 
\begin{align*}
    a:=\int_{\s^{n-1}}h_1-h_2\in\z.
\end{align*}
Consider the following first order differential system on $\s^{n-1}$:
\begin{align*}
    \begin{cases}
        d^*\omega=*(h_1-h_2)-a\delta_q=:F,\\
        d\omega=0.
    \end{cases}
\end{align*}
Since $p\in\big(1,n/(n-1)\big)$, $F\in\mathcal{L}(W^{1,p'}(\s^{n-1}))$. Moreover, $\langle F,1\rangle=0$. Hence, by Lemma \ref{appendix lemma d* e d}, we know that the previous differential system has a solution $\alpha\in\Omega_p^1(\s^{n-1})$ and the statement follows. 
\end{rem}
\begin{rem}
\label{remark equivalence distances}
As observed above, it is clear that $d(h_1,h_2)\le\tilde d(h_1,h_2)$, for every $h_1,h_2\in Z$. We claim that actually $d(h_1,h_2)=\tilde d(h_1,h_2)$, for every $h_1,h_2\in Z$. In order to prove the remaining inequality, fix any $h_1,h_2\in Z$ and let $\{\alpha_k\}_{k\in\n}\subset\Omega_p^1(\s^{n-1})$, $\{I_k\}_{k\in\n}\subset\R_1(\s^{n-1})$ be such that 
\begin{align*}
    \begin{cases}
        *(h_1-h_2)=d^*\alpha_k+\partial I_k+a\delta_q,\quad\forall\,k\in\n,\\
        \lvert\lvert\alpha_k\rvert\rvert_{L^p}\rightarrow d(h_1,h_2) \mbox{ as } k\rightarrow\infty,
    \end{cases}
\end{align*}
with
\begin{align*}
    a:=\int_{\s^{n-1}}h_1-h_2.
\end{align*}
By Corollary \ref{appendix corollary laplace}, the linear differential equation
\begin{align*}
    \Delta u=*(h_1-h_2)-a\delta_q
\end{align*}
has a weak solution $\psi\in\dot W^{1,p}(\s^{n-1})$. Let $\omega_k:=d\psi-\alpha_k$, for every $k\in\n$. Notice that
\begin{align*}
    d^*\omega_k=\partial I_k, \qquad\forall\, k\in\n.
\end{align*}
By Theorem \ref{Theorem: strong approximation for general domains}, for every $k\in\n$ there exists a sequence $\{\omega_{k}^j\}_{j\in\n}\subset\Omega_{p,R}^1(\s^{n-1})$ such that $r_k^j:=\omega_k-\omega_k^j\rightarrow 0$ strongly in $L^p$ as $j\rightarrow\infty$. By construction, it follows that
\begin{align*}
    *(h_1-h_2)=d^*(\alpha_k+r_k^j)+d^*\omega_k^j+a\delta_q,\qquad\forall\,k,j\in\n.
\end{align*}
We observe that by Proposition \ref{appendix proposition existence minimal connection finitely many singularities}  for every $k,j\in\n$ there exist $I_k^j\in\R_1(\s^{n-1})\cap\mathcal{N}_1(\s^{n-1})$ such that $d^*\omega_k^j=\partial I_k^j$. This implies that 
\begin{align*}
    \tilde d(h_1,h_2)\le\lvert\lvert\alpha_k+r_k^j\rvert\rvert_{L^p}\le\lvert\lvert\alpha_k\rvert\rvert_{L^p}+\lvert\lvert r_k^j\rvert\rvert_{L^p},\qquad\forall\,k,j\in\n.
\end{align*}
By letting first $j\rightarrow\infty$ and then $k\rightarrow\infty$ in the previous inequality, our claim follows. 
\end{rem}
\begin{prop}
\label{d is a metric}
$(Z,d)$ is a metric space. 
\begin{proof}
We need to check symmetry, triangular inequality and non-degeneracy.

\textit{Symmetry}. This is clear since both the $L^p$-norm and the space $\R_1(\s^{n-1})$ are invariant under sign change.

\textit{Triangular inequality}. Let $h_1,h_2,h_3\in Z$. By definition of infimum, for every $\varepsilon>0$ we can write
\begin{align*}
    \begin{cases}
        \displaystyle{*(h_1-h_2)=d^*\alpha_{\varepsilon}+\partial I_{\varepsilon}+\bigg(\int_{\s^{n-1}}h_1-h_2\bigg)\delta_{q}},\\
        \displaystyle{*(h_2-h_3)=d^*\alpha_{\varepsilon}'+\partial I_{\varepsilon}'+\bigg(\int_{\s^{n-1}}h_2-h_3\bigg)\delta_{q}},
    \end{cases}
\end{align*}
with $\alpha_{\varepsilon},\alpha_{\varepsilon}'\in\Omega_p^1(\s^{n-1})$ and $I_{\varepsilon},I_{\varepsilon}'\in\R_1(\s^{n-1})$ satisfying 
\begin{align*}
    \begin{cases}
        \lvert\lvert\alpha_{\varepsilon}\rvert\rvert_{L^p}\le d(h_1,h_2)+\varepsilon\\
        \lvert\lvert\alpha_{\varepsilon}'\rvert\rvert_{L^p}\le d(h_2,h_3)+\varepsilon.
    \end{cases}
\end{align*}
We notice that
\begin{align*}
    *(h_1-h_3)=d^*(\alpha_{\varepsilon}+\alpha_{\varepsilon}')+\partial(I_{\varepsilon}+I_{\varepsilon}')+\bigg(\int_{\s^{n-1}}h_1-h_3\bigg)\delta_{q}, \qquad\forall\,\varepsilon>0.
\end{align*}
Then, by definition of $d$, we have
\begin{align*}
    d(h_1,h_3)\le\lvert\lvert\alpha_{\varepsilon}+\alpha_{\varepsilon}'\rvert\rvert_{L^p}\le\lvert\lvert\alpha_{\varepsilon}\rvert\rvert_{L^p}+\lvert\lvert\alpha_{\varepsilon}'\rvert\rvert_{L^p}\le d(h_1,h_2)+d(h_2,h_3)+2\varepsilon, \qquad\forall\,\varepsilon>0.
\end{align*}
By letting $\varepsilon\rightarrow 0^+$ in the previous inequality, we get our claim. 

\textit{Non-degeneracy}. Assume that $d(h_1,h_2)=0$, for some $h_1,h_2\in Z$. Let 
\begin{align*}
    a:=\int_{\s^{n-1}}h_1-h_2\in\z.
\end{align*}
Then, since $d=\tilde d$ (see Remark \ref{remark equivalence distances}) and by definition of $\tilde d$, there exist $\{\alpha_k\}_{k\in\n}\subset\Omega_p^1(\s^{n-1})$ and $\{I_k\}_{k\in\n}\subset\R_1(\s^{n-1})\cap\mathcal{N}_1(\s^{n-1})$ such that
\begin{align*}
    *(h_1-h_2)=d^*\alpha_k+\partial I_k+a\delta_q
\end{align*}
and $\alpha_k\rightarrow 0$ strongly in $L^p$ as $k\to\infty$. Observe that
\begin{align*}
    \partial I_k\to \ast(h_1-h_2)-a\delta_q\quad\text{ in }(W^{1,\infty}(\s^{n-1}))^\ast.
\end{align*}
Now for any $k\in \mathbb{N}$, $\partial I_k$ can be represented as
\begin{align*}
    \partial I_k=\sum_{j=1}^{J_k}(\delta_{p_j^k}-\delta_{n_j^k})
\end{align*}
for some $J_k\in \mathbb{N}$ and points $p_j^k$, $n_j^k$ in $\s^{n-1}$.
But the space of distributions of the form 
\begin{align}\label{Equation: Distribution as sum of dipoles}
    \sum_{j\in J}(\delta_{p_j}-\delta_{n_j})\,\text{ such that }\,\sum_{j\in J}\lvert p_j- n_j\rvert<\infty
\end{align}
(for a countable set $J$ and points $p_j$, $n_j$ in $\s^{n-1}$) is closed with respect to the (strong) topology of $({W}^{1,\infty}(\s^{n-1}))^\ast$ (see Proposition A.1 in \cite{Ponce2004}), thus there exists a distribution $T$ as in \eqref{Equation: Distribution as sum of dipoles} such that
\begin{align}\label{Equation: Contradiction, Lp equal sum of deltas}
    \ast\,(h_1-h_2)=T-a\delta_q.
\end{align}
But this implies that $\ast\,(h_1-h_2)=0$, since the left-hand-side in \eqref{Equation: Contradiction, Lp equal sum of deltas} is in $L^p$ whilst the right-hand-side is in $L^p$ if and only if it is equal to zero. 
\end{proof}
\end{prop}
\begin{rem}
Notice that the proof above relies on the fact that $d=\tilde{d}$, which was proved using Corollary \ref{Corollary: strong approximation on boundary of cube}. As the proof of the Corollary \ref{Corollary: strong approximation on boundary of cube} was rather cumbersome, we remark here that there is a way to skip that passage. Indeed, in the proof of the non-degeneracy of $d$ it is not necessary to assume that $\{I_k\}_{k\in \mathbb{N}}$ lies in $\mathcal{N}_1(\s^{n-1})$. In fact, it follows from Corollary \ref{Corollary: representation of partial I} that $\partial I_k$ is of the form (\ref{Equation: Distribution as sum of dipoles}), and thus the limit of $\{I_k\}_{k\in \mathbb{N}}$ in $(W^{1,\infty}(\s^{n-1}))^\ast$ will also be of that form.
\end{rem}
\begin{prop}
\label{d metrizes the weak topology on bounded subsets}
Let $\{h_k\}_{k\in\n}\subset Z$ and $h\in Z$. Then the following are equivalent:
\begin{enumerate}
    \item $\{h_k\}_{k\in\n}\subset Z$ is uniformly bounded w.r.t. the $L^p$-norm and $d(h_k,h)\rightarrow 0$ as $k\rightarrow\infty$;
    \item $h_k\rightharpoonup h$ weakly in $L^p$ as $k\rightarrow\infty$.
\end{enumerate}
\begin{proof}
We prove separately the two implications. 

\textit{$2\Rightarrow 1$}. Pick any subsequence of $\{h_k\}_{k\in\n}$ (not relabelled). For every $k\in\n$, let 
\begin{align*}
    a_k:=\langle h_k-h,1\rangle=\int_{\s^{n-1}}h_k-h.
\end{align*}
Since $h_k\rightharpoonup h$ weakly in $L^p$ as $k\rightarrow\infty$, it follows that $a_k\rightarrow 0$ as as $k\rightarrow\infty$. Since $\{a_k\}_{k\in\n}\subset\z$, there exists $K\in\n$ such that $a_k=0$ for every $k\ge K$. Fix any $k\ge K$. By
Lemma \ref{appendix lemma d* e d}, the linear differential system 
\begin{align*}
    \begin{cases}
        d^*\omega=*(h_k-h)\\
        d\omega=0,
    \end{cases}
\end{align*}
respectively (if $n=2$)
\begin{align*}
    \begin{cases}
        d^*\omega=*(h_k-h)\\
        d\omega=0\\
        \displaystyle{\int_{\s^1}\omega=0},
    \end{cases}
\end{align*}
has a unique weak solution $\alpha_k\in\Omega_p^1(\s^{n-1})$.
By Remark \ref{appendix remark sobolev forms} we have
\begin{align*}
    \lvert\lvert\alpha_k\rvert\rvert_{W^{1,p}}\le C\big(\lvert\lvert d\alpha_k\rvert\rvert_{L^p}+\lvert\lvert d^*\alpha_k\rvert\rvert_{L^p}\big)=C\lvert\lvert h_k-h\rvert\rvert_{L^p}.
\end{align*}
Since $\{h_k\}_{k\in\n}$ is weakly convergent, we know that it is also uniformly bounded w.r.t the $L^p$-norm. Then $\{\alpha_k\}_{k\ge K}$ is uniformly bounded w.r.t the $W^{1,p}$-norm. Hence, by weak compactness in $W^{1,p}$, there exists a subsequence $\{\alpha_{k_l}\}_{l\in\n}\subset\{\alpha_k\}_{k\ge K}$ and a one-form $\alpha\in \Omega^1_{W^{1,p}}(\s^{n-1})$ such that $\alpha_{k_l}\rightharpoonup\alpha$ weakly in $W^{1,p}$. By Rellich-Kondrakov theorem, it follows that $\alpha_{k_l}\rightarrow \alpha$ strongly in $L^p$. We claim that $\alpha=0$. Indeed,
\begin{align*}
    \langle\alpha,\omega\rangle_{L^p-L^{p'}}&=\lim_{l\rightarrow\infty}\langle\alpha_{k_l},d\varphi+d^*\beta\rangle_{L^p-L^{p'}}\\
    &=\lim_{l\rightarrow\infty}-\langle d^\ast \alpha_{k_l}, \varphi\rangle_{L^p-L^{p'}}\\
    &=\lim_{l\rightarrow\infty}-\int_{\s^{n-1}} (h_{k_l}-h)\wedge \varphi=0, \qquad\forall\, \omega=d\varphi+d^*\beta\in\Omega^1(\s^{n-1}),
\end{align*}
respectively (if $n=2$)
\begin{align*}
    \langle\alpha,\omega\rangle_{L^p-L^{p'}}&=\lim_{l\rightarrow\infty}\langle\alpha_{k_l},d\varphi+d^*\beta+\eta\rangle_{L^p-L^{p'}}\\
    &=\lim_{l\rightarrow\infty}\langle\alpha_{k_l},d\varphi+d^*\beta\rangle_{L^p-L^{p'}}\\
    &=\lim_{l\rightarrow\infty}-\langle d^\ast \alpha_{k_l}, \varphi\rangle_{L^p-L^{p'}}\\
    &=\lim_{l\rightarrow\infty}-\int_{\s^{n-1}} (h_{k_l}-h)\wedge \varphi=0, \qquad\forall\, \omega=d\varphi+d^*\beta+\eta\in\Omega^1(\s^1),
\end{align*}
where $\eta\in\Omega^1(\s^1)$ is a harmonic $1$-form on $\s^1$ (hence a constant $1$-form) and the second equality follows because $\alpha_{k_l}$ is distributionally closed.

Hence, we have shown that $\alpha_{k_l}\rightarrow 0$ strongly in $L^p$ as $l\to\infty$. As $\ast(h_{k_l}-h)=d^*\alpha_{k_l}$, for every $l\in\n$, we have 
\begin{align*}
    d(h_{k_l},h)\le\lvert\lvert\alpha_{k_l}\rvert\rvert_{L^p}\rightarrow 0, \qquad\mbox{ as } l\rightarrow\infty.
\end{align*}
We have just proved that any subsequence of $\{h_k\}_{k\in\n}$ has a further subsequence converging to $h$ with respect to $d$, therefore $1.$ follows. 

\textit{$1\Rightarrow 2$}. Pick any subsequence of $\{h_k\}_{k\in\n}$ (not relabelled). Since $\{h_k\}_{k\in\n}\subset Z$ is uniformly bounded w.r.t. the $L^p$-norm, by weak $L^p$-compactness there exists a subsequence $\{h_{k_l}\}_{l\in\n}$ of $\{h_k\}_{k\in\n}$ and a $h_w\in Z$ such that $h_{k_l}\rightharpoonup h_w$ weakly in $L^p$. Since we have just shown that $2\Rightarrow 1$, we know that $d(h_{k_l},h_w)\rightarrow 0$ as $l\rightarrow\infty$. By uniqueness of the limit, we get $h_w=h$. We have just proved that any subsequence of $\{h_k\}_{k\in\n}\subset Z$ has a further subsequence converging to $h$ weakly in $L^p$, hence, $2.$ follows. 
\end{proof}
\end{prop}
\subsection{Slice distance on $\partial Q^n_1(0)$}
Let $Q_1(0)\subset\r^n$ be the unit cube in $\r^n$ centered at the origin and let $\Psi:\s^{n-1}\rightarrow\partial Q_1(0)$ be a bi-Lipschitz homeomorphism. We let $Y$ be the linear subspace of $\Omega_{p}^{n-1}(\partial Q_1(0))$ given by
\begin{align*}
    Y:=\bigg\{h\in\Omega_p^{n-1}(\partial Q_1(0)) \mbox{ s.t. } \int_{\partial Q_1(0)}h\in\z\bigg\}.
\end{align*}
\begin{rem}
Notice that $h\in Y$ if and only if $\Psi^*h\in Z$. Indeed, given any $h\in Z$ we have
\begin{align*}
    C_{\Psi}^{-1}\int_{\partial Q_1(0)}\lvert h\rvert^p\, d\H^{n-1}\le\int_{\s^{n-1}}\lvert\Psi^*h\rvert^p\, d\H^{n-1}\le C_{\Psi}\int_{\partial Q_1(0)}\lvert h\rvert^p\, d\H^{n-1},
\end{align*}
with $C_{\Psi}:=\left(\max\{\lVert d\Psi\rVert_{L^{\infty}},\lVert d\Psi^{-1}\rVert_{L^{\infty}}\}\right)^{(n-1)(p-1)}$, and
\begin{align*}
    \int_{\s^{n-1}}\Psi^*h=\int_{\partial Q_1(0)}h.
\end{align*}
\end{rem}
Thus, the functions $d_{\Psi},\tilde d_{\Psi}:Y\times Y\rightarrow[0,+\infty)$ given by
\begin{align*}
    d_{\Psi}(h_1,h_2)&:=d(\Psi^*h_1,\Psi^*h_2) \qquad\forall\,h_1,h_2\in Y,\\
    \tilde d_{\Psi}(h_1,h_2)&:=\tilde d(\Psi^*h_1,\Psi^*h_2) \qquad\forall\,h_1,h_2\in Y,
\end{align*}
are well-defined and coincide on $Y\times Y$ by Remarks \ref{remark the distances are finite} and \ref{remark equivalence distances}. Moreover, $(Y,d_{\Psi})$ is a metric space as a direct consequence of Proposition \ref{d is a metric} and the following statement is a corollary of Proposition \ref{d metrizes the weak topology on bounded subsets}.
\begin{cor}
\label{d_Psi metrizes the weak topology on bounded subsets}
Let $\{h_k\}_{k\in\n}\subset Y$ and $h\in Y$. Then, the following are equivalent:
\begin{enumerate}
    \item $\{h_k\}_{k\in\n}\subset Y$ is uniformly bounded w.r.t. the $L^p$-norm and $d_{\Psi}(h_k,h)\rightarrow 0$ as $k\rightarrow\infty$;
    \item $h_k\rightharpoonup h$ weakly in $L^p$ as $k\rightarrow\infty$.
\end{enumerate}
\end{cor}
\begin{rem}
Let $\Psi_1,\Psi_2:\s^{n-1}\rightarrow\partial Q_1(0)$ be bi-Lipschitz homeomorphisms. We claim that the distances $d_{\Psi_1}$ and $d_{\Psi_2}$ induced on $Y$ by $\Psi_1$ and $\Psi_2$ respectively are equivalent. Indeed notice that given any bi-Lipschitz map $\Lambda:\s^{n-1}\rightarrow\s^{n-1}$ we have 
\begin{align}
\label{estimate Psi}
    d(\Lambda^*h_1,\Lambda^*h_2)\le\lVert d\Lambda\rVert_{L^{\infty}}^{n-2}\lVert d\Lambda^{-1}\rVert_{L^{\infty}}^{\frac{n-1}{p}}d(h_1,h_2), \qquad\forall\, h_1,h_2\in Z.
\end{align}
To see this notice that if $\alpha\in\Omega_p^1(\s^{n-1})$ is a competitor in the definition of $d(h_1,h_2)$ then the form given by $(-1)^{n-2}\ast\Lambda^*(\ast\alpha)\in\Omega_p^1(\s^{n-1})$ is a competitor in the definition of $d(\Lambda^*h_1,\Lambda^*h_2)$. Hence
\begin{align*}
    d(\Lambda^*h_1,\Lambda^* h_2)\le\lVert\ast\Lambda^*(\ast\alpha)\rVert_{L^p}\le\lVert d\Lambda\rVert_{L^{\infty}}^{n-2}\lVert d\Lambda^{-1}\rVert_{L^{\infty}}^{\frac{n-1}{p}}\lVert\alpha\rVert_{L^p},
\end{align*}
for every competitor $\alpha$ in the definition of $d(h_1,h_2)$. By taking the infimum on all the competitors in the previous inequality, \eqref{estimate Psi} follows.
By applying \eqref{estimate Psi} we obtain
\begin{align*}
    d_{\Psi_2}(h_1,h_2)&=d(\Psi_2^*h_1,\Psi_2^*h_2)=d\big((\Psi_1^{-1}\circ\Psi_2)^*\Psi_1^*h_1,(\Psi_1^{-1}\circ\Psi_2)^*\Psi_1^*h_2\big)\\
    &\le C_{\Psi_1\Psi_2}d(\Psi_1^*h_1,\Psi_1^*h_2)=C_{\Psi_1\Psi_2}d_{\Psi_1}(h_1,h_2) \qquad\forall\,h_1,h_2\in Y,
\end{align*}
Analogously, we get 
\begin{align*}
    d_{\Psi_1}(h_1,h_2)&\le C_{\Psi_1\Psi_2}d_{\Psi_2}(h_1,h_2) \qquad\forall\,h_1,h_2\in Y,
\end{align*}
with $C_{\Psi_1\Psi_2}:=\max\left\{\lVert d(\Psi_2^{-1}\circ\Psi_1)\rVert_{L^{\infty}}, \lVert d(\Psi_1^{-1}\circ\Psi_2)\rVert_{L^{\infty}}\right\}^{n-2+\frac{n-1}{p}}$. 
\end{rem}
\subsection{Slice functions and their properties}
\begin{dfn}[Slice functions]
	Let $F\in\Omega_{p,\z}^{n-1}(Q_1(0))$. Given any arbitrary $x_0\in Q_1(0)$, we let $\rho_0:=2\dist_{\infty}(x_0,\partial Q_1(0))$. 
	
	We call the \textit{slice function of $F$ at $x_0$} the map $s:\dom(s)\subset (0,\rho_0)\rightarrow Y$ given by
	\begin{align*}
		s(\rho):=(x\mapsto\rho x+x_0)^*i_{\partial Q_{\rho}(x_0)}^*F, \qquad\forall\, \rho\in\dom(s),
	\end{align*}
	where $\dom(s)$ is the subset of $(0,\rho_0)$ defined as follows: $\rho\in\dom(s)$ if and only if the following conditions hold:
	\begin{enumerate}
		\item $\H^{n-1}$-a.e. point in $\partial Q_{\rho}(x_0)$ is a Lebesgue point for $F$,
        \item $\lvert F\rvert\in L^p(\partial Q_{\rho}(x_0),\H^{n-1})$,
		\item $\rho$ is a Lebesgue point for the $L^p$-function 
		\begin{align*}
		(0,\rho_0)\ni\rho\mapsto\int_{\partial Q_{\rho}(x_0)}i_{\partial Q_{\rho}(x_0)}^*F,
		\end{align*}
		\item $(x\mapsto\rho x+x_0)^*i_{\partial Q_{\rho}(x_0)}^*F\in Y$.
	\end{enumerate}
\end{dfn}
\begin{rem}\label{Remark: Lp norm of slice function}
Notice that $Dom(s)$ has $\L^1$ full measure in $(0,\rho_0)$. Moreover $s\in L^p\big((0,\rho_0);Y\big)$, in the following sense:
letting $j_\rho: x\mapsto \rho x+x_0$, we have
	\begin{align*}
	    \lVert s(\rho)\rVert_{L^p}^p=\int_{\partial Q_1(0)}\lvert j_\rho^\ast F\rvert^p \,d\mathscr{H}^{n-1}=\int_{\partial Q_\rho(x_0)}\lvert F\rvert^p\rho^{(n-1)(p-1)}\,d\mathscr{H}^{n-1}
	\end{align*}
	and thus
	\begin{align*}
		\int_0^{\rho_0}\lvert\lvert s(\rho)\rvert\rvert_{L^p}^p\, d\rho\leq \int_0^{\rho_0}\int_{\partial Q_\rho}\lvert F\rvert^p \rho^{(n-1)(p-1)}\, d \mathscr{H}^{n-1}d\rho\leq 2^n\int_{Q_{\rho_0}(x_0)}\lvert F\rvert^p d\mathscr{H}^n
	\end{align*}
\end{rem}
\begin{prop}
\label{the slice functions are (1-1/p)-holder continuous} 
Let $x_0\in Q_1(0)$ and set $\rho_0:=2\dist(x_0,\partial Q_1(0))$. Fix any $F\in\Omega_{p,\z}^{n-1}(Q_1(0))$ and let $s\in L^p\big((0,\rho_0),Y\big)$ be the slice function of $F$ at $x_0$. Let $K\subset (0,\rho_0)$ be compact.
Then, there exists a subset $E\subset K$ such that $\L^1(K\smallsetminus E)=0$ and a representative $\tilde s$ of $s$ defined pointwise on $E$ such that
\begin{align}
\label{holderianity estimate}
    d_{\Psi}\big(\tilde s(\rho_1),\tilde s(\rho_2)\big)\le C_{p,K,\Psi}\lvert\lvert F\rvert\rvert_{L^p}\lvert\rho_1-\rho_2\rvert^{\frac{1}{p'}}, \qquad\forall\,\rho_1,\rho_2\in E,
\end{align}
with
\begin{align*}
	C_{p,K,\Psi}:=C_{p,\Psi}\max_{\rho\in K}\rho^{1-n}
\end{align*}
\begin{proof}
Denote by $T_F\in\D_1(Q_1(0))$ the $1$-current on $Q_1(0)$ given by
\begin{align*}
    \langle T_F,\omega\rangle=\int_{Q_1(0)}F\wedge\omega, \qquad\forall\,\omega\in\D^1(Q_1(0)).
\end{align*}
Since $F\in\Omega_{p,\z}^{n-1}(Q_1(0))$, by Theorem \ref{Theorem: characterization of the integer valued fluxes class} there exists $I\in\R_1(Q_1(0))$ such that $\mathbb{M}(I)<+\infty$ and $\ast d F=\partial I$. By definition of integral $1$-current, there exist a locally $1$-rectifiable set $\Gamma\subset Q_1(0)$, a Borel measurable unitary vector field $\vec I$ on $\Gamma$ and a positive $\z$-valued $\H^1\res\Gamma$-integrable function $\theta\in L^1(\Gamma,\H^1)$ such that
\begin{align*}
    \langle I,\omega\rangle=\int_{\Gamma}\theta\langle\omega,\vec I\rangle\, d\H^1, \qquad\forall\,\omega\in\D^1(Q_1(0)).
\end{align*}
By the coarea formula, there exists $G\subset K$ such that $\L^1(K\smallsetminus G)=0$ and such that $\Gamma\cap\partial Q_{\rho}(x_0)$ is a finite set for every $\rho\in G$.

Let $\Psi:\s^{n-1}\to\partial Q_1(0)$ be the bi-Lipschitz map given by 
\begin{align*}
    \Psi(x):=\frac{x}{2\lVert x\rVert_{\infty}}, \qquad\forall\, x\in\s^{n-1}.
\end{align*}
Consider the map $\Phi:\s^{n-1}\times[0,\rho_0]\rightarrow\im(\Phi)=\overline{Q_{\rho_0}(x_0)}\subset \overline{Q_1(0)}$ given by
\begin{align*}
    \Phi(y,t):=x_0+t\Psi(y), \qquad\forall\,(y,t)\in\s^{n-1}\times[0,\rho_0].
\end{align*}
Notice that $\Phi\big\vert_{\s^{n-1}\times [\rho_1, \rho_2]}$ is a bi-Lipschitz homeomorphism onto its image for every $\rho_1,\rho_2\in (0,1)$.

We claim that estimate \eqref{holderianity estimate} holds on a full-measure subset of $G$. Indeed, fix any $\rho_1,\rho_2\in G$. Without loss of generality, assume that $\rho_2>\rho_1$. Let $\hat\Phi:=\Phi\big\vert_{\s^{n-1}\times [\rho_1,\rho_2]}$.
Define $\pi:=\text{pr}_1\circ\Phi^{-1}:\overline{Q_{\rho_0}(x_0)}\rightarrow\s^{n-1}$, where $\text{pr}_1:\s^{n-1}\times [0,\rho_0]\rightarrow\s^{n-1}$ is the canonical projection on the first factor, and notice that $\pi$ is a Lipschitz and proper map. Then, $\pi_*\big(T_F\res\im(\hat\Phi)\big)\in\D_1(\s^{n-1})$ can be expressed as follows: for any $\omega\in\Omega^1(\s^{n-1})$
\begin{align*}
    \langle\pi_*\big(T_F\res\im(\hat\Phi)\big),\omega\rangle&=\langle T_F\res\im(\hat\Phi),\pi^*\omega\rangle=\int_{\im(\hat\Phi)}F\wedge\pi^*\omega\\
    &=\int_{\im(\hat\Phi)}(\Phi^{-1})^*(\Phi^*F\wedge \Phi^*\pi^*\omega)\\
    &=\int_{\s^{n-1}\times [\rho_1,\rho_2]}\Phi^*F\wedge\text{pr}_1^*\omega\\
    &=\int_{\s^{n-1}\times [\rho_1,\rho_2]}\text{pr}_1^*\omega\wedge *(*\Phi^*F)\\
    &=\int_{\s^{n-1}\times [\rho_1,\rho_2]}\langle\text{pr}_1^*\omega,*\Phi^*F\rangle\, d\vol_{\s^{n-1}\times [\rho_1,\rho_2]}(y,t)\\
    &=\int_{\rho_1}^{\rho_2}\bigg(\int_{\s^{n-1}}\langle\text{pr}_1^*\omega,*\Phi^*F\rangle\, d\vol_{\s^{n-1}}(y)\bigg)\, dt\\
    &=\int_{\rho_1}^{\rho_2}\bigg(\int_{\s^{n-1}}\langle\omega,i_{\s^{n-1}\times\{t\}}^**\Phi^*F\rangle\, d\vol_{\s^{n-1}}(y)\bigg)\, dt\\
    &=\int_{\rho_1}^{\rho_2}\bigg(\int_{\s^{n-1}}\omega\wedge\big(*i_{\s^{n-1}\times\{t\}}^**\Phi^*F\big)\bigg)\, dt\\
    &=\int_{\s^{n-1}}\omega\wedge\bigg(\int_{\rho_1}^{\rho_2}\big(*i_{\s^{n-1}\times\{t\}}^**\Phi^*F\big)\, dt\bigg)\\
    &=(-1)^{n-1}\int_{\s^{n-1}}\omega\wedge\alpha,
\end{align*}
where 
\begin{align*}
    \alpha:=(-1)^{n-1}\int_{\rho_1}^{\rho_2}\big(*i_{\s^{n-1}\times\{t\}}^**\Phi^*F\big)\, dt\in\Omega_p^{n-2}(\s^{n-1}).
\end{align*}
In particular,
\begin{align*}
    \langle\partial\pi_*\big(T_F\res\im(\hat\Phi)\big),\varphi\rangle&=\langle\pi_*\big(T_F\res\im(\hat\Phi)\big),d\varphi\rangle\\
    &=(-1)^{n-1}\int_{\s^{n-1}}d\varphi\wedge\alpha=\langle d^*(*\alpha),\varphi\rangle, \qquad\forall\,\varphi\in C^{\infty}(\s^{n-1}).
\end{align*}
Recall that the restriction of an integral current to a measurable set is still an integral current. Moreover, the push-forward of an integral current through a Lipschitz and proper map remains an integral current (see \cite[Chapter 7, \S 7.5]{krantz_parks-geometric_integration_theory}). Then, $\tilde I:=-\pi_*\big(I\res\im(\hat\Phi)\big)\in\R_1(\s^{n-1})$.\\
So far, we have shown that
\begin{align*}
     \partial\pi_*\big((T_F-I)\res\im(\hat\Phi)\big)=\partial\pi_*\big(T_F\res\im(\hat\Phi)\big)-\partial\pi_*\big(I\res\im(\hat\Phi)\big)=d^*(*\alpha)+\partial\tilde I.
\end{align*}
Let $\zeta\in C^\infty_c((-1,1))$ such that $\displaystyle{\int_\mathbb{R}\zeta=1}$. For any $\varepsilon\in (0,\min\{\rho_1,\rho_0-\rho_2\}) $ set $\displaystyle{\zeta_\varepsilon=\frac{1}{\varepsilon}\zeta\left(\frac{\cdot}{\varepsilon}\right)}$ and let $\rchi_\varepsilon$ be the unique solution of
\begin{align*}
    \begin{cases}
        \rchi_\varepsilon'(x)=\zeta_\varepsilon(x-\rho_1)-\zeta_\varepsilon(x-\rho_2)\\
        \rchi_\varepsilon(0)=0.
    \end{cases}
\end{align*} 
Let $\psi\in C^\infty(\s^{n-1}\times[0,\rho_0])$ and let $\text{pr}_2:\s^{n-1}\times[0,\rho_0]\to[0,\rho_0]$ be the projection on the second factor.
We compute
\begin{align*}
    \langle(\Phi^{-1})_\ast I, \psi\,d(\rchi_\varepsilon\circ \text{pr}_2)\rangle&=\int_{\Phi^{-1}(\Gamma)}\theta_{(\Phi^{-1})_\ast I}\psi(\rchi_\varepsilon'\circ \text{pr}_2)\langle d \text{pr}_2, \vec I_{(\Phi^{-1})_\ast I}\rangle\,d\mathscr{H}^1\\
    &=\int_{\rho_1-\varepsilon}^{\rho_1+\varepsilon}\zeta_\varepsilon(t)\bigg(\int_{\Phi^{-1}(\Gamma)\cap(\s^{n-1}\times\{t\})}\psi \tilde \theta\,d\mathscr{H}^0\bigg)\,d\L^1(t)\\
    &\quad-\int_{\rho_2-\varepsilon}^{\rho_2+\varepsilon}\zeta_\varepsilon(t)\bigg(\int_{\Phi^{-1}(\Gamma)\cap(\s^{n-1}\times\{t\})}\psi \tilde \theta\,d\mathscr{H}^0\bigg)\,d\L^1(t)
\end{align*}
with $\tilde{\theta}=\theta_{(\Phi^{-1})_\ast I}\text{sgn}(\langle d \text{pr}_2, \vec I_{(\Phi_{-1})_\ast I}\rangle)\in L^1(\Phi^{-1}(\Gamma),\z)$.
Moreover
\begin{align*}
    \langle (\Phi^{-1})_\ast T_F, \psi\,d(\rchi_\varepsilon\circ \text{pr}_2)\rangle&=\int_{\s^{n-1}\times [0,\rho_0]}\psi(\rchi_\varepsilon'\circ\text{pr}_2){\Phi}^\ast F\wedge\,d\text{pr}_2\\
    &=\int_{\rho_1-\varepsilon}^{\rho_1+\varepsilon}\zeta_\varepsilon(t)\left(\int_{\s^{n-1}\times\{t\}}\psi {(\Phi\big\vert_{\s^{n-1}\times\{t\}})}^\ast F\right)\,d\L^1(t)\\
    &\quad-\int_{\rho_2-\varepsilon}^{\rho_2+\varepsilon}\zeta_\varepsilon(t)\left(\int_{\s^{n-1}\times\{t\}}\psi {(\Phi\big\vert_{\s^{n-1}\times\{t\}})}^\ast F\right)\,d\L^1(t).
\end{align*}
Now observe that
\begin{align*}
    \langle (\Phi^{-1})_\ast (T_F-I), (\rchi_\varepsilon\circ \text{pr}_2)\,d\psi\rangle\to\langle ((\Phi^{-1})_\ast (T_F-I))\res (\s^{n-1}\times[\rho_1,\rho_2]), d\psi\rangle
\end{align*}
as $\varepsilon\to 0^+$, by dominated convergence. On the other hand, since $\partial (T_F-I)=0$, we have
\begin{align*}
    \langle(\Phi^{-1})_\ast(T_F-I),(\rchi_\varepsilon\circ\text{pr}_2)\,d\psi\rangle=\langle (\Phi^{-1})_\ast I,\psi\,d(\rchi_\varepsilon\circ\text{pr}_2)\rangle-\langle (\Phi^{-1})_\ast T_F,\psi\,d(\rchi_\varepsilon\circ\text{pr}_2)\rangle.
\end{align*}
Therefore for almost every $\rho_1,\rho_2\in (0,\rho_0)$ (depending on $\psi$) we have
\begin{align}\label{Equation: boundary of pushforward}
    \nonumber
    \langle \partial (({\Phi^{-1}})_\ast (T_F-I)\res (\s^{n-1}\times[\rho_1,\rho_2])),\psi\rangle&=\int_{\Phi^{-1}(\Gamma)\cap(\s^{n-1}\times\{\rho_1\})}\psi \tilde \theta\,d\mathscr{H}^0\\
    &\quad-\int_{\Phi^{-1}(\Gamma)\cap(\s^{n-1}\times\{\rho_2\})}\psi \tilde \theta\,d\mathscr{H}^0\\
    \nonumber
    &\quad-\int_{\s^{n-1}\times\{\rho_1\}}\psi {(\Phi\big\vert_{\s^{n-1}\times\{\rho_1\}})}^\ast F\\
    \nonumber
    &\quad+\int_{\s^{n-1}\times\{\rho_2\}}\psi {(\Phi\big\vert_{\s^{n-1}\times\{\rho_2\}})}^\ast F.
\end{align}
Now let $\{\psi_k\}_{k\in \mathbb{N}}\subset C^\infty(\s^{n-1}\times[0,\rho_0])$ be a countable sequence dense in $C^1(\s^{n-1}\times [0,\rho_0])$. For every $k\in\n$, let $E_k\subset G$ be the set such that \eqref{Equation: boundary of pushforward} holds with $\psi=\psi_k$ (i.e. the set of the $\rho\in G$ wich are "$\zeta_\varepsilon$-Lebesgue points" of the integrands in \eqref{Equation: boundary of pushforward}, with $\psi=\psi_k$) and define 
\begin{align*}
    E:=\bigcap_{k\in\n} E_k.
\end{align*}
Then $\L^1(E)=\L^1(K)$ and for every $\rho_1,\rho_2\in E$ estimate \eqref{Equation: boundary of pushforward} holds with $\psi=\psi_k$ for every $k\in\n$. By density of $\{\psi_k\}_{k\in \mathbb{N}}$ in $C^1(\s^{n-1}\times [0,\rho_0])$, we can pass to the limit in \eqref{Equation: boundary of pushforward} and get that for any given couple of parameters $\rho_1,\rho_2\in\tilde E$ such estimate holds for every $\psi\in C^{\infty}(\s^{n-1}\times [0,\rho_0])$.
In particular, for every $\rho_1, \rho_2\in E$, $\varphi\in C^\infty(\s^{n-1})$ we have
\begin{align*}
    \langle \partial\pi_\ast ((T_F-I)\res{\text{Im}(\hat{\Phi})}),\varphi\rangle&=\sum_{x\in \Gamma_{\rho_1}}\tilde \theta(x,\rho_1)\varphi(x)-\sum_{x\in \Gamma_{\rho_2}}\tilde \theta(x,\rho_2)\varphi(x)\\
    &\quad-\int_{\s^{n-1}} \varphi\,\Psi^\ast s(\rho_1)+\int_{\s^{n-1}}\varphi\,\Psi^\ast s(\rho_2),
\end{align*}
where
\begin{align*}
    \Gamma_{\rho_1}&:=\text{pr}_1(\Phi^{-1}(\Gamma)\cap(\s^{n-1}\times\{\rho_1\}))\subset\s^{n-1},\\
    \Gamma_{\rho_2}&:=\text{pr}_1(\Phi^{-1}(\Gamma)\cap(\s^{n-1}\times\{\rho_2\}))\subset \s^{n-1}
\end{align*}
are finite set for any $\rho_1$, $\rho_2\in G$.

Gathering together what we have proved so far, we have 
\begin{align*}
    \ast\big(\Psi^*s(\rho_2)-\Psi^*s(\rho_1)\big)=d^*(*\alpha)+\partial I'+\bigg(\int_{\s^{n-1}}\Psi^*s(\rho_2)-\Psi^*s(\rho_1)\bigg)\delta_q,
\end{align*}
where $I'\in\R_1(\s^{n-1})$ is any rectifiable one-current of finite mass such that
\begin{align*}
    \partial I'= \sum_{x\in \Gamma_{\rho_2}}\tilde \theta(x,\rho_2)\delta_x-\sum_{x\in \Gamma_{\rho_1}}\tilde \theta(x,\rho_1)\delta_x+\partial \tilde I+ \left(\sum_{x\in\Gamma_{\rho_1}}\tilde \theta(x,\rho_1)-\sum_{x\in\Gamma_{\rho_2}}\tilde \theta(x,\rho_2)\right)\delta_q,
\end{align*}
i.e. $*\alpha$ is a competitor in the definition of $d\big(\Psi^*s(\rho_2),\Psi^*s(\rho_1)\big)$. Hence, in order to estimate $d\big(\Psi^*s(\rho_2),\Psi^*s(\rho_1)\big)$ we just need to find an upper bound for $\lvert\lvert*\alpha\rvert\rvert_{L^p}$.\\
Notice that $\lvert d\Phi\rvert\leq t\lvert d\Psi\rvert+\frac{\sqrt{n}}{2}$. Moreover since \begin{align*}
    \Phi^{-1}(x)=\left(\Psi^{-1}\left(\frac{x-x_0}{2\lVert x-x_0\rVert_\infty}\right), 2\lVert x-x_0\rVert_\infty\right)
\end{align*}
we have $\lvert J\Phi^{-1}(x)\rvert\leq 2^{-(n-2)}\left(\frac{\lVert d\Psi^{-1}\rVert_{L^\infty}}{\lVert x-x_0\rVert_\infty}\right)^{n-1}$.
Therefore
\begin{align*}
    \lVert \ast \alpha\rVert_{L^p}^p&\leq\int_{\s^{n-1}}\left\lvert\int_{\rho_1}^{\rho_2}\lvert \Phi^\ast F \rvert dt\right\rvert^pd\mathscr{H}^{n-1}\leq \lvert \rho_1-\rho_2\rvert^\frac{p}{p'}\int_{\s^{n-1}\times [\rho_1,\rho_2]}\lvert\Phi^\ast F\rvert^pd\mathscr{H}^{n-1}dt\\
    &\leq \left(2^{-(n-2)}\left(\lVert d\Psi\rVert_{L^\infty}+\frac{\sqrt{n}}{2}\right)^{p(n-1)}\frac{\lVert d\Psi^{-1}\rVert_{L^\infty}^{n-1}}{\rho_1^{n-1}}\right)\lvert \rho_1-\rho_2\rvert^\frac{p}{p'}\lVert F\rVert^p_{L^p(Q_1(0))}
\end{align*}
and our claim follows.

\end{proof}
\end{prop}
\subsection{Proof of Theorem \ref{Theorem: weak closure} for $Q^n_1(0)$}
For the proof of Theorem \ref{Theorem: weak closure} we need two technical Lemmata.
\begin{lem}
\label{useful lemma 1}
Let $\{f_k\}_{k\in\n}\subset L^1(0,1)$ be such that $\lvert\lvert f_k\rvert\rvert_{L^1}\le C$ for any $k\in \mathbb{N}$. Then there exist a sequence of compact subsets $\{W_h\}_{h\in\n_{\geq 2}}$ of $(0,1)$ such that for every $h\in\n_{\geq 2}$ the following properties hold:
\begin{enumerate}
    \item $\displaystyle{L^1(W_h)= 1-\frac{C+2}{h}}$;
    \item $W_h\subset (1/h,1)$;
    \item for almost every $\rho\in W_h$ and every $k\in\n$ there exists $k'>k$ such that $|f_{k'}(\rho)|\le h$. 
\end{enumerate}
\begin{proof}
Let $h\in\n_{\geq 2}$.
For any $l\in \mathbb{N}$ let
\begin{align*}
    A^h_l:=\bigcap_{k=l}^\infty f_k^{-1}([-h,h]^c).
\end{align*}
Notice that for any $l\in \mathbb{N}$ $A^h_l\subset A_{l+1}^h$ and set
\begin{align*}
    I_h=\bigcup_{l=1}^\infty A_l^h.
\end{align*}
Let $m\in\n$ and let $k\ge m$. Notice that
\begin{align*}
    C\ge\int_0^1|f_k(\rho)|\, d\rho\ge\int_{A_m^h}|f_k(\rho)|\, d\rho>h\L^1(A_m^h).
\end{align*}
By letting $m\rightarrow\infty$ in the previous inequality, be obtain 
\begin{align*}
    L^1(I_h)\leq\frac{C}{h}. 
\end{align*}
Then, by defining $E_h:=I_h^c\cap(1/h,1)$, we clearly get 
\begin{align*}
    \L^1(E_h)=\L^1((I_h\cup (0,1/h])^c)=1-\L^1(I_h\cup (0,1/h])\geq 1-\frac{C+1}{h}.
\end{align*}
Moreover, $E_h$ and any of its subsets satisfy the properties $2.$ and $3.$. 
Finally, since $E_h$ is measurable, we can find a compact set $W_h\subset E_h$ such that 
\begin{align*}
    \L^1(W_h)=1-\frac{C+2}{h}.
\end{align*}
By construction, $W_h$ satisfies 1, 2 and 3. 
\end{proof}
\end{lem}
Since we do not know if the space $(Y, d)$ is complete, we will also need the following lemma.
\begin{lem}
\label{useful lemma 2}
Let $K\subset [0,1]$ be compact and let $S\subset K$ be dense and countable. Let $\{f_k\}_{k\in\n}\subset C^0(K,Y)$ be such that 
\begin{enumerate}
    \item $\{f_k\}_{k\in\n}$ is uniformly Cauchy  from $K$ to $(Y,d)$ (i.e. it is a Cauchy sequence w.r.t. uniform convergence);
    \item for some $C>0$, 
        \begin{align*}
        \sup_{\rho\in S}\sup_{k\in\n}||f_k(\rho)||_{L^p}\le C;
        \end{align*}
    \item for some $A>0$ and some $\alpha\in(0,1]$, we have 
        \begin{align*}
            d\big(f_k(\rho),f_k(\rho')\big)\le A|\rho-\rho'|^{\alpha}, \qquad\forall\,\rho,\rho'\in S.
        \end{align*}
\end{enumerate}
Then, there exists $f\in C^0(K,Y)$ such that $f_k\rightarrow f$ uniformly.  
\begin{proof}
Fix any $\rho\in S$. By hypothesis 2, $\{f_k(\rho)\}_{k\in\n}$ is bounded in $L^p$ and therefore it has a subsequence converging weakly in $L^p$ to a limit $f(\rho)\in Y$ (recall that $Y$ is closed with respect to the weak $L^p$ convergence). By Corollary \ref{d_Psi metrizes the weak topology on bounded subsets}, such a subsequence converges in $(Y,d)$ to the same limit $f(\rho)$. Since $\{f_k\}_{k\in\n}$ is uniformly Cauchy from $K$ to $(Y,d)$, we have that $\{f_k(\rho)\}_{k\in\n}$ is Cauchy in $(Y,d)$, therefore $f_k(\rho)\xrightarrow{d}f(\rho)$ and, by Corollary \ref{d_Psi metrizes the weak topology on bounded subsets}, $f_k(\rho)\rightharpoonup(\rho)$ weakly in $L^p$.

Fix any $\rho\in K$ and let $\{\rho_i\}_{i\in\n}\subset S$ be such that $\rho_i\rightarrow\rho$. We claim that there exists $f_\rho\in Y$ such that $f(\rho_i)\to f_\rho$ w.r.t. $d$ and $f_\rho$ doesn't depend on the choice of the sequence $\{\rho_i\}_{i\in\n}$. Indeed 
by lower semicontinuity of the $L^p$ norm w.r.t. the weak $L^p$ convergence, we have 
\begin{align*}
    ||f(\rho_i)||_{L^p}&\le\liminf_{k\rightarrow\infty}||f_k(\rho_i)||_{L^p}\le C,
\end{align*}
for every $i\in\n$.
Then there exists a subsequence $\{\rho_{i_j}\}_{j\in \mathbb{N}}$ and a $f_\rho\in L^p$ such that $f(\rho_{i_j})\rightharpoonup f_\rho$ weakly in $L^p$. Since $Y$ is weakly closed in $L^p$ we have $f_\rho\in Y$. By Corollary \ref{d_Psi metrizes the weak topology on bounded subsets} we also have $f(\rho_{i_j}) \to f_\rho$ w.r.t. $d$.\\
Relabel the subsequence as $\{\rho_i\}_{i\in \mathbb{N}}$. To see that $f_\rho$ doesn't depend on the subsequence (and thus it doesn't depend on $\{\rho_i\}_{i\in \mathbb{N}}$), assume that $\{\tilde \rho_i\}_{i\in \mathbb{N}}$ is another sequence in $S$ with $\tilde \rho_i\to \rho$ and $f(\tilde \rho_i)\to \tilde f_{\rho}$ w.r.t. $d$.
To see that $f_\rho=\tilde f_{ \rho}$, first notice that by hypothesis $3.$ and by triangle inequality we have
\begin{align*}
    d\big(f(\rho_i),f(\tilde \rho_i)\big)&\le d\big(f(\rho_i),f_k(\rho_i)\big)+d\big(f_k(\rho_i),f_k(\tilde \rho_i)\big)+d\big(f_k(\tilde \rho_i),f(\tilde \rho_i)\big)\\
    &\le d\big(f(\rho_i),f_k(\rho_i)\big)+A|\rho_i-\tilde \rho_i|^{\alpha}+d\big(f_k(\tilde \rho_i),f(\tilde \rho_i)\big).
\end{align*}
for every $i\in\n$. Hence, passing to the limit as $k\rightarrow\infty$ in the previous inequality, we get 
\begin{align*}
    d\big(f(\rho_i),f(\tilde \rho_i)\big)\le A|\rho_i-\tilde \rho_i|^{\alpha}.
\end{align*}
Thus we finally obtain 
\begin{align*}
    d(f_\rho,\tilde f_{\rho})&\le d\big(f_\rho,f(\rho_i)\big)+d\big(f(\rho_i),f(\tilde \rho_i)\big)+d\big(f(\tilde \rho_i),\tilde f_{\rho}\big)\\
    &\le d\big(f_\rho,f(\rho_i)\big)+A|\rho_i-\tilde \rho_i|^{\alpha}+d\big(f(\tilde \rho_i),\tilde f_{\rho}\big)
\end{align*}
and, passing to the limit as $i\rightarrow\infty$, we get $d(f_\rho,\tilde f_{ \rho})=0$, i.e. $f_\rho=\tilde f_{ \rho}$.\\
For any $\rho\in K$ let
\begin{align*}
    f(\rho):=\lim_{i\rightarrow\infty}f(\rho_i),
\end{align*}
where $\{\rho_i\}_{i\in\n}\subset S$ is any sequence such that $\rho_i\rightarrow\rho$ and the limit is understood w.r.t. $d$. Then $f$ is a well-defined function on $K$.\\
To see that $f_k\to f$ uniformly, let $\varepsilon>0$ and let $K\in \mathbb{N}$ such that for any $m,n\in \mathbb{N}_{\geq K}$
\begin{align*}
    \sup_{\rho\in S} d(f_m(\rho), f_n(\rho))<\varepsilon.
\end{align*}
Let $\rho\in K$ and let $\{\rho_i\}_{i\in \mathbb{N}}$ be a sequence in $S$ such that $\rho_i\to \rho$.
Then for any $k\in \mathbb{N}_{\geq K}$ we have
\begin{align*}
    d(f_k(\rho), f(\rho))=\lim_{i\to\infty} d(f_k(\rho_i), f(\rho_i))=\lim_{i\to\infty} \lim_{m\to\infty} d(f_k(\rho_i), f_m(\rho_i))\leq \varepsilon.
\end{align*}
This shows that $f$ is the uniform limit of $\{f_k\}_{k\in \mathbb{N}}$ in $K$, with respect to $d$. As $f_k\in C^0(K,Y)$ for any $k\in \mathbb{N}$, we have that $f\in C^0(K,Y)$ (this also follows directly from the construction of $f$).
\end{proof}
\end{lem}
\begin{thm}[Weak closure for $Q_1^n(0)$]\label{weak closure for Q^n}
Fix any $n\in\n_{\geq 2}$ and assume that $p\in\big(1,n/(n-1)\big)$. Then $\Omega_{p,\z}^{n-1}(Q_1^n(0))$ is weakly sequentially closed.
\begin{proof}
Assume that $F\in\Omega_p^{n-1}(Q_1(0))$ belongs to the weak $L^p$-closure of $\Omega_{p,\z}^{n-1}(Q_1(0))$, i.e. there exists $\{F_k\}_{k\in\n}\subset\Omega_{p,R}^{n-1}(Q_1(0))$ such that $F_k\xrightharpoonup{L^p} F$. What we need to show is that $F\in\Omega_{p,\z}^{n-1}(Q_1(0))$, which amounts to saying that 
\begin{align}
\label{integer fluxes condition}
    \int_{\partial Q_{\rho}(x_0)}i_{\partial Q_{\rho}(x_0)}^*F\in\z, 
\end{align}
for every $x_0\in Q_1(0)$ and  for a.e. $\rho\in\big(0,2\dist_{\infty}(x_0,\partial Q_1(0))\big)$. Without losing generality, we will just show \eqref{integer fluxes condition} for $x_0=0$. 

\textbf{Step 1}. For any $k\in \mathbb{N}$ let $s_k$ be the slice function of $F_k$ at $0$. Fix any $h\in\n_{\geq 2}$ and let $W_h\subset (1/h,1)$ be the compact set given by applying Lemma \ref{useful lemma 1} with $f_k=\lVert s_k\rVert_{L^p}$ and $C=2^\frac{n}{p}\sup_{k\in \mathbb{N}}\lVert F_k\rVert_{L^p}$ (see Remark \ref{Remark: Lp norm of slice function}).
Let $E_h^k\subset W_h$ denote the subset associated to $W_h$ and $s_k$ by Proposition \ref{the slice functions are (1-1/p)-holder continuous}, let $\tilde E_h=\bigcap_{k\in\mathbb{N}}E_h^k$ and let $s_k$ denote its $\frac{1}{p'}$-H\"older representative on $\tilde{E}_h$, for any $k\in \mathbb{N}$.
By property $3.$ in Lemma \ref{useful lemma 1}, for almost every $\rho\in \tilde E_h$ we can find a subsequence $\{s_{k_{\rho}}(\rho)\}_{k_{\rho}\in\n}\subset\{s_k(\rho)\}_{k\in\n}$ such that $s_{k_{\rho}}(\rho)$ is uniformly bounded in $L^p$ by $h$. Denote by $E_h$ the set of all such $\rho$ and observe that $\L^1({E_h})=\L^1(W_h)=1-\frac{C+2}{h}$. Then for any $\rho\in  E_h$, $\{s_{k_{\rho}}(\rho)\}_{k_{\rho}\in\n}$ has a subsequence that converges weakly in $L^p$ or, equivalently, with respect to $d_\Psi$ (see Corollary \ref{d_Psi metrizes the weak topology on bounded subsets}).\\
Let $S_h\subset E_h$ be a countable dense subset. By a diagonal extraction argument, we find a subsequence $\{s_{k_l}\}_{l\in\n}$ such that $\{s_{k_l}(\rho)\}_{l\in\n}$ is convergent with respect to $d_\Psi$ and weakly in $L^p$, and is uniformly bounded in $L^p$ by $h$, for every $\rho\in S_h$. 

\textbf{Step 2}. Next we claim that for any $l\in \mathbb{N}$, $s_{k_l}$ can be extended to a $\frac{1}{p'}$-H\"older continuous function on $\overline{E_h}$ (with the same H\"older constant, which is bounded uniformly in $l$).\\
In fact let $f\in \{s_{k_l}\}_{l\in \mathbb{N}}$, let $\rho\in \overline{E_h}\smallsetminus E_h$, let $\{\rho_i\}_{i\in \mathbb{N}}$ be a sequence in $S_h$ such that $\rho_i\to \rho$ as $i\to \infty$ (observe that such a sequence exists, since $S_h$ is dense in $E_h$, which in turn is dense in $\overline{E_h}$).
Since $\{f(\rho_i)\}_{i\in\n}\subset Y$ is uniformly bounded  in $L^p$, by weak $L^p$-compactness there exists a subsequence $\{f(\rho_{i_j})\}_{j\in\n}$ of $\{f(\rho_i)\}_{i\in\n}$ such that $f(\rho_{i_j})\rightharpoonup f_{\rho}$ weakly in $L^p$ for some $f_{\rho}\in Y$.  By Corollary \ref{d_Psi metrizes the weak topology on bounded subsets}, we know that $d_{\Psi}\big(f(\rho_{i_j}),f_{\rho}\big)\rightarrow 0$.
Since $f$ is $\frac{1}{p'}$-H\"older continuous on $\overline {E_h}$, $f_\rho$ does not depend on the sequence $\{\rho_i\}_{i\in \mathbb{N}}$.
Hence, the function 
\begin{align*}
  	\tilde f(\rho):=\begin{cases}f(\rho) & \mbox{ if } \rho\in E_h\\ f_{\rho} & \mbox{ if }\rho\in \overline{E_h}\smallsetminus E_h\end{cases}
  \end{align*} 
is well-defined on $\overline{E_h}$ and satisfies \eqref{holderianity estimate} on $\overline{E_h}$.\\
In the following, in order to simplify the notation, we will denote again by $s_{k_l}$ the $\frac{1}{p'}$-H\"older extension of $s_{k_l}$ to $\overline{E_h}$, for any $l\in \mathbb{N}$.

\textbf{Step 3}. We show that $\{s_{k_l}\}_{l\in\n}$ converges uniformly on $\overline{E_h}$ to some $s\in C^0(\overline{E_h},Y)$.\\
Fix any $\varepsilon>0$. By 
Step 2. we know that the sequence $\{s_{k_l}\}_{l\in\n}$ is equicontinuous from $\overline{E_h}$ to $(Y,d_{\Psi})$. Therefore we can choose $\delta>0$ such that 
\begin{align*}
    d_{\Psi}\big(s_{k_l}(\rho),s_{k_l}(\rho'))<\varepsilon, \qquad\forall\,\rho,\rho'\in \overline{E_h} \mbox{ s.t. } \lvert\rho-\rho'\rvert<\delta \mbox{ and } \forall\,l\in\n.
\end{align*}
Notice that $\{(\rho-\delta,\rho+\delta)\}_{\rho\in S_h}$ is an open cover of $\overline{E_h}$. Since $\overline{E_h}$ is compact, we can find a finite set $\{\rho_1,...,\rho_m\}\subset S_h$ such that $\{(\rho_j-\delta,\rho_j+\delta)\}_{j=1,...,m}$ is a finite open cover of $\overline{E_h}$. 

Now let $\rho\in \overline{E_h}$. Observe that there exists a point $\rho_j\in \{\rho_1,...,\rho_m\}$ such that $\rho\in (\rho_j-\delta,\rho_j+\delta)$, i.e. $|\rho-\rho_j|<\delta$. By our choice of $\delta$, this implies $d\big(s_{k_l}(\rho),s_{k_l}(\rho_j))<\varepsilon$, for every $l\in\n$. By triangle inequality, we have 
\begin{align*}
    d_{\Psi}\big(s_{k_l}(\rho),s_{k_m}(\rho))&\le d_{\Psi}\big(s_{k_l}(\rho),s_{k_l}(\rho_j))+d_{\Psi}\big(s_{k_l}(\rho_j),s_{k_m}(\rho_j))+d_{\Psi}\big(s_{k_m}(\rho_j),s_{k_m}(\rho))\\
    &<2\varepsilon+d_{\Psi}\big(s_{k_l}(\rho_j),s_{k_m}(\rho_j)).
\end{align*}
But since $\rho_j\in S_h$, we know that there exists $L_j>0$ such that 
\begin{align*}
    d_{\Psi}\big(s_{k_l}(\rho_j),s_{k_m}(\rho_j))<\varepsilon, \qquad\forall\, l,m\ge L_j.
\end{align*}
Hence, by letting $\displaystyle{L:=\max_{j=1,...,m}{L_j}}$, we have that 
\begin{align*}
    d_{\Psi}\big(s_{k_l}(\rho),s_{k_m}(\rho))<3\varepsilon, \qquad\forall\, l,m\ge L, \forall\,\rho\in \overline{E_h}. 
\end{align*}
Here we have just proved that the sequence $\{s_{k_l}\}_{l\in\n}$ is uniformly Cauchy on $\overline{E_h}$. Since $\{s_{k_l}\}_{l\in\n}$ satisfies all the hypotheses of Lemma \ref{useful lemma 2}, we get that there exists $s\in C^0(\overline{E_h},Y)$ such that $s_{k_l}\rightarrow s$ uniformly on $\overline{E_h}$ w.r.t. $d_\Psi$.\\
Notice that since $\lVert s_{k_l}(\rho)\rVert_{L^p}\leq h$ for any $l\in \mathbb{N}$, for any $\rho\in E_h$, and since $s_{k_l}(\rho)\to s(\rho)$ w.r.t. $d_\Psi$ for any $\rho\in \overline{E_h}$, by Corollary \ref{d_Psi metrizes the weak topology on bounded subsets} there holds $s_{k_l}(\rho)\rightharpoonup s(\rho)$ in $L^p$ for any $\rho\in  E_h$.

\textbf{Step 4}. Let $s_0:E_h\rightarrow Y$ be the restriction to $E_h$ of the slice function of $F$ at $0$. We claim that $s=s_0$ a.e. in $E_h$. To show this we will prove that
\begin{align*}
    \int_{E_h}\varphi(\rho)\int_{\partial Q_1(0)}\psi\big(s_{k_l}(\rho)-s_0(\rho)\big)\, d\rho\rightarrow 0, \qquad\mbox{ as } l\rightarrow\infty, 
\end{align*}
for every $\varphi\in C^{\infty}_c((0,1))$ and for every $\psi\in\operatorname{Lip}(\partial Q_1(0))$. Indeed, an explicit computation gives 
\begin{align*}
    \int_{E_h}\varphi(\rho)&\int_{\partial Q_1(0)}\psi\big(s_{k_l}(\rho)-s_0(\rho)\big)\, d\rho\\
    &=\int_{E_h}\varphi(\rho)\bigg(\int_{\partial Q_{\rho}(0)}\psi\bigg(\frac{\cdot}{\rho}\bigg)i_{\partial Q_{\rho}(0)}^*F_{k_l}-\int_{\partial Q_{\rho}}\psi\bigg(\frac{\cdot}{\rho}\bigg)i_{\partial Q_{\rho}(0)}^*F\bigg)\, d\rho\\
    &=2^n\int_{Q_1(0)}\mathbbm{1}_{E_h}(2r)\varphi(2|\cdot|)\psi\bigg(\frac{\cdot}{\lvert\,\cdot\,\rvert}\bigg)dr\wedge (F_{k_l}-F)\rightarrow 0\text{ as } l\to\infty,
\end{align*}
since 
\begin{align*}
   \mathbbm{1}_{E_h}(2r)\varphi(2|\cdot|)\psi\bigg(\frac{\cdot}{\lvert\,\cdot\,\rvert}\bigg)dr\in\Omega_{p'}^1(Q_1(0))
\end{align*}
and $F_{k_l}\xrightharpoonup{L^p} F$.
On the other hand, since $s_{k_l}(\rho)\rightharpoonup s(\rho)$ in $L^p$ for any $\rho\in E_h$,
\begin{align*}
    \int_{E_h}\varphi(\rho)\int_{\partial Q_1(0)}\psi\big(s_{k_l}(\rho)-s(\rho)\big)\, d\rho\rightarrow 0 \qquad\mbox{ as } l\rightarrow\infty
\end{align*}
for every $\varphi\in C^{\infty}_c((0,1))$ and for every $\psi\in\operatorname{Lip}(\partial Q_1(0))$, we obtain
\begin{align*}
    \int_{E_h}\varphi(\rho)\int_{\partial Q_1(0)}\psi\big(s(\rho)-s_0(\rho)\big)\, d\rho=0, \qquad\forall\,\varphi\in C^{\infty}_c((0,1)),\forall\,\psi\in\operatorname{Lip}(\partial Q_1(0)).
\end{align*}
This means that $s_0(\rho)=s(\rho)\in Y$ for a.e. $\rho\in E_h$.

\textbf{Step 5}. Finally we show that \eqref{integer fluxes condition} holds for almost any $\rho\in (0,1)$.\\
In fact for any $\rho\in E_h$ such that $s_0(\rho)=s(\rho)$ we have 
\begin{align*}
    \int_{\partial Q_\rho(0)}i^\ast_{\partial Q_\rho}F=\int_{\partial Q_1(0)}s_0(\rho)=\int_{\partial Q_1(0)}s(\rho)=\lim_{k_l\to\infty}\int_{\partial Q_1(0)}s_{k_l}(\rho)\in \mathbb{Z},
\end{align*}
since $s_{k_l}(\rho)\rightharpoonup s(\rho)$ in $L^p$. Thus \eqref{integer fluxes condition} holds for $\L^1$-a.e. $\rho\in E_h$. Since the previous step can be repeated for any $h\in \mathbb{N}_{\geq 2}$, and since
\begin{align*}
    \lim_{h\to +\infty}\L^1(E_h)=\lim_{h\to +\infty} 1-\frac{C+1}{h}=1,
\end{align*}
we conclude that \eqref{integer fluxes condition} holds for $\L^1$-a.e. $\rho\in (0,1)$.
\end{proof}
\end{thm}
\begin{rem}\label{Remark: weak closure for bounded domains}
Let $D\subset\r^n$ be any open and bounded domain which is bi-Lipschitz equivalent to $Q_1(0)$. From Theorem \ref{weak closure for Q^n} follows that $\Omega_{p,\z}^{n-1}(D)$ (see Definition \ref{Definition: integer valued fluxes on generic domain}) is a weakly sequentially closed subspace of $\Omega_p^{n-1}(D)$. Indeed, let $\varphi:Q_1(0)\rightarrow D$ be any bi-Lipschitz homeomorphism and let $\{F_k\}_{k\in\n}\subset\Omega_{p,\z}^{n-1}(D)$ be such that $F_k\xrightharpoonup{L^p}F$ on $D$. Then, by Lemma \ref{Lemma: bi-Lipschitz maps preserve approximability properties} we have $\{\varphi^*F_k\}_{k\in\n}\subset\Omega_{p,\z}^{n-1}(Q_1(0))$ and as $\varphi$ is bi-Lipschitz we have $\varphi^*F_k\xrightharpoonup{L^p}\varphi^*F$ on $Q_1(0)$. By the Weak Closure Theorem (Theorem \ref{weak closure for Q^n}), $\varphi^*F\in \Omega_{p,\z}^{n-1}(Q_1(0))$. Thus $F\in\Omega_{p,\z}^{n-1}(D)$ (again by Lemma \ref{Lemma: bi-Lipschitz maps preserve approximability properties}). This shows that Theorem \ref{Theorem: weak closure} holds true.
\end{rem}
Observe that Theorem \ref{Theorem: weak closure} does not hold if $n=1$. In fact in this case the following holds.
\begin{lem}
\label{Lemma: weak closure for n=1}
Let $I$ be a bounded connected interval in $\mathbb{R}$. Let $p\in [1,\infty)$.\\
The weak closure of $L^p_\mathbb{Z}(I)$ in $L^p(I)$ is $L^p(I)$.
\end{lem}
\begin{proof}
Since $C^0(\overline{I})$ is dense in $L^p(I)$, it is enough to show that any function $f\in C^0(\overline{I})$ can be approximated weakly in $L^p(I)$ by functions in $L^p_\mathbb{Z}(I)$.
Without loss of generality we can assume that $I=[0,1)$.
For any $n\in \mathbb{N}$ let's define $f_n: I\to \mathbb{R}$ as follows:\\
for any $k\in \{1,..., 2^n\}$ let $I_k^n:=[\frac{k-1}{2^n}, \frac{k}{2^n})$, for any $k\in \{1,..., 2^n\}$ let $c_k:=\fint_{I_k^n}f(x)dx$ and for any $x\in I_k^n$ set
\begin{align*}
    f_n(x):=\begin{cases}
        \lceil c_k\rceil&\text{ if }x-\frac{k-1}{2^n}\leq \frac{c_k}{2^n\lceil c_k\rceil}\\
        0&\text{ otherwise}.
    \end{cases}
\end{align*}
Then $f_n\in L^p(I, \mathbb{Z})$ and $\int_{I_k^n} f_n(x)dx=\int_{I_k^n}f(x)dx$ for any $k\in \{1,..., 2^n\}$, for any $n\in \mathbb{N}$.\\
Moreover notice that since $f$ is bounded, the sequence $(f_n)_{n\in \mathbb{N}}$ is bounded in $L^p(I)$. Therefore if $p>1$ $(f_n)_{n\in \mathbb{N}}$ converges weakly in $L^p(I)$, up to a subsequence, to a function $\tilde{f}\in L^p(I)$. Testing against continuous functions on $\overline{I}$ it is easy to check that $f=\tilde{f}$.\\
If $p=1$ we have to check that, up to a subsequence,
\begin{align}
\label{Equation: weak convergence for L1 explicit}
    \lim_{n\to\infty}\int_{I} f_n g=\int_{I}fg\quad \forall\, g\in L^\infty(I). 
\end{align}
Since $L^{\infty}(I)\subset L^q(I)$ for any $q>1$, (\ref{Equation: weak convergence for L1 explicit}) follows from the case $p>1$ (with $p=q'$).
\end{proof}
\begin{rem}
Let $n\ge 2$ and let $D\subset\r^n$ be any open, bounded and Lipschitz domain in $\r^n$. It is still unknown if the space $L_{\z}^1(D)$ is weakly sequentially closed. Surely it is not weakly-$\ast$ closed, a proof this fact can be easily achieved by generalising the arguments in \cite[Section 8]{petrache-riviere-abelian}.
\end{rem}
\newpage
\appendix
\section{Minimal connections for forms with finitely many integer singularities}
Throughout Appendix A, we will denote by $M^m\subset\r^n$ some arbitrary embedded Lipschitz and connected $m$-dimensional submanifold of $\r^n$. Let $p\in [1,\infty]$. We will denote by $\Omega^1_{p,R,\infty}(M)$ the space introduced in Definition \ref{definition: Lipschitz forms}.
\begin{lem}
	\label{appendix lemma existence same boundary finitely many singularities}
	Let $F\in\Omega_{p,R,\infty}^{m-1}(M)$. Then, there exists a connection for $F$. 
	\begin{proof}
		Throughout the following proof, given any couple of points $x,y\in M$ we will denote by $(x,y)$ an arbitrarily chosen oriented Lipschitz curve with finite length joining $x$ and $y$. By assumption, it holds that 
		\begin{align*}
		*dF=\sum_{j=1}^Nd_j\delta_{x_j}, \qquad \mbox{ for some } d_1,...,d_N\in\z\smallsetminus\{0\} \mbox{ and } x_1,...,x_N\in M.
		\end{align*}
		We define 
		\begin{align*}
		\{i_1,...,i_p\}&:=\big\{j\in\{1,...,N\} \mbox{ s.t. } d_j>0\big\},\\
		\{j_1,...,j_q\}&:=\big\{j\in\{1,...,N\} \mbox{ s.t. } d_j<0\big\},\\
		d&:=\sum_{j=1}^Nd_j\in\z.
		\end{align*}
		We build a family $\mathscr{F}=\{I_{\alpha}\}_{\alpha\in A}$ of oriented Lipschitz curves in $M$ as follows. If there is no point $x_j$ such that $d_j<0$, then we set $\mathscr{F}=\emptyset$. Else, we start from $x_{i_1}$ and we add to the family $\mathscr{F}$ the curves $\big(x_{j_1},x_{i_1}\big),...,\big(x_{j_{k_1}},x_{i_1}\big)$, until we reach the condition $k_1=q$ or the condition 
		\begin{align*}
			r_1:=d_{i_1}+\sum_{l=1}^{k_1}d_{j_l}\le 0.
		\end{align*}
		If $k_1=q$, then we stop. Else, we move to the point $x_{i_2}$. 
  
        If $r_1=0$, then we add to $\mathscr{F}$ the segments $\big(x_{j_{k_1+1}},x_{i_2}\big),...,\big(x_{j_{k_2}},x_{i_2}\big)$, where $k_2\in\{1,...,q\}$ is the smallest value such that
		\begin{align*}
			r_2:=d_{i_2}+\sum_{l=k_1+1}^{k_2}d_{j_l}\le 0.
		\end{align*}
        If there is no $k\in\{1,...,q\}$ such that
        \begin{align*}
			d_{i_2}+\sum_{l=k_1+1}^{k}d_{j_l}\le 0,
		\end{align*}
        then we add to $\mathscr{F}$ the segments $\big(x_{j_{k_1+1}},x_{i_2}\big),...,\big(x_{j_{q}},x_{i_2}\big)$ and we set $k_2=q$. 
        
		If $r_1<0$, then we add to $\mathscr{F}$ the segments $\big(x_{j_{k_1}},x_{i_2}\big)$ and $\big(x_{j_{k_1+1}},x_{i_2}\big),...,\big(x_{j_{k_2}},x_{i_2}\big)$, where $k_2\in\{1,...,q\}$ is the smallest value such that
		\begin{align*}
			r_2:=d_{i_2}+r_1+\sum_{l=k_1+1}^{k_2}d_{j_l}\le 0.
		\end{align*}
        If there is no $k\in\{1,...,q\}$ such that
        \begin{align*}
			d_{i_2}+r_1+\sum_{l=k_1+1}^{k}d_{j_l}\le 0,
		\end{align*}
        then we add to $\mathscr{F}$ the segments $\big(x_{j_{k_1}},x_{i_2}\big)$ and $\big(x_{j_{k_1+1}},x_{i_2}\big),...,\big(x_{j_{q}},x_{i_2}\big)$ and we set $k_2=q$. 
        
        We proceed iteratively in this way, moving on to the subsequent points $x_{i_s}$ until $k_s=q$ or $s=p$. Then, the construction of the family $\mathscr{F}$ is complete. We let $x_{i_h}$ be the last node that is visited before the iteration stops and, for every $I_{\alpha}=(x_j,x_i)\in\mathscr{F}$, we define its multiplicity $m_{\alpha}$ as 
		\begin{align*}
			m_{\alpha}:=\begin{cases}\lvert d_j\rvert-\lvert r_l\rvert & \mbox{ if } i=i_l \mbox{ and } j=i_{k_l},\\ \min\{\lvert d_i\rvert,\lvert r_{l-1}\rvert\} & \mbox{ if } i=i_l \mbox{ and } j=i_{k_{l-1}},\\ \min\{|d_j|,|d_i|\} & \mbox{ else.}\end{cases}
		\end{align*}
		Finally, we divide three cases:
		\begin{enumerate}
			\item \textit{Case} $d=0$. Notice that this is always the case if $M$ has no boundary. We define the integer $1$-current 
			$I\in\mathcal{R}_1(M)$ given by
			\begin{align*}
				\langle I,\omega\rangle:=\sum_{\alpha\in A}m_{\alpha}\int_{I_{\alpha}}\omega, \qquad \mbox{ for every } \omega\in\mathcal{D}^1(M).
				\end{align*}
			\item \textit{Case} $d>0$. We fix a point $x_0\in\partial M$ and we let $I_s^p:=(x_0,x_{i_s})$, for every $s=h,...,p$. We define the integer $1$-current 
			$I\in\mathcal{R}_1(M)$ given by
			\begin{align*}
				\langle I,\omega\rangle:=\sum_{\alpha\in A}m_{\alpha}\int_{I_{\alpha}}\omega+r_h\int_{I_h^b}\omega+\sum_{s=h+1}^pd_{i_s}\int_{I_s^b}\omega,\qquad \mbox{ for every } \omega\in\mathcal{D}^1(M).
			\end{align*}
			\item \textit{Case} $d<0$. We fix a point $x_0\in\partial M$ and we let $I_s^n:=(x_{j_s},x_0)$, for every $s=k_h,...,q$ We define the integer $1$-current $I \in\mathcal{R}_1(M)$ given by
			\begin{align*}
				\langle I,\omega\rangle:=\sum_{\alpha\in A}m_{\alpha}\int_{I_{\alpha}}\omega+|r_h|\int_{I_h^b}\omega+\sum_{s=k_h+1}^q|d_{i_s}|\int_{I_s^b}\omega,\qquad \mbox{ for every } \omega\in\mathcal{D}^1(M).
			\end{align*}
		\end{enumerate}
		By direct computation, we verify that $I$ has the desired properties and the statement follows.
	\end{proof}
\end{lem}
\begin{lem}
	\label{appendix lemma duality}
	Let $F\in\Omega_{p,\mathbb{Z}}^{m-1}(M)$ (see Definition \ref{Definition: integer valued fluxes on generic domain}). Then, 
	\begin{align}\label{equation: estimates for connections}
	\inf_{\substack{T\in\mathcal{D}_1(M), \\ \partial T=*dF}}\mathbb{M}(T)=\inf_{\substack{T\in\mathcal{M}_1(M), \\ \partial T=*dF}}\mathbb{M}(T)=\sup_{\substack{\varphi\in W_0^{1,\infty}(M), \\ \lvert\lvert d\varphi\rvert\rvert_{L^{\infty}}\le 1}}\int_{M}F\wedge d\varphi<+\infty,
	\end{align}
	where $\mathcal{M}_1(M)$ denotes the set of all the $1$-currents with finite mass on $M$. Moreover, the infimum on the left-hand-side of the previous chain of equalities is achieved.
	\begin{proof}
		By definition, there exists an integer $1$-current $I\in\mathcal{R}_1(M)$ with finite mass such that $\partial I=*dF$. Hence
		\begin{align*}
			\inf_{\substack{T\in\mathcal{D}_1(M), \\ \partial T=*dF}}\mathbb{M}(T)\le\m(I)<+\infty.
		\end{align*}
		The first equality in \eqref{equation: estimates for connections} is clear.\\
		Notice that for every $T\in\mathcal{M}_1(M)$ such that $\partial T=*dF$ it holds that
		\begin{align*}
			\int_{M}F\wedge d\varphi=\langle*dF,\varphi\rangle=\langle\partial T,\varphi\rangle=\langle T,d\varphi\rangle\le\mathbb{M}(T)\lvert\lvert d\varphi\rvert\rvert_{L^{\infty}(M)}, \qquad\forall\,\varphi\in W_0^{1,\infty}(M).
		\end{align*}
		Hence, 
		\begin{align}
		\label{appendix equation lemma duality}
			\inf_{\substack{T\in\mathcal{M}_1(M), \\ \partial T=*dF}}\mathbb{M}(T)\ge\sup_{\substack{\varphi\in W_0^{1,\infty}(M), \\ \lvert\lvert d\varphi\rvert\rvert_{L^{\infty}}\le 1}}\int_{M}F\wedge d\varphi.
		\end{align}
		To prove that the former inequality is actually an equality, it suffices to show that the supremum on its right-hand-side is greater than the mass of some $1$-current with finite mass $T$ on $M$ such that $\partial T=\ast dF$. Define the vector subspace $X\subset\Omega_{\infty}^1(M)$ given by
		\begin{align*}
			X:=\{\omega\in\Omega_{\infty}^1(M) \mbox{ s.t. } \omega=d\varphi, \mbox{ for some } \varphi\in W_0^{1,\infty}(M)\}.
		\end{align*}
		Consider the linear functional $\phi:X\subset(\Omega^1(M),\lVert\cdot\rVert_{L^\infty})\rightarrow\r$ given by
		\begin{align*}
			\langle\phi,\omega\rangle=\int_{M}F\wedge\omega, \qquad\forall\,\omega\in X.
		\end{align*}
		By \eqref{appendix equation lemma duality} we get that $\phi$ is continuous on $X$, i.e. 
		\begin{align*}
			\lvert\lvert\phi\rvert\rvert_{\L(X)}=\sup_{\substack{\omega\in X, \\ \lvert\lvert \omega\rvert\rvert_{L^{\infty}}\le 1}}\int_{M}F\wedge\omega=\sup_{\substack{\varphi\in W_0^{1,\infty}(M), \\ \lvert\lvert d\varphi\rvert\rvert_{L^{\infty}}\le 1}}\int_{M}F\wedge d\varphi\le \inf_{\substack{T\in\mathcal{M}_1(M), \\ \partial T=*dF}}\mathbb{M}(T)<+\infty.
		\end{align*}
		By Hahn-Banach theorem,
		we can extend $\phi$ to a linear functional $T:\Omega^1(M)\rightarrow\r$ such that 
        \begin{align*}    \lvert\lvert T \rvert\rvert_{\L(\Omega_{\infty}^1(M))}=\lvert\lvert\phi\rvert\rvert_{\L(X)}=\sup_{\substack{\varphi\in W_0^{1,\infty}(M), \\ \lvert\lvert d\varphi\rvert\rvert_{L^{\infty}}\le 1}}\int_{M}F\wedge d\varphi
		\end{align*}
		But then, $T$ is a $1$-current on $M$ having finite mass and such that 
		\begin{align*}
			\m(T)\le\lvert\lvert T\rvert\rvert_{\L(\Omega_{\infty}^1(M))}=\sup_{\substack{\varphi\in W_0^{1,\infty}(M), \\ \lvert\lvert d\varphi\rvert\rvert_{L^{\infty}}\le 1}}\int_{M}F\wedge d\varphi.
		\end{align*}
        Moreover, 
        \begin{align*}
            \langle\partial T,\varphi\rangle=\langle T,d\varphi\rangle&=\langle\phi,d\varphi\rangle=\int_{M}F\wedge d\varphi=\langle\ast dF,\varphi\rangle, \qquad\forall\,\varphi\in W_0^{1,\infty}(M).
        \end{align*}
        Hence, 
        \begin{align*}
	       \m(T)=\inf_{\substack{\tilde T\in\mathcal{D}_1(M), \\ \partial \tilde T=*dF}}\mathbb{M}(\tilde T)=\inf_{\substack{\tilde T\in\mathcal{M}_1(M), \\ \partial \tilde T=*dF}}\mathbb{M}(\tilde T)=\sup_{\substack{\varphi\in W_0^{1,\infty}(M), \\ \lvert\lvert d\varphi\rvert\rvert_{L^{\infty}}\le 1}}\int_{M}F\wedge d\varphi<+\infty
	   \end{align*}
     and the statement follows.
	\end{proof}
\end{lem}
\begin{prop}
	\label{appendix proposition existence minimal connection finitely many singularities}
	Let $F\in\Omega_{p,R}^{m-1}(M)$. Then, there exists an integer $1$-current $L\in\mathcal{R}_1(M)$ such that $\partial L=*dF$ and
	\begin{align*}
		\m(L)=\inf_{\substack{T\in\mathcal{D}_1(M), \\ \partial T=*dF}}\mathbb{M}(T)=\sup_{\substack{\varphi\in W_0^{1,\infty}(M), \\ \lvert\lvert d\varphi\rvert\rvert_{L^{\infty}}\le 1}}\int_{M}F\wedge d\varphi.
	\end{align*}
	In particular, 
	\begin{align*}
	\m(L)\le C\lvert\lvert F\rvert\rvert_{L^p}.
	\end{align*}
	\begin{proof}
		Notice that by \cite[Chapter 1, Section 3.4, Theorem 8]{giaquinta-modice-soucek_cartesan_currents_2}, we have 
		\begin{align*}
		\inf_{\substack{T\in\mathcal{R}_1(M), \\ \partial T=*dF}}\mathbb{M}(T)=\inf_{\substack{T\in\mathcal{D}_1(M), \\ \partial T=*dF}}\mathbb{M}(T).
		\end{align*}
		Since the mass $\m(\,\cdot\,)$ is lower semicontinuous with respect to the weak convergence in $\D_1(M)$ and since $\m$-bounded subsets of the competition class $\mathcal{R}_1(M)\cap\{T\in\mathcal{D}_1(M) \mbox{ s.t. } \partial T=*dF\}$ are weakly sequentially compact (for a reference, see e.g. \cite[Equation (7.5), Theorem 7.5.2]{krantz_parks-geometric_integration_theory}), by the direct method of calculus of variations we conclude that there exists an integer $1$-current $L\in\mathcal{R}_1(M)$ such that $\partial L=*dF$ and 
		\begin{align*}
		\mathbb{M}(L)=\inf_{\substack{T\in\mathcal{R}_1(M), \\ \partial T=*dF}}\mathbb{M}(T)=\inf_{\substack{T\in\mathcal{D}_1(M), \\ \partial T=*dF}}\mathbb{M}(T)=\sup_{\substack{\varphi\in W_0^{1,\infty}(M), \\ ||d\varphi||_{L^{\infty}}\le 1}}\int_{M}F\wedge d\varphi,
		\end{align*}
		where the last equality follows from Lemma \ref{appendix lemma duality}. The statement follows.
	\end{proof}
\end{prop}
\section{Laplace equation on spheres}
Let $n\in\n$ be such that $n\ge 2$ and fix any $p\in(1,+\infty)$. We let 
\begin{align*}
    \dot W^{1,p}(\s^{n-1}):=\bigg\{u\in W^{1,p}(\s^{n-1}) \mbox{ s.t. } \bar u:=\int_{\s^{n-1}}u\, d\vol_{\s^{n-1}}=0\bigg\}.
\end{align*}
We can endow the space $\dot W^{1,p}(\s^{n-1})$ with the usual $W^{1,p}$-norm induced by $W^{1,p}(\s^{n-1})$, given by 
\begin{align*}
    \lvert\lvert u\rvert\rvert_{W^{1,p}}:=\lvert\lvert u\rvert\rvert_{L^p}+\lvert\lvert du\rvert\rvert_{L^p}, \qquad\forall\, u\in \dot W^{1,p}(\s^{n-1}).
\end{align*}
\begin{lem}[Poincar\'e inequality on $\dot W^{1,p}$]
\label{appendix poincare inequality}
There exists a constant $C>0$ such that 
\begin{align*}
    \int_{\s^{n-1}}|u|^p\,d\vol_{\s^{n-1}}\le C\int_{\s^{n-1}}|du|^p\,d\vol_{\s^{n-1}}, \qquad\forall\, u\in \dot W^{1,p}(\s^{n-1}).
\end{align*}
\begin{proof}
By contradiction, assume that for every $k>0$ there exists $u_k\in \dot W^{1,p}(\s^{n-1})$ such that $\lvert\lvert u_k\rvert\rvert_{L^p}=1$ and 
\begin{align*}
    1>k\int_{\s^{n-1}}|du_k|^p\,d\vol_{\s^{n-1}}. 
\end{align*}
This implies immediately that $\lvert\lvert du_k\rvert\rvert_{L^p}\rightarrow 0$ as $k\rightarrow\infty$. In particular, the sequence $\{u_k\}_{k\in\n}$ is bounded with respect to the $W^{1,p}$-norm. Hence, by weak compactness of $W^{1,p}(\s^{n-1})$, there exists a subsequence $\{u_{k_j}\}_{j\in\n}$ of $\{u_k\}_{k\in\n}$ such that $u_{k_j}\rightharpoonup u$ in $W^{1,p}(\s^{n-1})$. Moreover, by Rellich-Kondrakov theorem, we have $u_{k_j}\rightarrow u$ strongly in $L^p(\s^{n-1})$. Since $du_{k_j}\rightarrow 0$ strongly in $L^p$ we get $du=0$. Then, $u$ is constant on $\s^{n-1}$. Since $u_{k_j}\rightharpoonup u$ in $L^p(\s^{n-1})$, it follows that
\begin{align*}
    0=\lim_{j\rightarrow\infty}\int_{\s^{n-1}}u_{k_j}\, d\vol_{\s^{n-1}}=\int_{\s^{n-1}}u\, d\vol_{\s^{n-1}}
\end{align*}
and this leads to $u=0$. But this is absurd, since by strong $L^p$-convergence of $\{u_{k_j}\}_{j\in\n}$ to $u$ we obtain $\lvert\lvert u\rvert\rvert_{L^p}=1$. 
\end{proof}
\end{lem}
\begin{rem}
By Lemma \ref{appendix poincare inequality}, we conclude that we can endow $\dot W^{1,p}$ with the following much more convenient norm:
\begin{align*}
    \lvert\lvert u\rvert\rvert_{\dot W^{1,p}}:=||du||_{L^p}, \qquad\forall\, u\in \dot W^{1,p}(\s^{n-1}).
\end{align*}
Moreover, such a norm is equivalent to $W^{1,p}$-norm.
\end{rem}
\begin{rem}
\label{appendix remark continuity}
Notice that a linear functional on $W^{1,p}(\s^{n-1})$ restricts to an element of $(\dot W^{1,p'}(\s^{n-1}))^*$ if and only if it is $W^{1,p}$-continuous and $\langle F,1\rangle=0$. 
\end{rem}
\begin{lem}
\label{appendix lemma d* e d}
Let $F\in(W^{1,p'}(\s^{n-1}))^*$ be such that $\langle F,1\rangle=0$. Then, the following facts hold.
\begin{enumerate}
    \item If $n\ge 3$ the linear differential system
        \begin{align*}
            \begin{cases}
                d^*\omega=F\\
                d\omega=0
            \end{cases}
        \end{align*}
        has a unique weak solution $\alpha\in\Omega_p^1(\s^{n-1})$.
    \item If $n=2$ the linear differential system
        \begin{align*}
            \begin{cases}
                d^*\omega=F\\
                d\omega=0\\
                \displaystyle{\int_{\s^1}}\omega=0
            \end{cases}
        \end{align*}
        has a unique weak solution $\alpha\in\Omega_p^1(\s^1)$.
\end{enumerate}
In both cases, $\alpha$ satisfies the following estimate:
\begin{align*}
    ||\alpha||_{L^p}\le C\lvert\lvert F\rvert\rvert_{\mathcal{L}(\dot W^{1,p'}(\s^{n-1}))},
\end{align*}
for some constant $C>0$ depending only $n$.
\begin{proof}
Observe that, by Remark \ref{appendix remark continuity}, $F$ restricts to an element of $(\dot W^{1,p'}(\s^{n-1}))^*$. Consider the linear functional $\phi:\Omega^1(\s^{n-1})\rightarrow\r$ given by 
\begin{align*}
    \langle\phi,\omega\rangle=\langle F,u\rangle, \qquad\forall\,\omega=du+d^*\beta+\eta\in\Omega^1(\s^{n-1}), 
\end{align*}
where $\eta$ is a harmonic $1$-form on $\s^{n-1}$. 
Since $\phi$ is continuous and linear on $\Omega^1(\s^{n-1})$ w.r.t. the $L^{p'}$-norm, by Hahn-Banach theorem there exists a unique (recall that $L^p$-spaces are strictly convex) extension $\Phi\in(\Omega_{p'}^1(\s^{n-1}))^*$ of $\phi$ such that 
\begin{align*}
    ||\Phi||_{(\Omega_{p'}^1(\s^{n-1}))^*}=||\phi||_{(\Omega^1(\s^{n-1}))^*}\le C\lvert\lvert F\rvert\rvert_{\mathcal{L}(\dot W^{1,p'}(\s^{n-1}))}.
\end{align*}
By Riesz representation theorem, there exists a unique $\alpha\in\Omega_{p}^1(\s^{n-1})$ such that 
\begin{align}
\label{appendix riesz representation}
    \langle\alpha,\omega\rangle_{L^p-L^{p'}}:=\int_{\s^{n-1}}\alpha\wedge*\omega=\langle\Phi,\omega\rangle, \qquad\forall\,\omega\in\Omega_{p'}^1(\s^{n-1})
\end{align}
and 
\begin{align*}
     ||\alpha||_{L^p}=||\Phi||_{(\Omega_{p'}^1(\s^{n-1}))^*}\le C\lvert\lvert F\rvert\rvert_{(\dot W^{1,p'}(\s^{n-1}))^*}. 
\end{align*}
Finally, by applying equation \eqref{appendix riesz representation} with get 
\begin{align*}
    \langle\alpha,du\rangle_{L^p-L^{p'}}&=\langle\Phi,du\rangle=\langle\phi,du\rangle=\langle F,u\rangle, \qquad\forall\,u\in C^{\infty}(\s^{n-1}),
\end{align*}
and
\begin{align*}
    \langle\alpha,d^*\beta\rangle_{L^p-L^{p'}}&=\langle\Phi,d^*\beta\rangle=\langle\phi,d^*\beta\rangle=0, \qquad\forall\,\beta\in\Omega^2(\s^{n-1}).
\end{align*}
The two previous equations are exactly the weak forms of the equations $d^*\alpha=F$ and $d\alpha=0$ respectively. Moreover, in case $n=2$, we have 
\begin{align*}
    \int_{\s^1}\alpha=\int_{\s^1}\alpha\wedge\ast 1=\langle\alpha,\ast 1\rangle_{L^p-L^{p'}}=\langle\Phi,\ast 1\rangle=\langle F,1\rangle =0.
\end{align*}
This concludes about the existence of a solution to the differential systems given in points 1 and 2. For what concerns uniqueness, assume that $\alpha$ and $\alpha'$ are two solutions of the differential system given in point 1 (resp. 2) and define $\beta=\alpha-\alpha'$. Then, we distinguish the two cases:

\textit{Case $n\ge 3$}. In this case, $\beta$ satisfies
\begin{align*}
    \begin{cases}
        d^*\beta=0\\
        d\beta=0.
    \end{cases}
\end{align*}
Hence, $\beta$ is a harmonic $1$-forms on $\s^{n-1}$ for $n\ge 1$, which implies $\beta=0$.

\textit{Case $n=2$}. In this case, $\beta$ satisfies
\begin{align*}
    \begin{cases}
        d^*\beta=0\\
        d\beta=0\\
        \displaystyle{\int_{\s^1}}\beta=0.
    \end{cases}
\end{align*}
Hence, $\beta$ is a harmonic $1$-forms on $\s^1$, which implies $\beta=c\vol_{\s^1}$ for some $c\in\r$. But since $\beta$ has vanishing integral on $\s^1$, we get $\beta=0$. 
\end{proof}
\end{lem}
\begin{dfn}[Sobolev spaces of differential forms]
\label{appendix definition sobolev forms}
Fix any $k\in\n\smallsetminus\{0\}$. We define the \textit{Sobolev space of $W^{1,p}$-regular differential $k$-forms on $\s^{n-1}$} by
\begin{align*}
    \Omega_{W^{1,p}}^k(\s^{n-1}):=\big\{\omega\in\Omega_{p}^k(\s^{n-1}) \mbox{ s.t. } d\omega,d^*\omega\in L^p\big\}.
\end{align*}
We endow such space with the norm 
\begin{align*}
    \lvert\lvert\omega\rvert\rvert_{W^{1,p}}:=\lvert\lvert\omega\rvert\rvert_{L^p}+\lvert\lvert d\omega\rvert\rvert_{L^p}+\lvert\lvert d^*\omega\rvert\rvert_{L^p}, \qquad\forall\,\omega\in\Omega_{W^{1,p}}^k(\s^{n-1}).
\end{align*}
\end{dfn}
\begin{rem}
\label{appendix remark sobolev forms}
It can be shown (see \cite[\S 3 and \S 4]{scott-L^p_theory_of_differential_forms_on_manifolds}) that such Sobolev spaces are completely equivalent to the usual ones, namely the space of $k$-forms having local coefficients in $W^{1,p}$. Moreover, in case $n\ge 3$ there exists $C>0$ such that 
\begin{align}\label{poincare on forms}
    \lvert\lvert\omega\rvert\rvert_{W^{1,p}}\le C\big(\lvert\lvert d\omega\rvert\rvert_{L^p}+\lvert\lvert d^*\omega\rvert\rvert_{L^p}\big), \qquad\forall\,\omega\in\Omega_{W^{1,p}}^1(\s^{n-1}).
\end{align}
Indeed, let 
\begin{align*}
    X&:=\{d\alpha \mbox{ s.t. } \alpha\in\Omega_{W^{1,p}}^1(\s^{n-1})\},\\
    Y&:=\{d^*\beta \mbox{ s.t. } \beta\in\Omega_{W^{1,p}}^1(\s^{n-1})\}.
\end{align*}
By \cite[Proposition 7.1]{scott-L^p_theory_of_differential_forms_on_manifolds}, both $X$ and $Y$ are closed linear subspaces  respectively of $\Omega_p^2(\s^{n-1})$ and $\dot W^{1,p}(\s^{n-1})$. Then, $X\oplus Y$ is a Banach space with respect to the standard norm on the direct sum of two Banach spaces. We claim that the liner operator $T:\Omega_{W^{1,p}}^1(\s^{n-1})\rightarrow X\oplus Y$ given by $T\omega=(d\omega,d^*\omega)$, for every $\omega\in\Omega_{W^{1,p}}^1(\s^{n-1})$ is a continuous linear bijection between Banach spaces. Indeed, the fact that $T$ is injective follows form the fact that there no non-zero harmonic forms on $\s^{n-1}$ for $n\ge 3$. Hence, we just need to show that $T$ is surjective. Pick any $(d\alpha,d^*\beta)\in X\oplus Y$. By Lemma \ref{appendix lemma d* e d}, the linear differential system 
\begin{align*}
    \begin{cases}
        d^*\omega=d^*(\beta-\alpha)\\
        d\omega=0
    \end{cases}
\end{align*}
has a unique weak solution $\tilde\omega\in\Omega_p^1(\s^{n-1})$. Since by construction we have $d\tilde\omega, d^*\tilde\omega\in L^p$, we conclude that $\tilde\omega\in\Omega_{W^{1,p}}^1(\s^{n-1})$. Then, by letting $\omega:=\tilde\omega+\alpha\in\Omega_{W^{1,p}}^1(\s^{n-1})$ we see that 
\begin{align*}
    T\omega=(d\omega,d^*\omega)=(d\tilde\omega+d\alpha,d^*\tilde\omega+d^*\alpha)=(d\alpha,d^*\beta)
\end{align*}
and we have proved our claim. This proves that $T$ has a continuous inverse and the statement follows with $C=\lvert\lvert T^{-1}\rvert\rvert_{\mathcal{L}(X\oplus Y,\Omega_{W^{1,p}}^1(\s^{n-1}))}$. 

In case $n=2$, the estimate \eqref{poincare on forms} still holds for every $\omega\in\Omega_{W^{1,p}}^1(\s^1)$ such that
\begin{align*}
    \int_{\s^1}\omega=0.
\end{align*}
The proof is completely analogous.
\end{rem}
\begin{rem}[$L^p$-Hodge decomposition]
\label{appendix remark hodge}
Let 
\begin{align*}
    X&:=\{d\varphi \mbox{ s.t. } \varphi\in\dot W^{1,p}(\s^{n-1})\},\\
    Y&:=\{d^*\beta \mbox{ s.t. } \beta\in\Omega_{W^{1,p}}^2(\s^{n-1})\},\\
    Z&:=\{\eta\in\Omega^1(\s^{n-1}) \mbox{ s.t. } \Delta\eta=(dd^*+d^*d)\eta=0\}.
\end{align*}
Then, as a particular consequence of the $L^p$-Hodge decomposition theorem (see e.g. \cite[Proposition 6.5]{scott-L^p_theory_of_differential_forms_on_manifolds}), the operator $T:X\oplus Y\oplus Z\rightarrow\Omega_p^1(\s^{n-1})$ given by \begin{align*}
    T(d\varphi,d^*\beta,\eta):=d\varphi+d^*\beta+\eta
\end{align*} 
is a continuous and linear isomorphism between Banach spaces. Hence, $T$ has a continuous inverse. We let
\begin{align*}
    C_H:=\lvert\lvert T^{-1}\rvert\rvert_{\mathcal{L}(\Omega_p^1(\s^{n-1}),X\oplus Y)}.
\end{align*}
We conclude that for every $\omega\in\Omega_p^1(\s^{n-1})$ there exist $\varphi\in\dot W^{1,p}(\s^{n-1})$, $\beta\in\Omega_{W^{1,p}}^2(\s^{n-1})$ and $\eta\in Z$ such that $\omega=d\varphi+d^*\beta+\eta$ and 
\begin{align}
\label{appendix equation useful estimate}
\lvert\lvert d\varphi\rvert\rvert_{L^p}+\lvert\lvert d^*\beta\rvert\rvert_{L^p}+\lvert\lvert\eta\rvert\rvert_{L^p}\le C_H\lvert\lvert\omega\rvert\rvert_{L^p}.
\end{align}
\end{rem}
\begin{lem}[A weak version of Poincarè lemma]
\label{appendix lemma sobolev poincare}
Let $n\ge 3$ and let $\alpha\in\Omega_{p}^1(\s^{n-1})$ be such that $d\alpha=0$ weakly on $\s^{n-1}$. Then, there exists a Sobolev function $\varphi\in \dot W^{1,p}(\s^{n-1})$ such that $d\varphi=\alpha$ weakly on $\s^{n-1}$.
\begin{proof}
We follow the notation of Remark \ref{appendix remark hodge} and we notice that, since $n\ge 3$ we have $Z=\{0\}$. Hence, we write $\alpha=d\varphi+d^*\beta$, for $\varphi\in\dot W^{1,p}(\s^{n-1})$ and $\beta\in\Omega_{W^{1,p}}^2(\s^{n-1})$. We observe that
\begin{align*}
    \langle d^*\beta,\omega\rangle_{L^p-L^{p'}}&=\langle d^*\beta,d\psi+d^*\gamma\rangle_{L^p-L^{p'}}\\
    &=\langle d^*\beta,d^*\gamma\rangle_{L^p-L^{p'}}\\
    &=\langle\alpha-d\varphi,d^*\gamma\rangle_{L^p-L^{p'}}\\
    &=\langle \alpha,d^*\gamma\rangle_{L^p-L^{p'}}=0, \qquad\forall\,\omega=d\psi+d^*\gamma\in\Omega^1(\s^{n-1}).
\end{align*}
This implies $d^*\beta=0$ and the statement follows. 
\end{proof}
\end{lem}
\begin{cor}[Laplace equation on spheres]
\label{appendix corollary laplace}
Let $F\in\mathcal{L}(W^{1,p'}(\s^{n-1}))$ such that $\langle F,1\rangle=0$. Then, the linear differential equation
\begin{align*}
   \Delta u = F 
\end{align*}
has a unique weak solution $\varphi\in\dot W^{1,p}(\s^{n-1})$ satisfying 
\begin{align*}
    ||\varphi||_{W^{1,p}}\le C\lvert\lvert F\rvert\rvert_{(\dot W^{1,p'}(\s^{n-1}))^*}.
\end{align*}
\begin{proof}
First, we face the case $n\ge 3$. By Lemma \ref{appendix lemma d* e d} we can find $\alpha\in\Omega_p^1(\s^{n-1})$ satisfying 
\begin{align*}
    \begin{cases}
        d^*\alpha=F\\
        d\alpha=0
    \end{cases}
\end{align*}
and 
\begin{align*}
    ||\alpha||_{L^p}\le C\lvert\lvert F\rvert\rvert_{(\dot W^{1,p'}(\s^{n-1}))^*}.
\end{align*}
Since $d\alpha=0$, by Lemma \ref{appendix lemma sobolev poincare} there exists $\varphi\in\dot W^{1,p}(\s^{n-1})$ such that $\alpha=d\varphi$. Hence, we get
\begin{align*}
    \Delta\varphi=d^*d\varphi=d^*\alpha=F. 
\end{align*}
Moreover, by Lemma \ref{appendix poincare inequality}, we have 
\begin{align*}
    ||\varphi||_{W^{1,p}}\le C\lvert\lvert d\varphi\rvert\rvert_{L^p}=C||\alpha||_{L^p}\le C\lvert\lvert F\rvert\rvert_{(\dot W^{1,p'}(\s^{n-1}))^*}.
\end{align*}
This concludes the proof in case $n\ge 3$.

If $n=2$, then by Lemma \ref{appendix lemma d* e d} we can find $\alpha\in\Omega_p^1(\s^{n-1})$ satisfying 
\begin{align*}
    \begin{cases}
        d^*\alpha=F\\
        d\alpha=0\\
        \displaystyle{\int_{\s^1}\alpha=0}
    \end{cases}
\end{align*}
and 
\begin{align*}
    ||\alpha||_{L^p}\le C\lvert\lvert F\rvert\rvert_{(\dot W^{1,p'}(\s^{n-1}))^*}.
\end{align*}
By setting $\varphi:=\ast\alpha$, the statement follows. 
\end{proof}
\end{cor}
\section{Some technical lemmata}
In this section we will make use of the following notation:
let $T$ be an $m$-rectifiable current in $\mathbb{R}^n$, then $T$ can be represented as follows:
\begin{align*}
    \langle T,\omega\rangle=\int_{\mathbb{R}^n}\theta\langle \omega, \xi\rangle d\mathscr{H}^m\big|_\Sigma\quad\forall \omega\in \Omega^m(\mathbb{R}^n),
\end{align*}
where $\Sigma$ is a locally $m$-rectifiable set,
\begin{align*}
    \theta: \Sigma\to \mathbb{Z}
\end{align*}
is a locally $\mathscr{H}^m$-integrable, non-negative function
and
\begin{align*}
    \xi:\Sigma\to\Lambda_m\mathbb{R}^n
\end{align*}
is an $\mathscr{H}^m$-measurable function such that for $\mathscr{H}^m$-almost every point $x\in \Sigma$, $\xi(x)$ is a simple unit $m$-vector in $T_x\Sigma$.\\
In this case we write
\begin{align*}
    T=\tau(\Sigma, \theta, \xi).
\end{align*}
\begin{lem}
\label{appendix Lemma closure integer currents}
For any $k\in \mathbb{N}$ let
\begin{align*}
    T_k=\tau( \Sigma_k, \theta_k, \xi_k)
\end{align*}
be an $m$-rectifiable current on $\mathbb{R}^n$ of finite mass.
Assume that $(T_k)_{k\in \mathbb{N}}$ is a Cauchy sequence with respect to the convergence in mass.\\
Then there exists an $m$-rectifiable current
\begin{align*}
    T=\tau( \Sigma, \theta, \xi)
\end{align*}
such that
\begin{align*}
    T_k\to T\quad (k\to\infty)\quad\text{in mass.}
\end{align*}
\end{lem}
\begin{proof}
    Replacing the original sequence by a subsequence if necessary, we may assume that for any $k\in \mathbb{N}$
    \begin{align*}
        \mathbb{M}(T_k- T_{k+1})\leq 2^{-k}.
    \end{align*}
    Now for any $k\in \mathbb{N}$ let
    \begin{align*}
        \tilde{T}_k:=\begin{cases}
            T_1&\text{ if }k=1\\
            T_k-T_{k-1}&\text{ if }k>1.
        \end{cases}
    \end{align*}
    Then for any $k\in \mathbb{N}$ we have
    \begin{align*}
        T_k=\sum_{i=1}^k\tilde{T}_i.
    \end{align*}
    For any $k\in \mathbb{N}$ write
    \begin{align*}
        \tilde{T}_k=\langle \tilde{\Sigma}_k, \tilde{\theta}_k, \tilde{\xi}_k\rangle.
    \end{align*}
    Notice that
    \begin{align*}
        \sum_{i=1}^k\tilde\theta_i\tilde\xi_i=\theta_k\xi_k \qquad\H^m\mbox{- a.e. on } \Sigma,
    \end{align*}
    for every $k\in\n$.
    Set
    \begin{align*}
        \Sigma:=\bigcup_{k\in \mathbb{N}}\big(\tilde{\Sigma}_k\smallsetminus\tilde\theta_k^{-1}(0)\big).
    \end{align*}
    Then $\Sigma$ is $m$-rectifiable as countable union of $m$-rectifiable sets.
    Moreover $\mathscr{H}^m(\Sigma)<\infty$.
    In fact
    \begin{align*}
        \mathscr{H}^m(\Sigma)\leq\sum_{k\in\mathbb{N}}\mathscr{H}^m\big(\tilde{\Sigma}_k\smallsetminus\tilde\theta_k^{-1}(0)\big)\leq\sum_{k\in \mathbb{N}}\int_{\mathbb{R}^n}\lvert\tilde{\theta}_k\rvert\, d\mathscr{H}^m\res{\tilde{\Sigma}_k}=\sum_{k\in \mathbb{N}}\mathbb{M}(\tilde{T}_k)<\infty.
    \end{align*}
    Next let
    \begin{align*}
        \theta=\sum_{k\in \mathbb{N}}\tilde{\theta}_k,
    \end{align*}
where $\tilde\theta_k$ is extended by zero on $\Sigma\smallsetminus\tilde\Sigma_k$ for any $k\in \mathbb{N}$ .

By Beppo-Levi Theorem
\begin{align}\label{Equation: theta is integrable}
\nonumber
    \int_{\r^n}\lvert\theta\rvert\, d\H^m\res\Sigma&=\sum_{k\in \mathbb{N}}\int_{\mathbb{R}^n}\lvert\tilde{\theta}_k\rvert\, d\mathscr{H}^m\res\Sigma\\
    &=\sum_{k\in \mathbb{N}}\int_{\mathbb{R}^n}\lvert\tilde\theta_k\rvert\, d\mathscr{H}^m\res \tilde{\Sigma}_k=\sum_{k\in \mathbb{N}}\mathbb{M}(\tilde{T}_k)<\infty.
\end{align}
    Therefore $\theta$ is finite $\mathscr{H}^m$-a.e. in $\Sigma$, i.e. for $\mathscr{H}^m$-a.e. $x\in \Sigma$ there are only finitely many $k\in \mathbb{N}$ so that $\theta_k(x)\neq 0$.
    In particular the sum
    \begin{align*}
        \sum_{k\in \mathbb{N}}\tilde\theta_k(x)\xi_k(x)
    \end{align*}
    is well defined and finite for $\mathscr{H}^m$-a.e. $x\in \Sigma$ (again $\xi_k$ is extended by zero on $\Sigma\smallsetminus\tilde\Sigma_k$ for any $k\in \mathbb{N}$) and we can write
    \begin{align*}
        \theta(x)\xi(x)=\sum_{k\in \mathbb{N}}\tilde\theta_k(x)\xi_k(x)
    \end{align*}
    for some $\theta(x)\in\mathbb{Z}_{\geq 0}$ and for some simple unit $m$-vector $\xi(x)$ in $T_x\Sigma$, for $\H^m$-a.e. $x\in\Sigma$.
    
    Observe that $\theta$ is an $\mathscr{H}^m$-measurable function on $\Sigma$ as the absolute value of the a.e.-limit of $\mathscr{H}^m$-measurable functions. Analogously, $\xi$ is an $\mathscr{H}^m$-measurable map on $\Sigma$ as the a.e.-limit of $\mathscr{H}^m$-measurable maps.
    We set
    \begin{align*}
        T:=\tau(\Sigma, \theta, \xi)
    \end{align*}
    and we claim that
    \begin{align*}
        T_k\to T\quad(k\to\infty)\quad\text{in mass.}
    \end{align*}
    In fact we know that since the space of $m$-currents is complete under the convergence in mass (as a dual space), there exists an $m$-current $T'$ such that
    \begin{align*}
        T_k\to T'\quad(k\to\infty)\quad\text{in mass.}
    \end{align*}
    To see that $T=T'$ observe that for any $\omega\in\D^m(\mathbb{R}^n)$
    \begin{align*}
        \theta_k\langle \omega, \xi_k\rangle=\sum_{i=1}^k\tilde\theta_i\langle\omega,\tilde\xi_i\rangle\to\theta\langle \omega, \xi\rangle\quad \mathscr{H}^m\text{-a.e. in }\Sigma,
    \end{align*}
    thus by (\ref{Equation: theta is integrable}) and Dominated Convergence Theorem we conclude that
    \begin{align*}
        \langle T_k,\omega\rangle\to \langle T,\omega\rangle\quad(k\to\infty).
    \end{align*}
    In particular $T=T'$.
\end{proof}
\begin{lem}\label{lemma: continuity of the dilation operator}
Let $\alpha\in(1,+\infty)$, $q\in(-\infty,1]$, $\eps\in(0,1)$ and let $\Omega\subset Q_{1-\eps}^n(0)$ be open, Lipschitz and bounded. For any $p\in[1,+\infty)$ and $\mu:=f\L^n$ with $f=(\frac{1}{2}-\lVert\,\cdot\,\rVert_{\infty})^q$, consider the continuous linear operator $P_{\alpha}:L^p(\Omega,\mu;\r^n)\rightarrow L^p(\Omega,\mu;\r^n)$ given by
\begin{align*}
    (P_{\alpha}V)(x):=\begin{cases}\alpha^{n-1}V(\alpha x)& \mbox{ if }\,x\in\alpha^{-1}\Omega,\\
    0 &\mbox{ on } \Omega\smallsetminus\alpha^{-1}\Omega.\end{cases}
\end{align*}
Then:
\begin{enumerate}
    \item For every $\alpha\in (1,+\infty)$ such that $\lvert 1-\alpha^{-1}\rvert\le\eps$ holds that
    \begin{align*}
        \lVert P_{\alpha}V\rVert_{L^p(\mu)}\le C\alpha^{n-1-\frac{n}{p}}\lVert V\rVert_{L^p(\mu)},
    \end{align*}
    for some constant $C>0$ depending only on $q$ and $p$. 
    \item For every $V\in L^p(\Omega,\mu;\r^n)$ we have that $P_{\alpha}V\rightarrow V$ in $L^p(\Omega,\mu;\r^n)$ as $\alpha\rightarrow 1^+$.
\end{enumerate}
\begin{proof}
Fist we prove 1. Fix any $V\in L^p(\Omega,\mu;\r^n)$ and compute
\begin{align*}
    \int_{\Omega}\lvert P_{\alpha}V\rvert^p\, d\mu&=\alpha^{p(n-1)}\int_{\alpha^{-1}\Omega}\lvert V(\alpha x)\rvert^p\, d\mu(x)\\
    &\le\alpha^{p(n-1)-n}\bigg(\int_{\Omega}\lvert V(y)\rvert^p\, d\mu(y)+\int_{\Omega}\lvert V(y)\rvert^p\frac{f(\alpha^{-1}y)-f(y)}{f(y)}\, d\mu(y)\bigg).
\end{align*}
As in \eqref{estimate on the difference measure} we can estimate
\begin{align*}
    \left\lvert\frac{f(\alpha^{-1}y)-f(y)}{f(y)}\right\rvert\leq q\left(\frac{1}{2}-\lVert y\rVert_\infty\right)^{-1}\lVert y\rVert_\infty(1-\alpha^{-1})\leq C
\end{align*}
for any $y\in Q_{1-\varepsilon}(0)$ and any $\alpha\geq 1$ such that $\lvert 1-\alpha^{-1}\rvert\leq \varepsilon$, for some constant $C$ depending only on $q$. Therefore
\begin{align*}
    \int_{\Omega}\lvert P_{\alpha}V\rvert^p\, d\mu\le (C+1)\alpha^{p(n-1)-n}\int_{\Omega}\lvert V(y)\rvert^p\, d\mu(y).
\end{align*}
Hence $1.$ follows.
We are left to prove $2.$. Fix any $\delta>0$ and let $V_{\delta}\in C_c^{0}(\Omega;\r^n)$ be such that 
\begin{align*}
    \lVert V_{\delta}-V\rVert_{L^p(\mu)}\le\delta.
\end{align*}
By 1, we have 
\begin{align*}
    \lVert P_{\alpha}V-V\rVert_{L^p(\mu)}&\le\lVert P_{\alpha}(V-V_{\delta})\rVert_{L^p(\mu)}+\lVert P_{\alpha}V_{\delta}-V_{\delta}\rVert_{L^p(\mu)}+\lVert V_{\delta}-V\rVert_{L^p(\mu)}\\
    &\le (C\alpha^{n-1-\frac{n}{p}}+1)\delta+\lVert P_{\alpha}V_{\delta}-V_{\delta}\rVert_{L^p(\mu)},
\end{align*}
for every $\alpha\in(1,+\infty)$ such that $\lvert 1-\alpha^{-1}\rvert\le\varepsilon$. Since $V_{\delta}$ is continuous and compactly supported, it follows from dominated convergence that $\lVert P_{\alpha}V_{\delta}-V_{\delta}\rVert_{L^p(\mu)}\rightarrow 0$ as $\alpha\rightarrow 1^+$. Hence, by letting $\alpha\rightarrow 1^+$ in the previous inequality we get
\begin{align*}
    \limsup_{\alpha\rightarrow 1^+}\lVert P_{\alpha}V-V\rVert_{L^p(\mu)}&\le (C+1)\delta.
\end{align*}
As $\delta>0$ was arbitrary, 2 follows. 
\end{proof}
\end{lem}
\printbibliography
\Addresses
\end{document}